\documentclass[10pt]{article}
\usepackage[a4paper, left=2.5cm, right=2.5cm, top=2.5cm, bottom=2.5cm]{geometry}

\usepackage{amsmath, amssymb, amsthm}
\usepackage{graphicx}
\usepackage{algorithm, algpseudocode, bm, verbatim}
\usepackage{xcolor}
\usepackage{dsfont}
\usepackage{hyperref}
\usepackage[authoryear]{natbib} 

\newcommand{\cl}{\mathcal}
\newcommand{\ds}{\mathds}
\newcommand{\set}[1]{\ensuremath{\mathcal{#1}}}
\newcommand{\nt}[1]{\textcolor{black}{#1}}

\theoremstyle{plain}
\newtheorem{thm}{Theorem}
\newtheorem*{thm*}{Theorem}
\newtheorem{prp}{Proposition}[section]
\newtheorem{lmm}[prp]{Lemma}

\algnewcommand{\Inputs}[1]{%
	\State \textbf{Inputs:}
	\Statex \hspace*{\algorithmicindent}\parbox[t]{.8\linewidth}{\raggedright #1}
}
\algnewcommand{\Initialize}[1]{%
	\State \textbf{Initialize:}
	\Statex \hspace*{\algorithmicindent}\parbox[t]{.8\linewidth}{\raggedright #1}
}

\title{Spectrum Estimation through Kirchhoff Random Forests}

\author{%
	S. Barthelmé\thanks{CNRS, Univ. Grenoble Alpes, Grenoble-INP, GIPSA-lab, Grenoble, France; \texttt{simon.barthelme@gipsa-lab.fr}} 
	\and
	F. Castell\thanks{Aix-Marseille Univ., CNRS, Centrale Marseille, I2M, Marseille, France; \texttt{fabienne.castell@univ-amu.fr}}
	\and
	A. Gaudillière\thanks{Aix-Marseille Univ., CNRS, Centrale Marseille, I2M, Marseille, France; \texttt{alexandre.gaudilliere@math.cnrs.fr}}
	\and
	C. Mélot\thanks{Aix-Marseille Univ., CNRS, Centrale Marseille, I2M, Marseille, France; \texttt{clothilde.melot@univ-amu.fr}}
	\and
	M. Quattropani\thanks{Sapienza Università di Roma, Italy; \texttt{matteo.quattropani@uniroma1.it}}
	\and
	N. Tremblay\thanks{CNRS, Univ. Grenoble Alpes, Grenoble-INP, GIPSA-lab, Grenoble, France and UiT the Arctic University of Norway, Department of Mathematics \& Statistics, Norway; \texttt{nicolas.tremblay@cnrs.fr}}%
}

\date{} 

\begin{document}
	\maketitle
	
	\begin{abstract}
		Given a non-oriented edge-weighted graph,
		we show how to make some estimation of the associated Laplacian eigenvalues
		through Monte Carlo evaluation of spectral quantities
		computed along Kirchhoff random rooted spanning forest trajectories.
		The sampling cost of this estimation
		is only linear in the node number,
		up to a logarithmic factor.
		By associating a double cover of such a graph with any symmetric real matrix,
		we can then perform spectral estimation in the same way for the latter.
	\end{abstract}
	

\section{Introduction}
Exact eigendecomposition of very large matrices is infeasible, due to its high
computational cost, which scales as $O(n^3)$, for a $n \times n$
matrix. In many applications, matrices are indeed very large, but exact spectral
information is not required. For instance, in applications to physics
\cite{weisse2006kernel}, it can be enough to estimate the distribution
function of eigenvalues rather than each exact eigenvalue. This problem of
estimating the \emph{spectral cumulative distribution} is the topic of this paper. 

We introduce a method for spectral estimation of large matrices, which consists
in exploiting a stochastic process to gain information on the spectrum of a
\emph{graph Laplacian}. Graph Laplacians are a subset of symmetric matrices, but
we shall see that the method can be extended to general symmetric matrices via an
embedding. The stochastic process we use is a \emph{Kirchhoff forest}. Kirchhoff
forests are a type of random spanning forests that have been used successfully
in various applications which involve linear-algebraic computations in large
graphs \cite{avena2020intertwining,jaquard2023smoothing,pilavci2021graph}.

Briefly, our approach follows the same pattern used in classical methods: we
estimate (generalised) moments of the spectral density, and combine these
estimates in order to reconstruct the cumulative spectral density.

Classical methods like the kernel polynomial method \cite{weisse2006kernel} or the stochastic Lanczos quadrature~\cite{chen_analysis_2021} 
estimate standard (polynomial) moments of the spectral distribution. Here, we
construct a set of observables which can be used to estimate the expectation of
certain \emph{rational} functions under the spectral density. We show how many
such expectations can be estimated efficiently by coupling Kirchhoff forests at
different intensities. In addition, we show that our estimators have \emph{relative}
error in $O\left(\frac{1}{\sqrt{n}}\right)$, and cost linear in $n$. 

Since we do not estimate classical but rational moments, we require a bespoke
reconstruction method. This method relies on (pointwise) transformation to a classical
truncated moment problem. We explain how to get bounds on the cumulative
spectral density, and how to perform a maximum-entropy reconstruction.

Our method yields promising results, as shown in the experiment section \ref{madeleine}. A
comparison to state-of-the-art methods can be found in a related conference
paper \cite{BCGQT}. Here, the focus is on introducing and studying the relevant
observables, on implementing the coupling process introduced in \cite{AG}, and
on formulating a reconstruction process. \nt{The end result is a spectral estimation
method with complexity linear in $n$, unlike almost all classical methods that have a 
complexity in the number of non-zero entries of the matrix (as they are based \textit{a minima} on matrix-vector multiplications). There are some recent methods, such as~\cite{braverman2022sublinear}, that propose density estimation in time linear in $n$. These methods are however based on approximating the matrix-vector multiplication by subsampling uniformly at random $\mathcal{O}(n)$ entries of the matrix, which is very different in spirit from the highly correlated sampling of random forests. A current limitation of our method is that reconstructing the cumulative density function from rational, rather than polynomial, moment estimates makes it difficult to obtain precise theoretical control over the density’s estimation error --unlike what is available for other approaches, for example in terms of the Wasserstein-1 distance in~\cite{chen_analysis_2021} or~\cite{braverman2022sublinear}.}

\subsection{Results}

To give a high-level overview of our results, we need to introduce some notation.
We focus on a specific set of symmetric matrices called \emph{graph Laplacians}
(the general case of symmetric matrices is in appendix \ref{olga}).
Consider a non-oriented edge-weighted graph without self-loops,
with vertex set $\cl X$ of finite size $n$ and with non negative weights
$$
w(x, y) = w(y, x) \geq 0,
\qquad x, y \in \cl X,
$$
using the convention that an edge connects $x$ and $y$ if and only if $w(x, y) >
0$. With such an edge we associate two distinct oriented {\it arrows\/} $(x, y)$
and $(y, x)$. We say that $(x, y)$ {\it leaves\/} $x$ and {\it reaches\/} (or
{\it points to\/}) $y$. We define the total weight of a node as
\begin{equation}
	\label{eq:total-weight}
	w(x) = \sum_{y \neq x} w(x,y)
\end{equation}
which corresponds to the degree of $x$ in the unweighted case. We note $m$ the
total number of edges.

The graph Laplacian is a symmetric matrix $L$ with entries
\begin{equation}
	\label{eq:laplacian-def}
	L_{x,y} =
	\begin{cases}
		w(x,y)\ & x \neq y \\
		- w(x) \ & x = y
	\end{cases}
\end{equation}
Graph Laplacians are diagonally-dominant matrices, and are known to have
favourable computational properties with respect to (approximate) inversion
\cite{spielman2004nearly}. Note that in our definition $L$ is
negative-definite, for consistency with the definition of the classical Laplace
operator, but many sources in spectral graph theory use the opposite convention. 

The \emph{spectral measure} of $L$ is defined as:
\begin{equation}
	\label{eq:spectral-measure}
	\sigma = {1 \over n} \sum_{j < n} \delta_{\lambda_j}.
\end{equation}
where $0 = \lambda_0 \leq \lambda_1 \leq   \dots \leq  \lambda_{n-1}$
are the eigenvalues of $-L$, counted with multiplicity.
We will seek to estimate the cumulative spectral distribution
\begin{equation}
	\label{eq:cumulative-distribution}
	F(q) = \sigma([0,q]) = \frac{1}{n}  \left| \bigl\{j < n : \lambda_j  \leq q \bigr\} \right|
\end{equation}
for different values of $q$. Defining $\alpha$ as
\[ \alpha = \max_{x} w(x) \]
a straightforward application of the Gershgorin circle theorem gives the
well-known bound $\lambda_{n-1} \leq  2\alpha$. The spectral measure has support
on $[0,2\alpha]$.

We will do so by estimating certain expectations
under $\sigma$, specifically:
\begin{equation}
	\label{eq:expectations}
	m_k(q) = \frac{1}{n} \sum_{j < n} \left( \frac{q}{q+\lambda_j} \right)^k,
	\qquad q > 0,
\end{equation}
which up to trivial rescalings correspond to the Stietjes transform of the
spectral measure (for $k=1$) and its $(k - 1)$-th derivatives, evaluated on the
positive real axis.

Our main result is that expectations of this form can be estimated efficiently,
with a control on the relative error. Since our estimators involve running random walks,
some preprocessing may be required in order to ensure that a random neighbour can be
chosen with probability $\frac{w(x,y)}{w(x)} $ in time $O(1)$. In unweighted
graphs this is trivial. In weighted graphs, if the weights are not too variable
rejection sampling can be used, otherwise we may need to perform the
preprocessing necessary to run Walker's alias algorithm at each node
\cite{walker1977efficient}, which costs $O(m)$ in total. This is no more than
the cost of reading the matrix.

Combining Theorems \ref{clo} and
\ref{ben}, we have:
\begin{thm}[Informal] Given a preprocessing step of cost at most $O(m)$, we can
	obtain an unbiased estimator of
	$m_k(q)$ for $k \leq l$ and any $r$ values of $q \in [q_{\min},q_{\max}]$ at
	cost $O\left( n l \left( \frac{\nt{\beta}}{q_{\min}}
	+ \log \frac{\alpha}{q_{\min}} + r  \right) \right)$
	where \nt{$\beta = \frac{1}{n} \sum_{x} \emph{\text{max}}_{y}\, w(y,x)$}.
	This estimator has pointwise \emph{relative} error in $O\left(\frac{1}{\sqrt{n}}\right)$.
\end{thm}
In large graphs we obtain in practice very reliable estimates, which can be used for
effective reconstruction, as we show. 

\subsection{Outline of the paper}
The remainder of the paper is organized as follows. In Section~\ref{MC.sec}
we describe the main object underlying our analysis, i.e., \emph{Kirchhoff
	forests}, and provide some insight on the way they are related to the spectral
distribution of the graph they live in. In Section~\ref{sampling.sec} we recall
the well-known \emph{Wilson's algorithm} and present some modification of the
classical algorithm which will turn out to be useful in sampling efficiently a
whole continuum of forests. Section \ref{spectralcdf.sec} is devoted to a
description of the techniques we develop to extract information about the
spectral distribution from the statistics collected by the randomized algorithm.
Section~\ref{bertrand} gathers the proofs of Theorems \ref{ben} and \ref{clo}. 
In Section~\ref{madeleine} we present the outcome of our algorithm
on some benchmark examples and comment on related numerical issues.
Appendix~\ref{olga} explains
how coupled Kirchhoff forest replicas can also be used
to perform spectral estimation with any symmetric real matrix.
We eventually give a pseudocode version
of our coupled forest algorithm in Appendix~\ref{sacha},
before collecting in Appendix~\ref{daniele} some basics, sometimes in a different form,
of the fundamental work of Krein and Nudel'man~\cite{KN} 
after Markov's work on the moment problem. These basics are needed to get bounds 
on the spectral cumulative  distribution function in Section \ref{spectralcdf.sec}.

\section{Monte-Carlo estimates using Kirchhoff forests}
\label{MC.sec}

In this section, we define the Kirchhoff forest,
and describe the observables we use to obtain Monte-Carlo estimates
of the Stieltjes transform of the spectral measure $\sigma$,
and of its derivatives.
The algorithms we use to sample the Kirchhoff forest,
and the sampling computational costs will be discussed in section \ref{sampling.sec}.

\subsection{Kirchhoff forest}

A {\bf rooted spanning forest\/} (r.s.f.)  is a collection $\phi$
of arrows such that:
\begin{itemize}
	\item[i)] for each node in $\cl X$, there is one arrow at most that leaves it;
	\item[ii)] one cannot close any loop by following arrows
	from the node they leave to the node they reach.
\end{itemize}
When no arrow in $\phi$ leaves a given node $x$ in $\cl X$,
we say that $x$ is a {\bf root\/} of the forest $\phi$.
The set of roots of $\phi$ is denoted by $\rho(\phi)$.

For any $q > 0$, we call {\bf Kirchhoff forest with rate $q$}
a random rooted spanning forest $\Phi_q$
such that for all rooted spanning forest $\phi$, 
$$
\ds P\bigl(\Phi_q = \phi\bigr)
= {q^{|\rho(\phi)|} \prod_{(x, y) \in \phi} w(x, y) \over Z(q)} 
$$
where $|\rho(\phi)|$ is the cardinality of $\rho(\phi)$
and $Z(q)$ is the normalizing constant
$$
Z(q) = \sum_{\phi \, {\rm r.s.f}}  q^{|\rho(\phi)|} \prod_{(x, y) \in \phi} w(x, y).
$$
When $q = 1$ and in the case of an unweighted graph
---i.e. when the $w(x, y)$ can only take the values 0 and~1---
Kirchhoff forest is the uniformly distributed rooted spanning forest.
For smaller or larger values of $q$ and non-uniform weights $w(x, y)$
the random forest $\Phi_q$ is sometimes still referred
as the ``uniform'' rooted spanning forest.
Because of its connection with the Gaussian free field with positive mass
it is also referred in~\cite{AEM} as the ``massive spanning forest''.
We prefer to give it another name
and call it Kirchhoff forest in view 
of his 1847's theorem that identifies the combinatorial partition function $Z$
with the characteristic polynomial of the graph Laplacian $L$.

\bigskip\par\noindent
\bf Kirchhoff's theorem~\cite{Ki}.
\it For all $q > 0$, it holds\footnote{See~\cite{AG} for a probabilistic proof
	that follows~\cite{Ma} with the notation of the present paper.}
$$
Z(q) = \det\bigl(q{\rm Id} - L\bigr).
$$ 
\rm
\smallskip\par

Kirchhoff's theorem reads 
\begin{equation}\label{logZ.eq}
	\ln Z(q) = \sum_{j < n} \ln (q + \lambda_j),
	\qquad q > 0.
\end{equation}
Computing the first derivative in $q$ of both sides of \eqref{logZ.eq} leads to 
$$
\ds E\bigl[|\rho(\Phi_q)|\bigr]
= \sum_{j < n} {q \over q + \lambda_j}.
$$
And computing the second derivative in $q$ of both sides of \eqref{logZ.eq}, one obtains
\begin{equation}\label{edith}
	{\rm Var}\bigl(|\rho(\Phi_q)|)
	\leq \ds E\bigl[|\rho(\Phi_q)|\bigr].
\end{equation}
For $q_0 > 0$ given, we describe in  section \ref{sampling.sec}
a way to sample the coupled forest process
$$
q \ge q_0 \mapsto \Phi_q
$$
whithin a cost which is almost linear in the number of nodes. Therefore, we have
access to Monte Carlo estimations of the function
\begin{equation}\label{arthur}
	q \ge q_0
	\mapsto  \ds E\bigl[|\rho(\Phi_q)|]
	= \sum_{j < n} {q \over q + \lambda_j}. 
\end{equation}
Up to sign conventions and an extra factor $qn$,
\eqref{arthur} is the Stieltjes transform
\begin{equation}\label{gregory}
	\cl S : q \geq q_0
	\mapsto \cl S(q)
	= {1 \over n} \sum_{j < n} {1 \over q + \lambda_j}
	= {1 \over qn} \ds E\bigl[|\rho(\Phi_q)|\bigr]
\end{equation}
of the Laplacian spectral measure $ \sigma$.
Moreover, this estimate is quite sharp since by \eqref{edith},
its variance is bounded by its first, rather than second, moment.  

However, inverting the Stieltjes transform over the positive real domain
is extremely ill-conditioned. Therefore, 
small estimation errors on the Stieltjes transform
will produce huge errors after inversion 
for the eigenvalues we want to estimate.

\subsection{Replicas}

A first tentative strategy to overcome the aforementioned problem
would then consist in relating the {\it complex\/} Stieltjes transform 
to some Kirchhoff forest observables.
It holds for example, for all complex number $z = q e^{i\theta}$
with $q > 0$, $\theta \in [0, 2\pi[$
and $z \neq -\lambda_0$, $-\lambda_1$,~\dots
$$
{1 \over n} \sum_{j < n} {1 \over z + \lambda_j}
= {\ds E\Bigl[|\rho(\Phi_q)|e^{i\theta|\rho(\Phi_q)|}\Bigr]
	\over n z \ds E\bigl[e^{i\theta|\rho(\Phi_q)|}\bigr]}.
$$
But the relative precision needed by this formula
for these expectations of signed or complex variables
is completely out of reach for direct Monte Carlo estimations.
We will rather define new observables closely related
to the first derivatives of the Stieltjes transform in the positive domain,
and we will later deal with the conditioning issue
using a statistical point of view.

Consider $l$ replicas $\Phi_{k, q}$, $k= 0, \cdots, l-1$,
i.e. $l$ independent copies,
of the Kirchhoff forest $\Phi_q$.
The reader should think of very modest $l$
(we will take $l = 4$ in our applications).
For any $x$ in $\cl X$ and any rooted spanning forest $\phi$,
$x$ is covered by a single tree in $\phi$.
By writing $\rho_x(\phi)$ 
for the root of this tree
and setting $R^0(x) = x$,
we inductively define
$$
R^{k + 1}(x) = \rho_{R^k(x)}(\Phi_{k, q}),
\qquad k < l.
$$
We set
\begin{equation}\label{felix}
	\xi_q^k
	= \bigl\{x \in \cl X : R^k(x) = x\bigr\} 
	= \bigl\{
	x \in R^k(\cl X)
	: R^k(x) = x
	\bigr\},
	\qquad k = 0, \dots, l,
\end{equation}
so that $\xi_q^0 = \cl X$,
$\xi_q^1 = \rho(\Phi_{0, q})$
and 
$\xi_q^{k + 1}\subset \rho(\Phi_{k, q})$
for all $k < l$. 
We will prove in Section~\ref{bertrand}
\begin{thm}\label{ben}
	For all $k \leq l$, it holds
	$$
	\ds E\bigl[|\xi_q^k|\bigr]
	= \sum_{j < n} \left(q \over q + \lambda_j\right)^k . 
	$$
	In addition, 
	$$
	\ds P\bigl(x, y \in \xi_q^k\bigr)
	\leq \ds P\bigl(x \in \xi_q^k\bigr) \ds P\bigl(y \in \xi_q^k\bigr),
	\qquad x \neq y \in \cl X,
	$$
	so that 
	\begin{equation} \label{edith.bis}
		{\rm Var}(|\xi_q^k|) \leq \ds E[|\xi_q^k|].
	\end{equation}
\end{thm}

Using  $s$ replicated forests,
i.e. $ls$ independent copies of $\Phi_q$, one can therefore estimate
$\ds E\bigl[|\xi_q^k|\bigr]$ by 
\begin{equation}
	\label{MC.eq}
	\hat{\xi}_q^{k} = \frac 1 s  \sum_{j=1}^s \left| \xi_q^{k,(j)}  \right|, 
\end{equation}
where $(\xi_q^{k,(j)} , j= 1, \cdots, s)$ are $s$ i.i.d. copies of $\xi_q^k$, 
each one obtained from $l$ independent copies 
of the Kirchhoff forest. This estimate has relative error 
$ ( \hat{\xi}_q^{k} - \ds E\bigl[|\xi_q^k|\bigr] ) / \ds E\bigl[|\xi_q^k|\bigr]$ of order 
$$
\frac{\sqrt{{\rm Var}(|\xi_q^k|)}}{\sqrt{s} \ds E\bigl[|\xi_q^k|\bigr]}
\le \frac{1}{\sqrt{s \ds E\bigl[|\xi_q^k|\bigr]}}
\le \frac{C}{\sqrt{sn}},
$$
as long as $\ds E\bigl[|\xi_q^k|\bigr]$
is of order $n$. This will be the case as soon as $q_0$ is chosen to be of the order of 
$\bar{\lambda} = \frac 1 n \sum_{j < n} \lambda_j$ since
$$
n \ge  \ds E\bigl[|\xi_q^k|\bigr]
\ge n  \left(\frac{q}{q + \bar{\lambda}} \right)^k  
\ge n \left(\frac{q_0}{q_0 + \bar{\lambda}} \right)^k.
$$
Up to simple and explicit multiplicative factors,
$(\ds E\bigl[|\xi_q^k|\bigr], k <  l)$ are the Stieltjes transform 
of the spectral measure $\sigma$ and its first derivatives:
$$
\cl S^{(k)}(q) = \frac{(-1)^k k!}{ n} \sum_{j < n} \frac{1}{(q+\lambda_j)^{k+1}} 
= \frac{(-1)^k k!}{ n q^k} \ds E\bigl[|\xi_q^{k+1}|\bigr]. 
$$

We close this section with a last comment on the bounds~\eqref{edith} and~\eqref{edith.bis}. 
Kirchhoff's theorem actually characterizes the law or $|\rho(\Phi_q)|$,
as the law of a sum of independent 0-1 Bernoulli random variables
with mean $q / (q + \lambda_j)$, $j < n$.
This identity in law can in turn be seen as a consequence
of the determinantality\footnote{%
	This means that it exists a matrix
	$K_q = \bigl(K_q(x, y)\bigr)_{x, y \in \cl X}$
	such that for all $A \subset \cl X$,
	$$
	\ds P\bigl(A \subset \rho(\Phi_q)\bigr)
	= \det\bigl(K_q(x, y)\bigr)_{x, y \in A}.
	$$
	Here we have $K_q = q(q{\rm Id} - L)^{-1}$
	(see~\cite{AG}).
}
of the point process $\rho(\Phi_q)$,
which derives from that of the arrow process $\Phi_q$
proved in~\cite{BP}.
Now, the determinantality of the point process $\rho(\Phi_q)$
leads directly to its negative correlations:
for all $x \neq y$ in $\cl X$, 
$$
\ds P\bigl(x, y \in \rho(\Phi_q)\bigr) 
\leq \ds P\bigl(x \in \rho(\Phi_q)\bigr)
\ds P\bigl(y \in \rho(\Phi_q)\bigr).
$$
And such negative correlations 
{\it by themselves\/} imply the bound~\eqref{edith}.
The  point processes $\xi_q^k$ are  {\it not\/}  determinantal
but {\it have\/} negative correlations, which is sufficient to get~\eqref{edith.bis} .

\section{Sampling Kirchhoff forests}
\label{sampling.sec}

In this section, we address several sampling issues:
\begin{itemize}
	\item We describe Wilson's algorithm whose output is a r.s.f. $\Phi_q$ for some given $q>0$.
	\item We explain next how to couple all  the $(\Phi_q , q>0)$
	using the stack representation of Diaconis \& Fullton.
	\item  We describe the process $t>0 \mapsto \Phi_{1/t}$
	as a non homogeneous Markov jump process
	on the space of rooted spanning forest.
	This description involves an auxiliary Markov process
	on the r.s.f.  with active and frozen trees.
	At each jump of $\Phi_{1/t}$,
	the new state will be some final configuration 
	of the auxiliary Markov process.
	\item We describe the actual algorithm we use to sample trajectories
	$t>0 \mapsto \Phi_{1/t}$.
	\item We bound the computational cost of this algorithm.
\end{itemize}
Except for the final point, the content of this section largely rephrases the material in \cite{Wi,AG}, with the aim of presenting the coupled forest process in a more pedagogical and accessible manner.

\subsection{Wilson's algorithm} 
\label{Wilson.sec}

For any given $q > 0$,
Wilson's algorithm builds a random rooted spanning forest 
in the following way:
\begin{itemize}
	\item[i)] starting from any node $x$, 
	it runs the continuous time random walk
	with generator $L$ up to an independent 
	random exponential time $T_q$ of rate $q$,
	and erases all the possibly created loops
	as soon as they appear.
	It produces in this way
	a first path $\Gamma^x_q$, which is oriented and self-avoiding,
	joining $x$ to the node where the random walk
	was stopped.
	Call $A$ the collection of arrows along this path
	and $B$ the set of the $|A| + 1$ nodes
	along $\Gamma^x_q$.
	\item[ii)] Starting from another node $x$ outside $B$, if any,
	the algorithm runs another loop-erased random walk
	up to the minimum of another independent exponential time of rate $q$
	and the hitting time of $B$.
	The algorithm produces in this way a new self-avoiding path $\Gamma^x_{q, B}$.
	It adds to $A$ the arrows along that path,
	and adds to $B$ the new nodes along it.
	\item[iii)] If $B = \cl X$
	the algorithm returns the rooted spanning forest $\Phi_q = A$; 
	if not it repeats steps ii)--iii).
\end{itemize}

\medskip\par\noindent
\bf Wilson's theorem~\cite{Wi}.
\it The output $\Phi_q$ of the algorithm is a Kirchhoff forest 
with rate $q$.
\rm
\medskip\par

It is useful to attribute a status to each node $x$ during the algorithm.
A node is said to be 
\begin{itemize}
	\item {\bf frozen} if it belongs to an edge of $A$;
	\item {\bf active} otherwise.
\end{itemize}
The current node of the loop-erased random walk is by definition active.
The algorithm is initialized by setting all the nodes to active
and by choosing a first current node.

\subsection{Stack representation and coupled forests.}
\label{fiorella}

We stress the crucial fact that Wilson's theorem implies
that the way the starting nodes are chosen in steps i) and ii)
is irrelevant to the law of $\Phi_q$ as constructed by the algorithm.
This feature of Wilson's algorithm is referred to as its {\bf abelianity\/}.

This is made clear by Wilson's proof that relies
on the stack representation of random walks introduced
in~\cite{DF}.
The latter associates with each node $x$ an infinite stack of arrows $(x, y)$ or stops  
that appear
with probability $w(x, y) / (q + w(x))$ and $q / (q + w(x))$, respectively,
at each level independently.
We can then build the continuous time loop-erased random walk
that is stopped at rate $q$
by reading and removing at rate $q + w(x)$
the first stop or arrow 
from the stack associated with the current node $x$
before processing this stop or arrow,
i.e.,  stopping after reading a stop 
or following it after reading an arrow. 

\begin{figure}[tbp]
	\centerline{\includegraphics[scale=0.4]{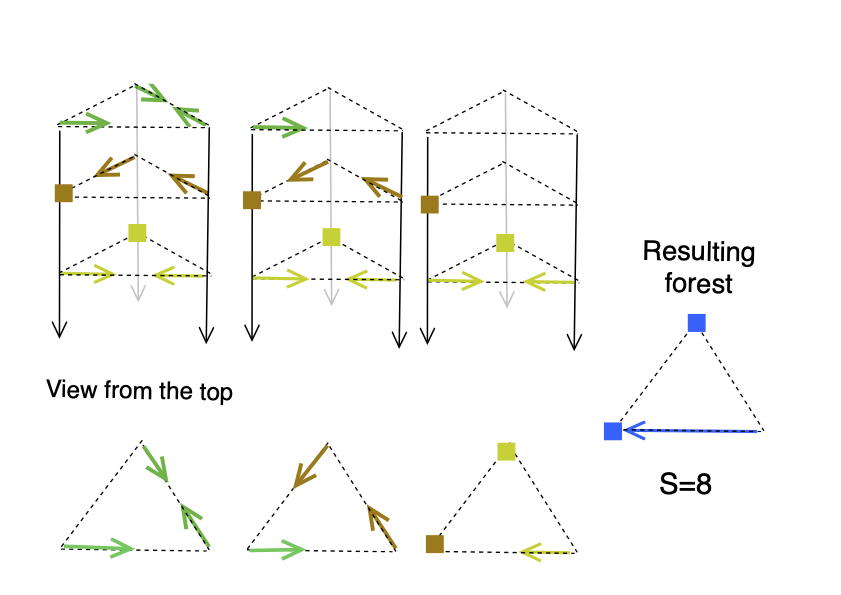}}
	\caption{Cycle erasure in the stacks viewed from the top, in a triangle.
		The squares in the stacks correspond to stops.
		$S$ is the number of arrows or stops that have to be read from the stacks
		to obtain a rooted forest.}
	\label{stack.fig}
\end{figure}

\begin{figure}[tbp]
	\hbox to \hsize{%
		\includegraphics[width = 2.77 in]{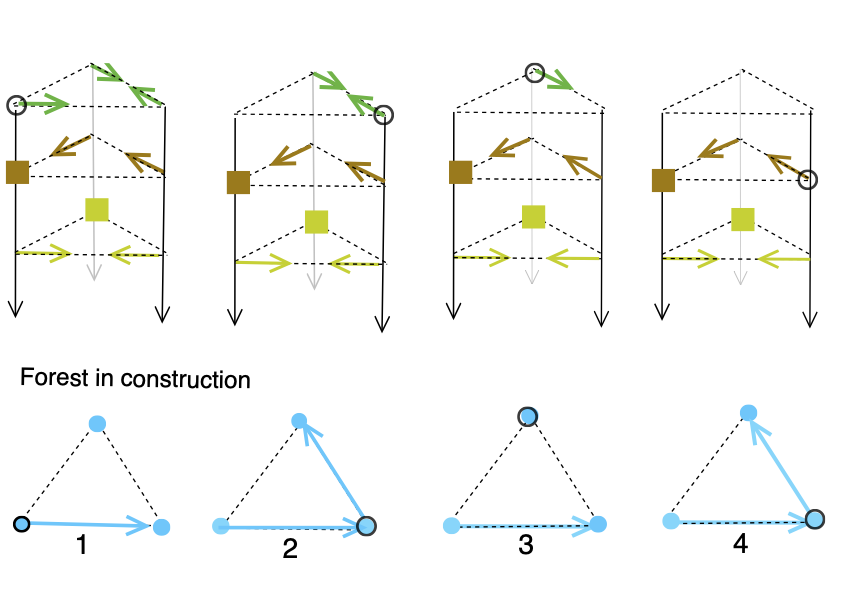}%
		\hfill%
		\includegraphics[width = 2.77 in]{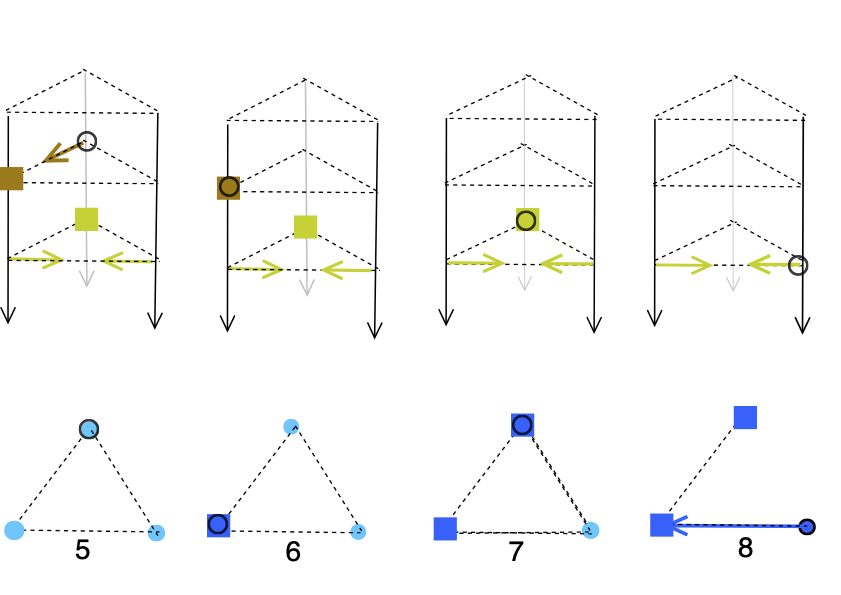}%
	}
	\caption{Constructing the forest by reading the stacks according to Wilson's order.
		The light blue nodes are active, while the dark blue ones are frozen.
		The current node is marked by a black circle.
		$S$ is now the number of steps, i.e. reading a stop or an arrow,
		before the rooted forest is entirely constructed.}
	\label{Wilson.fig}
\end{figure}

\begin{figure}[tbp]
	\hbox to \hsize{%
		\includegraphics[width = 2.77 in]{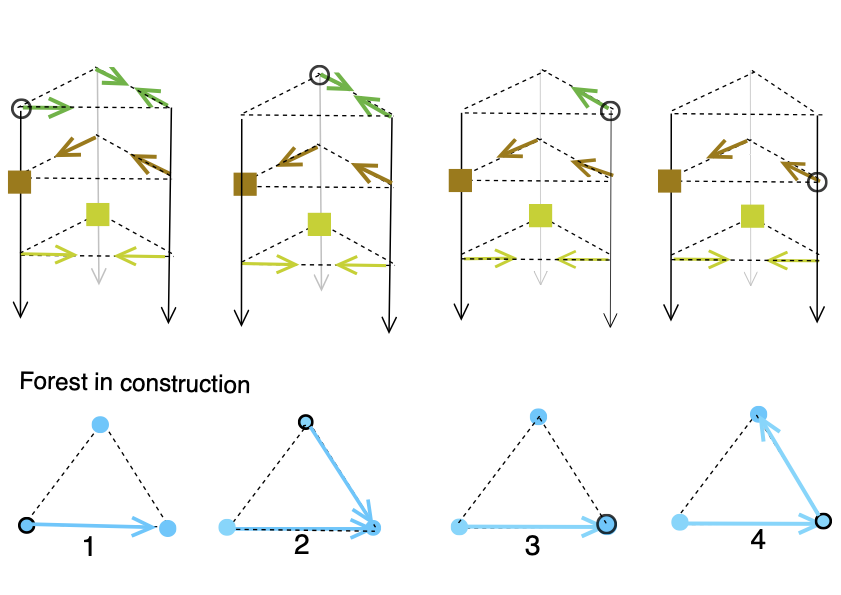}%
		\hfill%
		\includegraphics[width = 2.77 in]{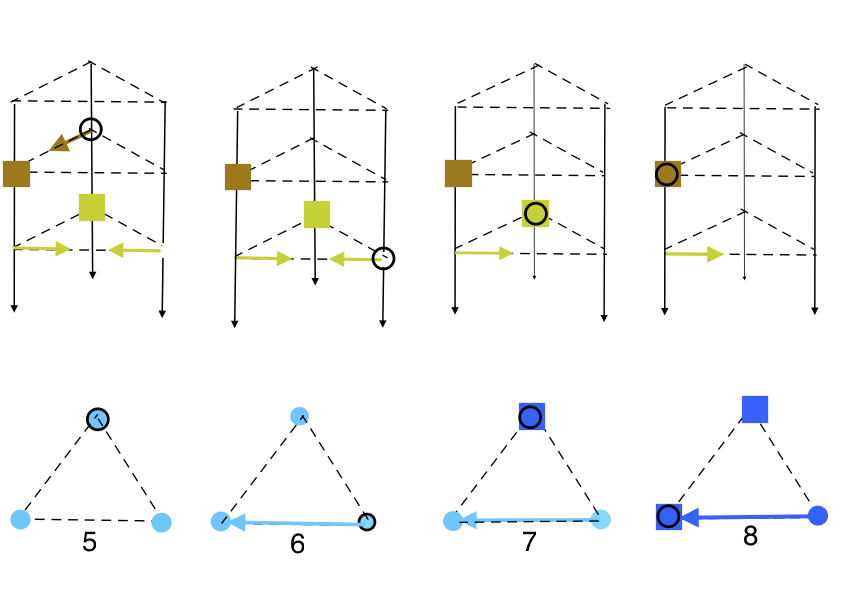}%
	}
	\caption{Constructing the forest by reading the stacks in random order.
		Unlike the Wilson's order,
		the current node is chosen randomly among the active roots.}
	\label{Random.fig}
\end{figure}

With such a construction,
whatever the chosen starting nodes,
the same loops will be erased by Wilson's algorithm,
the same arrows will eventually form $\Phi_q$,
and its roots will lie at each node
where a stop has been read.
Actually, as explained in \cite{Wi},
the forest obtained from the algorithm is the forest one would obtain by 
removing the cycles when looking at the stacks from the top,
as illustrated in Figure \ref{stack.fig} in the case of a triangle. 
Therefore, it is not even mandatory to follow the path of the random walk
to read the arrows in the stacks. 
We can as well choose as current node any active node which is a root of a tree.  
We will speak of {\bf Wilson's order} when the arrows in the stacks
are read according to the original loop-erased random walk paths.
Figure \ref{Wilson.fig} illustrates the forest construction
when the stacks are read according to Wilson's order, 
while in Figure \ref{Random.fig} the current node is chosen at random
among the active roots. 

As noted in~\cite{AG}, 
such a stack representation
allows for a coupling of all $\Phi_q$, $q >0$, together.
Instead of sampling stops {\it or\/} arrows $(x, y)$
at all levels of the stack associated with a node $x$,
we can independently sample
marks $U$ in $[0, 1]$ {\it and\/} arrows $(x, y)$,
the former being uniformly distributed on $[0, 1]$
and the latter appearing now with probability $w(x, y) / w(x)$. This alea 
does not depend on $q$. 
To build $\Phi_q$ it suffices to interpret (or replace) 
any arrow by a stop each time its associated mark $U$
is smaller that $q / (q + w(x))$.

It is  explained in~\cite{AG} and further develop in sections \ref{aurelie} and \ref{farah}
that this coupling leads to piecewise-constant random forest trajectories
$$
t \geq 0 \mapsto \Phi_{1 / t}
$$
that start from the deterministic empty forest $\Phi_\infty$
made of $n$ trees reduced to a simple root. This random forests process has 
the property that all arrows that are read to compute 
a Kirchhoff forest $\Phi_q$ at time $t = 1 / q$, 
will also be read to compute the Kirchhoff forest $\Phi_{q'}$
with $q' < q$, at a later time $t' = 1 / q'$.
The jumps of the random forest trajectories
can be described through an auxiliary random process
that we introduce in the next section.
The coupled forest process will be described
as a non-homogeneous Markov process in section \ref{farah}
and the sampling algorithm we actually implement is given in section \ref{simone}. 

\subsection{An auxiliary Markov process: frozen and active trees.}
\label{aurelie}
For fixed $q > 0$, we can see the Kirchhoff forest $\Phi_q$
produced by Wilson's algorithm
as the final configuration of a continuous time Markov process $(\Psi^q_s)_{s \geq 0}$
with a family of absorbing states or frozen configurations.
As we will discuss later, this is {\it not\/} convenient 
from a practical, algorithmic, point of view.
But this helps to describe the coupled forest process.

Consider a rooted spanning forest
and declare each tree to be either {\it active\/} or {\it frozen\/}.
Such forests $\phi$ with active and frozen trees
will be the configurations of our Markov process.
Each root is also declared to be active or frozen,
in accordance with the status of the tree it belongs to.
We will denote by $(\phi, A)$ the configuration
with active root set $A \subset \rho(\phi)$
in the forest $\phi$,
and by $\Psi^{q; \phi, A}$ the process $\Psi^q$
started at $(\phi, A)$.

Here are the evolution rules of $\Psi^q$:
\begin{itemize}
	\item frozen roots will remain such;
	\item each active root freezes at rate $q$,
	and so does the tree it belongs to;
	\item each active root in $x \in \cl X$
	also adds at rate $w(x, y)$ the arrow $(x, y)$
	to the current forest $\phi$
	with two possible outcomes:
	\begin{itemize}
		\item a cycle appears in $\phi \cup \{(x, y)\}$,
		i.e. if $(x, y)$ joins $x$
		to some node  $y$ in its tree. In this case,  each arrow in the cycle is erased,
		the tree is fragmented, new roots appear 
		at each site along the erased cycle, and $\phi$ is updated to this smaller forest
		(with less edges).
		\item otherwise $\phi' = \phi \cup \{(x, y)\}$
		is a rooted spanning forest
		and $\phi$ is updated to $\phi'$,
		so that the root $x$ disappears,
		and two trees coalesce.
	\end{itemize}
	In the fragmentation case,
	new trees are declared active trees.
	In the coalescence case,
	the status of the new tree
	is inherited from that of its root.
\end{itemize}
$\Psi^q$ can then evolve in four different manners only:
\begin{itemize}
	\item by freezing a tree;
	\item by fragmenting an active tree into smaller active trees;
	\item by grafting an active tree on another one
	to get a larger active tree;
	\item by grafting an active tree on a frozen tree
	to get a larger frozen tree.
\end{itemize}

\noindent\textbf{Building $\Phi_{q}$ from $\Psi^{q; \emptyset, \cl X}$---}
By starting $\Psi^q$ from the configuration with $n$ active roots,
(so that the initial forest $\phi$ is the empty one), 
running $\Psi^q$ amounts to discovering marks and arrows
of the stack representation in a random order determined by the successive times
when the active roots add a new arrow. Indeed, the difference with Wilson's order, 
where only the current node can add an edge or freeze,
is that here all active roots can add edges and freeze.
This is illustrated in Figure \ref{Random.fig}. 
The abelianity implies that $\Phi_q$
is the final configuration of $\Psi^{q; \emptyset, \cl X}$.

\noindent\textbf{Building $\Phi_{q'}$ from $\Psi^{q'; \Phi_q, W}$
	with $q > q'$ and $W$ a random subset of $\rho(\Phi_q)$---}
By conditioning on $\Phi_q$,
the stack distribution of the non-erased arrows and marks is
described with the same product measure as the original one
with an additional bias:
the first mark at each root $x\in\rho(\Phi_q)$
is uniformly distributed on $[0, q / (q + w(x))]$ instead of $[0, 1]$.
Hence, in order to build $\Phi_{q'}$ from $\Phi_q$,
we only have to sample a set  $W$ of roots to unfreeze
by independently choosing each root of $x\in\rho(\Phi_q)$
with a suitable probability $p_x$,
before running $\Psi^{q'; \Phi_q, W}$ up to freezing.
According to our original coupling,
the probability $p'_x$ that a frozen root $x$ of $\Phi_q$
remains such for our smaller $q'$ without reading any further arrow in $x$,
is the probability that a uniform variable on $[0, q / (q + w(x))]$
is smaller than $q' / (q' + w(x))$ :
$$
p'_x = {q' / (q' + w(x)) \over q / (q + w(x))}.
$$
By running $\Psi^{q'; \Phi_q, W}$ up to freezing,
this is achieved in two ways only:
by keeping $x$ out of $W$, or by unfreezing it before discovering 
that the first {\it unbiased} mark in $x$
is below the threshold $q' / (q' + w(x))$.
The probability $p_x$ must then solve
$$
(1 - p_x) + p_x {q' \over q' + w(x)}
= p'_x,
$$
so that
\begin{equation}\label{caroline}
	p_x = 1 - {q' \over q}.
\end{equation}
Note that $p_x$ does not depend on $x$.

\subsection{Coupled forests as a non homogeneous Markov process.}
\label{farah}
After observing that $(\Phi_\infty, \emptyset)$
is the deterministic final frozen configuration of $\Psi^{\infty; \emptyset, \cl X}$, 
we can finally describe $\Phi_{1 / \cdot}$
as a non-homogeneous Markov process on the frozen forests,
started at $\Phi_\infty = \emptyset$,  and indexed by the time $t = 1/ q$. 

In the limit $q \rightarrow q'_+$ and using
$$
-{dq \over q} = {dt \over t},
$$
Equation~\eqref{caroline}
gives that each root $x$ in $\rho(\Phi_{1 / t_-})$
wakes up at rate $1 / t$, and  $\Phi_{1 / t}$ is obtained as 
the final forest of $\Psi^{1 / t}$
started in $(\Phi_{1 / t_-}, \{x\})$.
The difficulty raised by the fact that rates $1 / t$ diverge
in the neighbourhood of 0 is removed by observing 
that any wake-up of a given root $x$ at a given time $t$
is ineffective if $x$ is frozen again
before putting any arrow to some neighbour $y \neq x$. This
occurs with probability
$$
{q \over q + w(x)}
= {1 \over 1 + t w(x)}
= 1 - {t w(x) \over 1 + t w(x)}.
$$ 
It follows that the effective wake-up rate in $x$ at time $t$
is $w(x) / (1 + t w(x))$, which is bounded.

By considering only effective wake-up events, 
we obtain the following description of the frozen forest process $\Phi_{1 / \cdot}$,
which starts from the empty forest with $n$ frozen roots:
\begin{itemize}
	\item each root is unfrozen at rate $w(x) / (1 + t w(x))$;
	\item when a root $x$ is unfrozen at time $t = 1 / q$, 
	an arrow $(x, y)$ is chosen 
	with probability $w(x, y) / w(x)$
	to be added to $\Phi_{1 / t_-}$
	before being processed with grafting or pruning;
	\item in the pruning case
	---when $y$ is rooted in $x$ in $\Phi_{1 / t_-}$---
	$\Phi_{1 / t}$ is the final configuration of $\Psi^{1/t; \phi, A}$,
	with $\phi$ the pruned forest obtained by removing from $\Phi_{1 / t_-} \cup \{(x, y)\}$ 
	all the arrows that form its only one cycle,  
	and with $A$ the foot set of those arrows;
	\item in the grafting case
	---when $y$ is not rooted in $x$ in $\Phi_{1 / t_-}$---
	$\Phi_{1 / t}$ is $\phi = \Phi_{1 / t_-} \cup \{(x, y)\}$,
	final and initial configuration of $\Psi^{1/t; \phi, \emptyset}$.
	
\end{itemize}

\subsection{Wilson's order}
\label{simone}

While relevant  to describe $(\Phi_{1 / t})_{t \geq 0}$
as a non-homogeneous Markov process,
the update rule account at unfreezing times $t = 1 / q$
through the Markov process $(\Psi^q_s)_{s \geq 0}$
is not convenient for practical sampling. Indeed, 
exploring the stacks in this random order, 
with many active trees evolving {\it in parallel\/}
with a same continuous time $s$,
would {\it a priori\/} imply to read many times
the already sampled arrows to decide
between fragmentation and coalescence
at each new arrow occurrence.
By running {\it successive\/} loop-erased random walks,
Wilson's exploring order
deals in contrast with only one non-trivial active tree
at each time, actually reduced to a single branch.
In this situation, no extra read is needed
to decide between coalescence or fragmentation
(either the new arrow points to an active node and the active tree is fragmented,
or the new arrow points to a frozen node and the tree coalesces to another frozen one).
Extra reads are possibly required for cycle erasure only,
in the fragmentation case\footnote{When sampling a single forest $\Phi_q$ for a given $q$,
	the last version of the algorithm in \cite{Wi} actually avoids any cycle erasure.
	We reproduce it in Appendix \ref{sacha}.}.

\begin{figure}[tbp]
	\centerline{\includegraphics[scale=0.5]{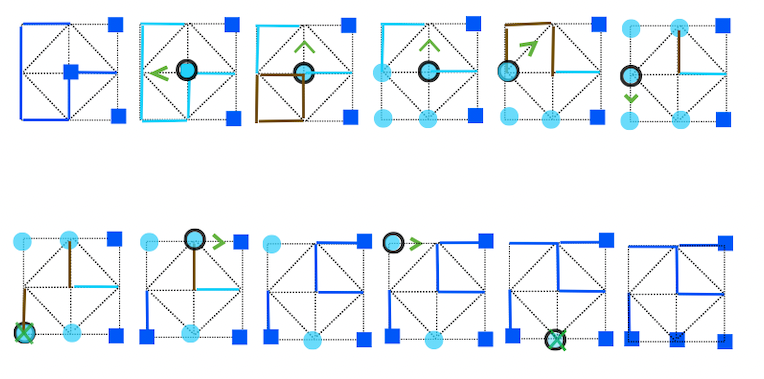}}
	\caption{Computing the extra reads in the coupled forest algorithm following Wilson's order. 
		The dark blue trees are frozen.
		The light blue ones are active.
		The root of the current tree is marked by a black circle, 
		and the next arrow or stop in the corresponding stack is indicated in green,
		a cross corresponding to a stop. 
		New sampled arrows are in brown,
		as are the arrows that we have to read again
		in order to decide between a cycle or a grafting to an active tree.
		Once in brown, the decision between ``cycle or grafting'' does not need more extra reads:
		if the current root points to a brown arrow, there is a cycle;
		if it points to a light blue one, the current tree coalesces.
		The number of reread arrows $R$ is the number of times one light blue arrow becomes brown 
		($R= 5$ in our example).}
	\label{Rereads.fig}
\end{figure}

For this reason we will stick to Wilson's order
in building $\Phi_{1 / t}$ from the configuration
$(\Phi_{1 / t_-}, \{x\})$ with one active tree only,
when some root $x$ wakes up at time $t$.
After adding and processing a new arrow $(x, y)$,
we will always freeze or add and process arrows 
from the root of the tree that covers $x$ 
in the current forest, until freezing
or grafting on a frozen tree.
We will then choose another active root $x$,
which will necessarily belong to the nodes
of the unfrozen tree at time $t$,
and proceed in the same way. 
In repeating this procedure up to complete freezing,
we may have to read again some previously sampled arrows,
used to build $\Phi_{1 / t_-}$,
while building $\Phi_{1 / t}$.
However, by proceeding however in this order,
any read again or re-sampled arrow
that might be followed
to decide between coalescence and fragmentation
would be part a cycle to be immediately erased.
As a consequence
and as detailed in the pseudo-code appendix~\ref{sacha},
we can tag each node where an arrow as been read again or re-sampled,
before re-sampling its mark and arrow each time we are led to this node again.
Each needed extra read,
which is associated with a previously sampled arrow
used to build $\Phi_{1 / t_-}$,
is performed at most once for each node initially covered
by the active tree we started with.
This is the key for the control 
in Theorem~\ref{clo} 
of the number of extra reads required,
and is illustrated in Figure \ref{Rereads.fig}.   

In Appendix~\ref{sacha} we give a pseudocode description
of this procedure to build the process
$(\Phi_q)_{q_{\max} \geq q \geq q_{\min}}$
for $q_{\max} > q_{\min} > 0$ in terms of two algorithms.
The first one follows Wilson's algorithm to sample $\Phi_{q_{\max}}$
in terms of a list of all node successors 
together with a priority queue $\mathcal{P}$ 
that orders the roots $x$ of $\Phi_{q_{\max}}$ 
according to their unfreezing time $1 / q_x$.
The second one extends this priority queue $\mathcal{P}$
to describe all unfreezing times $t = 1 / q$ between $1 / q_{\max}$ 
and $1 / q_{\min}$ together with the successor lists
that describe the associated forests $\Phi_q$.

\subsection{Sampling costs}
We turn now to the cost of the algorithm described in section \ref{simone}
to sample one trajectory $\Phi_{1/t}$ on the whole time interval $[0, 1 / q_0]$.
This cost is essentially made of two ingredients:
\begin{itemize}
	\item the sampling cost of the number $S$ of arrows to be sampled,
	in order to build this trajectory;
	\item the total number $R$ of extra reads of already sampled arrows,
	which are possibly kept in the new forest built at each unfreezing time.
\end{itemize} 
Theorem \ref{clo} gives estimates on the mean number $\ds E(S)$ and $\ds E(R)$.
It will be proved in Section~\ref{bertrand}. 

\begin{thm} \label{clo}
	Let $W$ denote the diagonal matrix whose entries are $(w(x), x \in \cl X)$.
	For all $q_0 > 0$,
	the sampled arrows number $S$
	and the extra reads number $R$
	needed to build one random forest trajectory
	$$
	t \in [0, 1 / q_0] \mapsto \Phi_{1 / t}
	$$
	satisfy 
	$$
	\ds E[S]
	= {\rm Tr}\bigl((q_0{\rm Id} - L)^{-1}(q_0 {\rm Id}  + W)\bigr)
	\leq n \left(1 + {\nt{\beta} \over q_0}\right) 
	$$
	and
	$$
	\ds E[R] 
	\leq n \ln\left(
	1 + {\alpha\over q_0}
	\right).
	$$
\end{thm}

Using Walker's alias algorithm to sample arrows at a given vertex at cost  $O(1)$, 
we get a total sampling cost for one trajectory
in $O(n (\nt{\beta}  q_0^{-1}  + \ln( \alpha q_0^{-1}))$. 
In the case when we start from the weight matrix,
this of course amounts to assuming that each non-zero entry
has been previously read once.
\nt{The total sampling cost of one coupled forest thus scales with the number of nodes. However, it grows rapidly as $q_0$ approaches zero. In our experiments, we find that setting $q_0$ to a small fraction of the mean degree\footnote{\nt{we call mean degree here the mean value of the diagonal of $W$: it is equal to the mean eigenvalue of $-L$, thus the notation $\bar\lambda$}} ($q_0=\epsilon \bar\lambda$ with $\epsilon\approx0.01$) already produces very good results.}

\section{Estimating the spectral cumulative distribution function.}
\label{spectralcdf.sec}

We turn now to  the inverse problem
of getting estimates of the cumulative distribution function of $\sigma$
$$
F : q \in [0, 2 \alpha]
\mapsto {1 \over n} \left| \bigl\{j < n : \lambda_j \leq q\bigr\}  \right|
$$
from the estimates of the Stieltjes transform
and its derivatives defined in section \ref{MC.sec}. 
We look at this problem through a moment problem.
Recall that $\sigma$ is indeed supported on $[0, 2\alpha]$ since
there exists a unique stochastic matrix $P$ satisfying
$$ L=\alpha(P-\rm{Id});$$
the spectrum of any symmetric stochastic matrix $P$
being  between $-1$ and $1$,
we get
$$
0 \leq \lambda_j \leq 2 \alpha,
\qquad j < n.
$$

\subsection{A moment problem}
\label{moment.sec}

Take a positive $q \leq 2 \alpha$ and consider the random variable 
$$
Y_q = {q \over q + \lambda_J} \in \left[{q \over q + 2 \alpha}, 1\right]
$$
with $J $ uniformly distributed on $\left\{ 0, \cdots, n-1 \right\}$. 
After rescaling, an estimate $\hat\xi_q^k$ of
$$
\ds E\bigl[|\xi_q^k|\bigr]
= \sum_{j < n} \left(q \over q + \lambda_j\right)^k,
$$
for $k \geq 1$,
is an estimate $\hat m_k$ of the $k^{\rm th}$ non-trivial moment of $Y_q$
$$
m_k
= E\bigl[Y_q^k\bigr]
= {1 \over n} \sum_{j < n} \left(q \over q + \lambda_j\right)^k.
$$
Since 
$$
F(q)
= {1 \over n} \left| \bigl\{j < n : \lambda_j  \leq q \bigr\}\right|
= P\left(Y_q \geq {1 \over 2}\right) 
$$
we are led to the problem of estimating 
the tail distribution function $G$ of a bounded-support random variable
$Y_q$ with values in $[a, b] = [q / (q + 2 \alpha), 1]$ at one precise point $y = 1 / 2$,
given its first $l$ non-trivial moments
$m_1$,~\dots, $m_l$.
This problem was solved by Markov 
who gave the sharpest possible lower and upper bounds on $G(y)$
given the constraints. We refer to~\cite{KN} 
and Appendix~\ref{daniele} in the present paper for details.

\begin{figure}[tbp] 
	\caption{\footnotesize%
		Cumulative distribution function (dashed line)
		in natural (left) and log-log (right) scales
		of the spectral measure
		for the bunny graph with 2053 nodes
		and mean degree 52.33
		together with lower (upward triangles) and upper (downward triangles) bounds
		computed from Monte Carlo estimation of $m_1$, \dots, $m_4$
		after sampling 400 replicated forest trajectories
		up to time $1 / q_0$ with $q_0 = \bar\lambda / 100$
		and with 1 (yellow), 2 (green), 3 (cyan) or 4 (blue)
		valid moment estimates.
		\label{brune}
	}
	\sbox0{\includegraphics{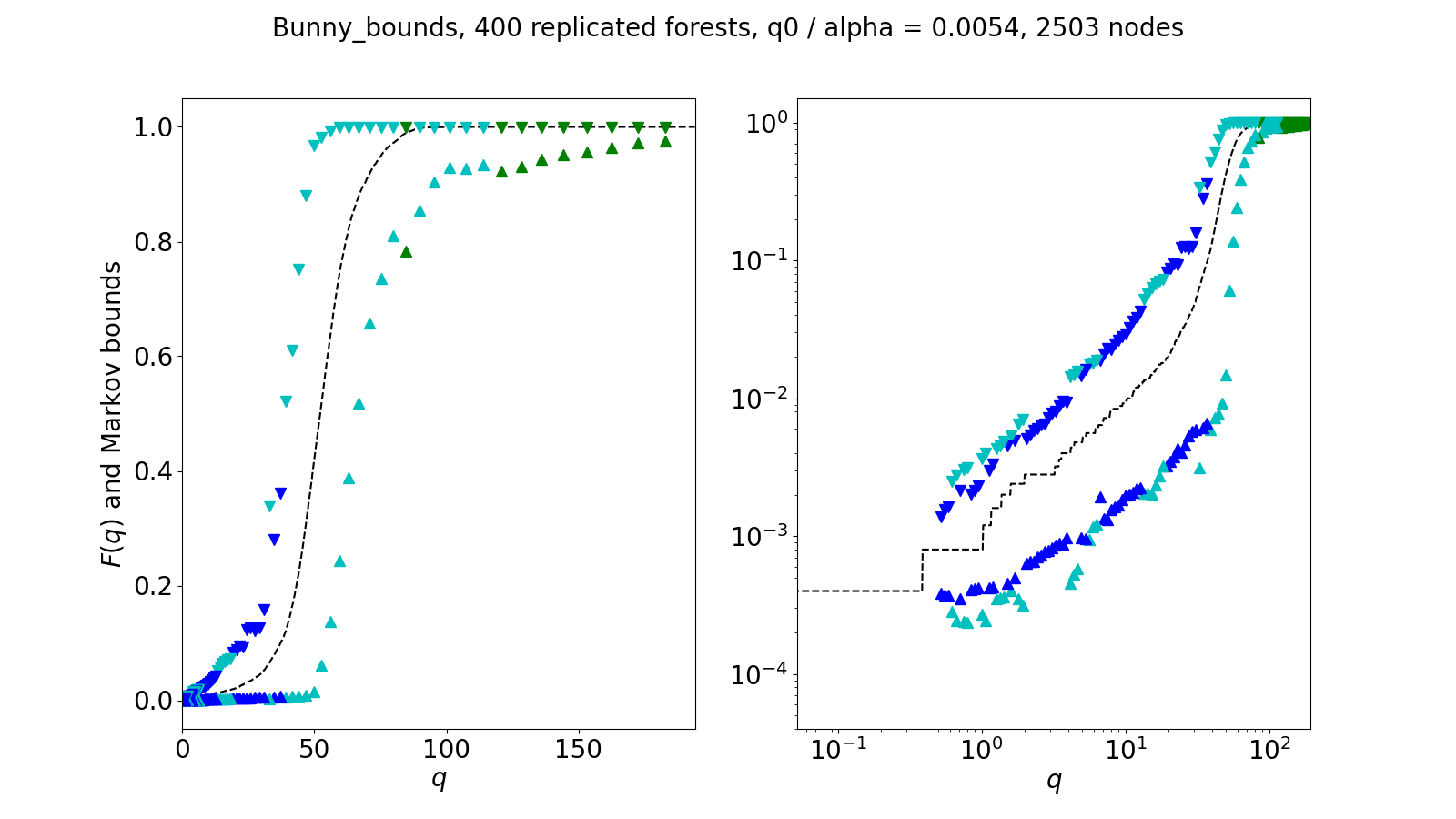}}
	\hbox to \hsize{\hfill\includegraphics[clip, trim=0 {.02\ht0} 0 {.06\ht0}, width=5.6in]{./pictures/Bunny_bounds.png}\hfill}
\end{figure}

In the particular case
of the unweighted ($w(x, y) \in \{0, 1\}$) Stanford bunny graph
---a subsampled surface of a ceramic bunny figuring, in which each point 
is connected to a few dozens of its nearest neighbours---
Figure~\ref{brune} shows Markov bounds,
computed from Monte Carlo estimation $\hat m_k$
of $m_k$ with $k \leq l = 4$,
along $s = 400$ forest trajectories with $l$ replicas
(a total of 1600 forest trajectories),
on \hbox{$q \in [q_0, 2\alpha] \mapsto F(q)$}
with $q_0 = \bar\lambda / 100$.
Colors show in each value of $q$
the number of ``valid'' moment estimates:
not all finite sequences $\hat m_1$, \dots, $\hat m_k$
are indeed admissible moment sequences
for a probability distribution on a given interval $[a, b]$
(see section \ref{daniele} for a precise definition) .
In first approximation, we say
that the number of valid moment estimates 
is the largest $k$ for which any sequence
$\hat m_1$, \dots, $\hat m_{k - 1}$, $\tilde m_k$
with $\tilde m_k$ in the $k^{\rm th}$
$95\%$ confidence interval
associated with our Monte Carlo estimates, 
is an admissible moment sequence.
We refer to Section~\ref{madeleine} for more details.

These are very bad, or quite bad, bounds
on the cumulative distribution function $F$
of the spectral measure $\sigma$.
In this example $F(50)$
is essentially estimated 
by the trivial bounds 0 and 1.
Only the median of the spectral measure
or the spectrum bulk location
can be roughly identified with Markov bounds of Figure 1:
these are around $q = 45$.
An explanation of this poor performance is that Markov bounds are reached
by very singular distributions,
supported by three atoms at most for $l = 4$,
while $Y_q$ can generically take $n$ distinct values.
We will go beyond Markov bounds
by using a standard equivalence
between canonical and microcanonical ensembles
to make a prediction on $F$,
instead of bounding it.

\subsection{Maximal entropy predictions}
\label{keenan}

Going from our Monte Carlo noisy estimates
made from $s$ samples $\xi^k_{q, i} \subset \cl X$, $i < s$,
$$
\hat m_k: q \geq q_0 \mapsto \hat m_k(q)
= {1 \over sn} \sum_{i < s} |\xi^k_{q, i}|,
\qquad k \leq l,
$$
of the moment functions
$$
q \geq q_0 \mapsto \ds E[Y^k_q] = {1 \over n}\sum_{j < n}
\left(q \over q + \lambda_j\right)^k,
\qquad k \leq l,
$$
to an estimate of the cumulative distribution function 
$$
F: q \geq 0 \mapsto {1 \over n} \left| \bigl\{j < n : \lambda_j  \leq q \bigr\} \right|
$$
can be made in many different ways
with many different designs for a statistical approach.
In this paper we simply point out 
that for each value of $q \geq q_0$
and for many different graphs,
the few estimates $\hat m_1(q)$,~\dots, $\hat m_l(q)$
(we take $l = 4$ in our applications)
can already lead through low cost maximal entropy principle
to sound predictions on $F(q)$
by themselves, i.e., disregarding all other information
contained in $\hat m_k$, $k \leq l$,
that are all the values of $\hat m_k(q')$ for $q' \neq q$.

To see this, fix again a positive $q \leq 2 \alpha$.
From $s$ samples of $l$ replicas,
we get estimates $\hat m_1(q)$,~\dots, $\hat m_l(q)$
of the first $l$ non-trivial moments $m_k$ of $Y_q$
with a relative Monte Carlo error of order $1 / \sqrt{sn}$.
This means that we can build small $\epsilon_1$,~\dots, $\epsilon_l$
in the large $n$ regime such that
$$
m_k \in [\hat m_k(q) (1 - \epsilon_k), \hat m_k(q) (1 + \epsilon_k)],
\qquad 1 \leq k \leq l,
$$
with probability 0.95 at least.
We claim that in a weak sense
and as a consequence of an equivalence of ensembles
---see Section~\ref{majid} for a precise statement---
if $\hat m_1(q)$,~\dots, $\hat m_l(q)$ form an admissible moment sequence,
then for all 
$$
y \in [a, b] = \left[{q \over q + 2 \alpha}, 1\right]
$$
the overwhelming majority in this large $n$ regime
of the atomic distributions
$$
\nu = {1 \over n} \sum_{j < n} \delta_{y_j},
$$
with $y_j \in [a, b]$ for all $j < n$ and
\begin{equation}\label{ahmed}
	\int_a^b y^k \nu(dy)
	= {1 \over n} \sum_{j < n} y_j^k
	\in [\hat m_k(q) (1 - \epsilon_k), \hat m_k(q) (1 + \epsilon_k)],
	\qquad 1 \leq k \leq l,
\end{equation}
have a tail probability $\nu([y, b])$
that is close to the tail probability $\nu^*([y, b])$
of the continuous probability distribution $\nu^*$ on $[a, b]$
that maximizes the entropy
$$
h(\nu) = -\int_a^b \nu(dx) \log {d\nu \over dx}
$$
under the constraints
\begin{equation}\label{julie}
	\int_a^b x^k \nu(dx) = \hat m_k(q),
	\quad k \leq l.
\end{equation}
Defining for any $\beta = (\beta_1, \dots, \beta_l) \in \ds R^l$
\begin{equation}\label{lucie}
	\nu_\beta(dy)
	= {1 \over \Xi_\beta} \exp\Biggl\{-\sum_{k = 1}^l \beta_k y^k\Biggr\} dy,
	\qquad
	\Xi_\beta
	= \int_a^b \exp\Biggl\{-\sum_{k = 1}^l \beta_k y^k\Biggr\} dy,
\end{equation}
it holds $\nu^* = \nu_{\beta^*}$
with $\beta^*$ the unique minimizer of the smooth and convex function
$$
\beta \in \ds R^l \mapsto \ln \Xi_\beta
+ \sum_{k = 1}^l \beta_k \hat m_k(q),
$$
which can be easily computed by Newton's method in small dimension $l$
(see for example~\cite{MP})
except when $\beta^*$ is large,
i.e., $\nu_{\beta^*}$ is strongly concentrated around a few points.
But in this case Markov bounds can already provides
a good estimation for $F(q)$.
And we can take $\nu_{\beta^*}([1/2, 1])$
as predicted value for $F(q)$ in the generic case.

\begin{figure}[tbp] 
	\caption{\footnotesize%
		Cumulative distribution function (dashed line)
		in natural (left) and log-log (right) scales
		of the spectral measure
		for the bunny graph with 2053 nodes
		and mean degree 52.33
		together with its lower (upward triangles)
		and upper (downward triangles) bounds
		and its predictions
		(yellow vertical segment, green cross, cyan or blue asterisks)
		computed from Monte Carlo estimation of $m_1$, \dots, $m_4$
		after sampling 400 replicated forest trajectories
		up to time $1 / q_0$ with $q_0 = \bar\lambda / 100$
		and with 1 (yellow), 2 (green), 3 (cyan) or 4 (blue)
		valid moment estimates.
		\label{jean-francois}
	}
	\sbox0{\includegraphics{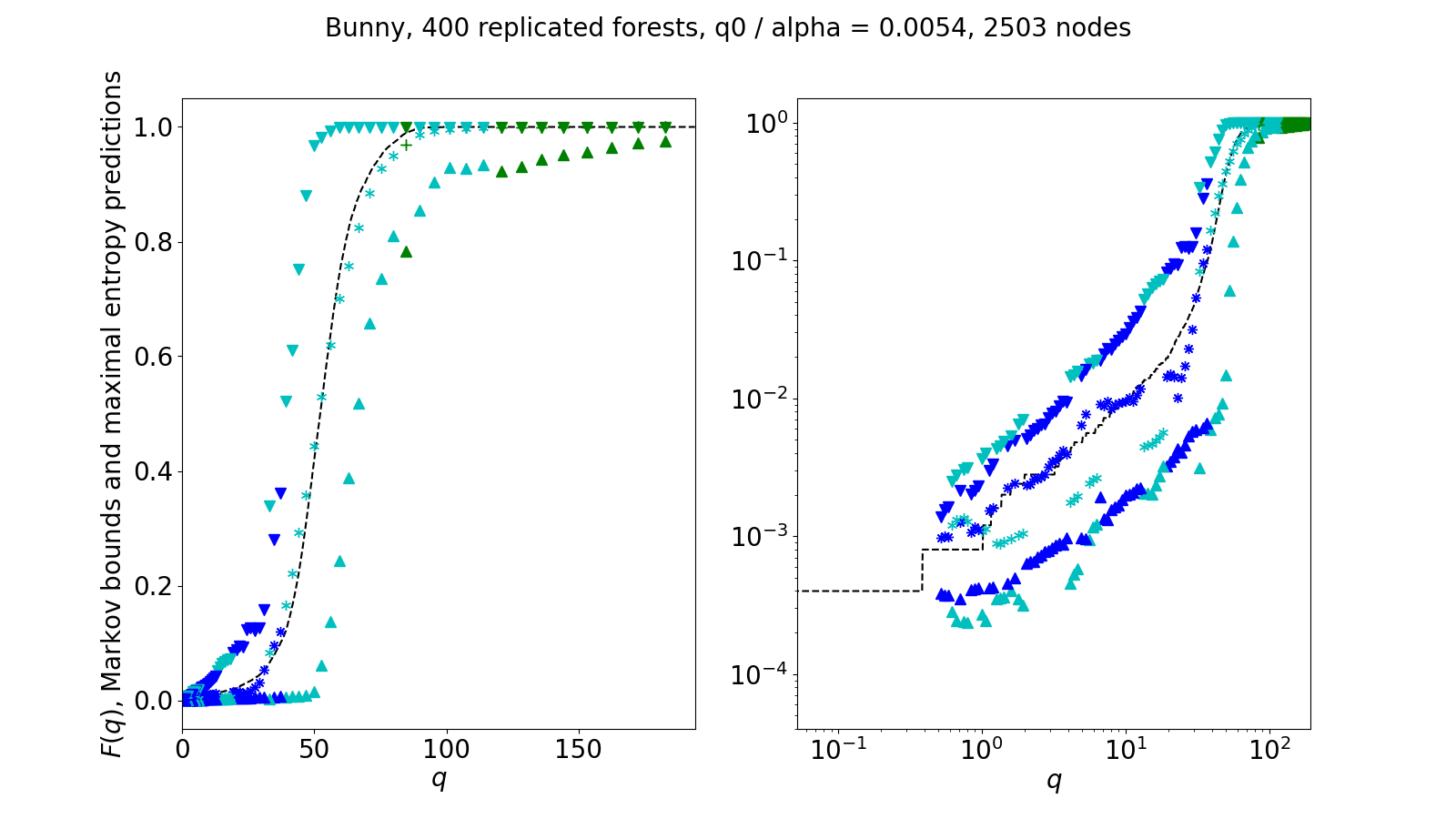}}
	\hbox to \hsize{\hfill\includegraphics[clip, trim=0 {.02\ht0} 0 {.06\ht0}, width=5.6in]{./pictures/Bunny.png}\hfill}
\end{figure}

Figure~\ref{jean-francois} adds to Figure~\ref{brune}
the predicted values $\nu_{\beta^*}([1 / 2, 1])$ for $F(q)$.
We found them much less disappointing
than the previous Markov bounds.
More numerical experiments are presented in Section~\ref{madeleine}.

We will argue in Section~\ref{madeleine}
that, to be consistent with the relative errors $\epsilon_k$
of order $1 / \sqrt{sn}$ for $s$ sampled trajectories,
the computational cost of these predictions
is of order $\sqrt{sn}$ only.
It is then negligible with respect to the numerical cost
for estimating $m_1$,~\dots, $m_l$.

\subsection{Equivalence of ensembles}\label{majid}
The material of this section is rather standard. It is intended to  make precise
our previous claim on the ``overwhelming majority''
of the atomic measures $\nu$ that satisfy the constraint~\eqref{ahmed}.
Let  $\bar V_n$ be the set 
of atomic measures on a interval $[a, b]$
$$
\nu = {1 \over n} \sum_{j < n} \delta_{y_j}
$$
that satisfy the constraint
\begin{equation}\label{salia}
	\int_a^b y^k \nu(dy) \in [m_k(1 - \epsilon_k), m_k(1 + \epsilon_k)],
	\qquad k = 1, \dots, l,
\end{equation}
for an admissible moment sequence $m_1$,~\dots, $m_l$
and small $\epsilon_1$,~\dots, $\epsilon_l > 0$. 

For making sense of the notion of  ``overwhelming majority of elements of $\bar V_n$'',
we first have to define a probability distribution 
on the set $\bar V_n$. 
i.e., to define a random atomic measure $\tilde \zeta$
that satisfies the moment constraints~\eqref{salia}.
We will then be able to check that, for all $\epsilon > 0$ and all $y \in [a,b]$, 
\begin{equation}\label{pablo}
	\lim_{n \rightarrow \infty} P\bigl(
	|\tilde\zeta([y, b]) - \nu^*([y, b])| > \epsilon
	\bigr) = 0.
\end{equation}

\noindent\textbf{Building $\tilde\zeta$ and checking~\eqref{pablo}.}
Let $\mu_{\beta^*}$ be the canonical Gibbs measure on $[a, b]^n$ defined by
$$
\mu_{\beta^*}(dz)
= {1 \over \Xi_{\beta^*}^n} \exp\Biggl\{-\sum_{k = 1}^l \beta^*_k
\sum_{j < n} z_j^k\Biggr\} dz,
$$
with $\Xi_{\beta^*}$ the partition function of $\nu_{\beta^*}$
defined in~\eqref{lucie}.
We denote by $\zeta = (\zeta_j)_{j < n}$
the random variable in $[a, b]^n$ with law $\mu_{\beta^*}$.
The coordinates $\zeta_j$, $j < n$,
are independent random variables with the same law $\nu_{\beta^*}$ on $[a, b]$.
With
$$
V_n = \Biggl\{
z \in [a, b]^n :
\forall 1 \leq k \leq l,\:
{1 \over n} \sum_{j < n} z_j^k
\in [m_k (1 - \epsilon_k), m_k (1 + \epsilon_k)]
\Biggr\},
$$
it holds then
$$
\lim_{n \rightarrow \infty} P(\zeta \in V_n) = 1
$$
by the weak law of large numbers.
By associating with each $z \in [a, b]^n$
a probability 
$$
\bar z = {1 \over n} \sum_{j < n} \delta_{z_j}.
$$
on $[a, b]$
we build the random atomic measure $\tilde\zeta$
by conditioning $\bar\zeta$ to
$$
\{\bar\zeta \in \bar V_n\} = \{\zeta \in V_n\}.
$$
While $\zeta$ conditioned to $\{\zeta \in V_n\}$
is not distributed according to the Lebesgue measure on $V_n$,
the microcanonical ensemble obtained 
by a further conditioning on $\sum_{j < n} \zeta_j^k$, $k \leq l$,
can be seen as a ``uniform distribution''
on our atomic measures with moments $\sum_{j < n} \zeta_j^k / n$, $k \leq l$.
Together with the fact that we are working with small $\epsilon_k$,
the law of $\tilde\zeta$ does then provide
a reasonable statistical model on $\bar V_n$.
Within this model, Equation~\eqref{pablo}
immediately follows from 
\begin{equation*}
	\begin{split}
		P\bigl(|\tilde\zeta([y, b]) - \nu_{\beta^*}([y, b])| > \epsilon\bigr)
		&= {P\bigl(|\bar\zeta([y, b]) - \nu_{\beta^*}([y, b])| > \epsilon,
			\zeta \in V_n\bigr)
			\over P\bigl(\zeta \in V_n\bigr)} \\
		&\leq {P\bigl(|\bar\zeta([y, b]) - \nu_{\beta^*}([y, b])| > \epsilon\bigr)
			\over P\bigl(\zeta \in V_n\bigr)}
	\end{split}
\end{equation*}
and the weak law of large number again.
\qed

\section{Proofs of Theorems~\ref{ben} and~\ref{clo}}
\label{bertrand}

\subsection{Proof of Theorem~\ref{ben}}

Recall that, for a given node $x$ and forest $\phi$,
$\rho_x(\phi)$ is the root of the tree that covers $x$ in $\phi$.
By using again the abelianity of Wilson's algorithm,
we can assume $x$ to be our first starting point
for the loop-erased random walks.
Consider the discrete time random walk $(\tilde X_k)_{k \geq 0}$
on $\cl X\cup \{\Delta \}$, where $\Delta$ is a cemetery state,
with transition probabilities
\begin{equation}
	\label{skeleton.eq}
	\tilde P(x,y)=\begin{cases}
		\frac{w(x,y)}{q+w(x)}&\text{if }x\neq y,\:x,y\neq \Delta\\
		\frac{q}{q+w(x)} &\text{if }x\neq\Delta ,  y=\Delta\\
		0&\text{if }x= y\neq\Delta\\
		1&\text{if }x=y=\Delta.
	\end{cases}
\end{equation}
Call $T_\Delta$ the absorption time of the chain in $\Delta$.
For all $x$ and $y$ in $\cl X$,
$$
\begin{aligned}
	\ds P\bigl(\rho_x(\Phi_q) = y\bigr) 
	& = \ds P_x \bigl(\tilde X(T_{\Delta}-1) =y \bigr) 
	= \sum_{k \geq 1}   \ds P_x  \bigl(\tilde X_{k-1} =y, \tilde X_k = \Delta \bigr)  
	\\
	& = {q \over q + w(y)} \sum_{k \geq 1} \tilde P_{\cl X}^{k-1}(x,y)
	=  {q \over q +  w(y)} \biggl({\rm Id} - \tilde P_{\cl X} \biggl)^{-1} (x, y)
	\\
	& = q \biggl({\rm Id} - \tilde P_{\cl X} \biggl)^{-1}
	\biggl(q {\rm Id} +  W  \biggl)^{-1}(x,y), 
\end{aligned}
$$
where $\tilde P_{\cl X}$ is the square submatrix of  $\tilde P$ restricted to $\cl X$.
Note that $\tilde P_{\cl X} = \biggl(q {\rm Id} +W\biggl)^{-1} (L + W)$,
so that 
$$
\biggl(q {\rm Id} +  W  \biggl)  \biggl({\rm Id} - \tilde P_{\cl X} \biggl)
= q {\rm Id} +  W - L - W = q{\rm Id} - L  . 
$$
Setting $K_q = q(q{\rm Id} - L)^{-1}$,
this gives $\ds P\bigl(\rho_x(\Phi_q) = y\bigr) = K_q(x,y)$
and therefore, for all $k \geq 0$,
$$
\ds P\bigl(R^k(x) = y\bigr) = K_q^k(x, y). 
$$
Hence, 
$$
\ds E\bigl[|\xi^k_q|\bigr]
= \sum_{x \in \cl X} \ds P\bigl(x \in \xi^k_q \bigr)
= \sum_{x \in \cl X} \ds P\bigl(R^k(x) = x \bigr)
= {\rm Tr}\bigl(K_q^k\bigr)
= \sum_{j < n} \biggl({q \over q + \lambda_j}\biggr)^k.
$$

Let us now think of the forests $\Phi_{k, q}$, $k < l$,
as forests on different layers that are $l$ distinct copies of $\cl X$.
Since negative correlations are relevant for positive $k$ only,
we choose $k \geq 1$ and observe that for checking
the occurrence of $\bigl\{x, y \in \xi^k_q\bigr\}$ with $x \neq y$,
we can first sample the multilayer path leading from $x$ to $R^k_q(x)$,
with one tree branch in each of the first $k$ layers,
before sampling the multilayer path leading from $y$ to $R^k_q(y)$.
Negative correlations simply follow from the fact that
conditionally on $\bigl\{x \in \xi^k_q\bigr\}$,
in order to have $\bigl\{y \in \xi^k_q\bigr\}$,
the second path must avoid the first one,
since otherwise we would have
$$
R^k_q(y) = R^k_q(x) = x \neq y.
$$
Using the negative correlations property, the variance bound is standard:
$$
\begin{aligned}
	{\rm Var}\bigl(|\xi_q^k|\bigr)
	&= \sum_x \ds P\bigl( x \in \xi_q^k)
	+ \sum_{x \neq y} \ds P\bigl(x \in \xi_q^k, y \in \xi_q^k \bigr) 
   - \sum_x \ds P\bigl( x \in \xi_q^k)^2
	- \sum_{x \neq y} \ds P\bigl(x \in \xi_q^k\bigr) 
	\ds P\bigl(y \in \xi_q^k \bigr) \\
	&   \leq \sum_x \ds P\bigl( x \in \xi_q^k) = \ds E\bigl[|\xi_q^k| \bigr] \, .                                     
\end{aligned}
$$
\qed

\subsection{Proof of Theorem~\ref{clo}}

Using the stack representation,
the arrows used in the stack of $x$ in order to construct $\Phi_q$
correspond to marks $U > \frac{q}{q+w(x)}$.
Since  for $q \ge q'$, $\frac{q}{q+w(x)} \ge \frac{q'}{q'+w(x)}$,
the same arrows are used to construct $\Phi_{q'}$.
Hence, the number $S$ of sampled arrows necessary to build the whole trajectory 
$$
t \in [0, 1 / q_0] \mapsto\Phi_{1 / t},
$$
is exactly the number of sampled arrows
necessary to build $\Phi_{q_0}$ with Wilson's algorithm.\\

We now prove the identity \nt{(we include it here for completeness even though it is a consequence of Proposition 1 of~\cite{Ma})}
\begin{equation}\label{tonino}
	\ds E[S]={\rm Tr}\bigl( (q_0 {\rm Id} -L)^{-1}(W+q_0 {\rm Id})\bigr).
\end{equation}
Consider again the discrete time random walk $(\tilde X_k)_{k \geq 0}$,
whose transition probabilities are
given in~\eqref{skeleton.eq}.
For every $x\in \cl X$,
$\ell(x)=\sum_{k=0}^{T_\Delta - 1} \ds 1_{\{\tilde X_k=x\}}$, 
is the local time spent at $x$ before absorption.
We denote by $M(x,y)$ the expectation of $\ell(y)$ 
when the starting node is $x$, i.e. 
$$
M(x,y)=\ds E_x[\ell(y)].
$$ 
If we start Wilson's algorithm at $x$,
the expected number of arrows in the stack of $x$
that we need to read is given by $\ds E_x[\ell(x)]$.
Therefore, by linearity of expectation, what we need to compute is 
$$
\ds E[S]={\rm Tr}(M).
$$
In order to prove \eqref{tonino}
it is then sufficient to show that
$$
M = (q_0 {\rm Id} - L)^{-1}(W+q_0 {\rm Id}).
$$
Let us notice that,
for every couple $x,y\in \cl X$ we have
$$
M(x,y) = \ds 1_{\{x = y\}}
+\sum_{z \in\cl X\setminus\{x\}}{w(x,z) \over q_0+w(x)}M(z,y),
$$
thanks to Markov property at time $1$.
Multiplying by $(q_0+w(x))$ both sides of the equality we get
$$
\begin{aligned}
	q_0 M(x,y)
	&= \ds 1_{\{x = y\}} (q_0 +w(x))
	- w(x)M(x,y)+\sum_{z\neq x}L(x,z)M(z,y) = \ds 1_{\{x=y\}} (q_0 +w(x))  + LM(x,y),
\end{aligned}
$$
which reads $(q_0 {\rm Id} -L)M =(q_0 {\rm Id} + W)$, 
from which \eqref{tonino} immediately follows.\\

\nt{We will now show the upper bound:
\begin{equation}\ds E[S] \leq n\left(1+\frac{\beta}{q_0}\right) ~~ \text{with} ~~ \beta=\frac{1}{n}\sum_x \text{max}_y w(y,x)
	\label{eq:ubS}
\end{equation}
Using the previous expression for $M$, we have:
\begin{align*}
	\ds E[S] &= {\rm Tr}(M) = {\rm Tr}\left((q_0 {\rm Id} - L)^{-1}(q_0 {\rm Id}+W)\right)= {\rm Tr}\left((q_0 {\rm Id} - L)^{-1}(q_0 {\rm Id}-L+L+W)\right)\\
	&= n + {\rm Tr}\left((q_0 {\rm Id} - L)^{-1}(L+W)\right)
\end{align*}
Looking into this last trace, and denoting by $K$ the matrix $K=q_0(q_0 {\rm Id} - L)^{-1}$, one obtains
\begin{align*}
{\rm Tr}\left((q_0 {\rm Id} - L)^{-1}(L+W)\right) &= \frac{1}{q_0}{\rm Tr}\left(K(L+W)\right)=\frac{1}{q_0}\sum_x \sum_y K(x,y) w(y,x)\\
&\leq \frac{1}{q_0}\sum_x \text{max}_{y} w(y,x) \sum_y K(x,y)
\end{align*}
as all entries of $K$ are non-negative ($K$ being an M-matrix). 
Now, as the vector of all-ones $1_n$ is an eigenvector of $K$ verifying\footnote{as $K$ is diagonalizable in the eigenbasis of $L$ and $L 1_n = 0$} $K 1_n = 1_n$, one has that for all $x$, $\sum_y K(x,y)=1$ yielding the claimed upper-bound \eqref{eq:ubS}.}\\

We now turn to $R$. 
Since each root $x$ of $\Phi_{1 / t}$
is unfrozen at rate $w(x) / (1 + w(x)t)$,
and when unfreezing $x$ a time $t$
the number extra reads is bounded
by the size $|\tau_x(\Phi_{1 / t})|$
of the tree that covers $x$ in $\Phi_{1 / t}$,
it holds
$$
\begin{aligned}
	\ds E[R]
	&\leq \ds E\left[\int_0^{1/q_0}
	\sum_{x \in \rho(\Phi_{1/t})} \left|\tau_x(\Phi_{1/t})\right|
	{w(x) \over 1 + w(x)t} {\rm dt} \right] \\ 	
	&\leq \ds E \left[\int_0^{1/q_0}  \frac{\alpha}{1 + \alpha t}  
	\sum_{x \in \rho(\Phi_{1/t})} \left|\tau_x(\Phi_{1/t})\right|
	{\rm dt} \right]= n \int_0^{1/q_0} \frac{\alpha}{1+ \alpha t} {\rm dt}
	= n \ln \left(1+ \frac{\alpha}{q_0} \right).
\end{aligned}
$$
\qed

\section{Numerical issues and results}
\label{madeleine}

\subsection{Subsampling the Stieltjes transform and its first derivatives}

Our coupled forest algorithm defined and described
in Section~\ref{fiorella}, Section~\ref{farah} and Appendix~\ref{sacha}
allows to build a single piecewise constant random forest trajectory
$(\Phi_q)_{2\alpha \geq q \geq \epsilon_0\bar\lambda}$
at a sampling cost in
$O\left( n \left(\frac{\nt{\beta}}{\epsilon_0\nt{\bar\lambda}}
+ \ln \left(\frac{\alpha}{\epsilon_0 \bar\lambda}\right)\right)\right)$. 

For each value of $q$ in $[\epsilon_0 \bar\lambda, 2\alpha]$,
the estimates of the Stieltjes transform and its derivatives given by \eqref{MC.eq}
are based on $s$ samples of $l$ replicas of such trajectories.
$l$ is typically small ($l=4$ in our experiments).
For each sample, $\eqref{MC.eq}$ requires
the computation of the size of $\xi^1_q$,~\dots, $\xi^l_q$ defined in~\eqref{felix}. 
This counting cost is in $O(n)$, for each sample and each $q$.

As a consequence we decided to make these estimates
in $\lceil 1 / \epsilon_0\rceil$ values of $q$,
naturally organized in a geometric progression
$q_i = q_0 r^i$, $i < 1 / \epsilon_0$,
with 
$$
r  = \exp\left\{\epsilon_0 \ln {2\alpha \over q_0}\right\},
$$
so that $q_0 r^{1 / \epsilon_0} = 2\alpha$.
\nt{The total sampling and counting cost
for all these values of $q$
is then in $O\left( n \left(\frac{\beta}{\epsilon_0\bar\lambda}
+ \ln \left(\frac{\alpha}{\epsilon_0 \bar\lambda}\right)\right)\right)$
for each replicated coupled forest trajectory. Let us further simplify this computation cost. 
First, since
$$
\alpha/n  \le \bar\lambda \le \alpha
$$
(recall that  $n \bar\lambda = \mbox{Tr}(-L) = \sum_{x \in \cl X} w(x)$),
one has $\ln \left(\frac{\alpha}{\epsilon_0 \bar\lambda}\right) \leq \ln \left(\frac{n}{\epsilon_0}\right)$. Second, one has $\beta\leq \bar\lambda$, such that $\frac{\nt{\beta}}{\epsilon_0\nt{\bar\lambda}}\leq \frac{1}{\epsilon_0}$. 
An upper-bound of the cost is thus $O\left(\frac{1}{\epsilon_0}n\ln n\right)$ for 
each replicated coupled forest trajectory. Note that this upper-bound can be quite crude. For instance, in an unweighted graph (for which $\beta=1$) that is additionally regular ($\alpha=\bar{\lambda}$), this cost can even drop to $O\left( \frac{n}{\epsilon_0\bar\lambda}+n\ln \frac{1}{\epsilon_0}\right)$. }

\subsection{Newton's descent algorithm,
	numerical integration and valid moment estimates}

At a total computational cost in $O(s\epsilon_0^{-1}n\ln n)$
and for all $q = q_i$, $i \leq 1 / \epsilon_0$,
we get the estimates $\hat m_k(q)$ of the moments
$$
m_k = {1 \over n} \sum_{j < n} \left(q \over q + \lambda_j\right)^k,
\qquad k \leq l,
$$
with a relative error of order $1 / \sqrt{sn}$
associated with a Monte Carlo confidence interval.
In order to find the maximal entropy estimator of $F(q)$
at a given $q$ in $[q_0, 2\alpha]$,
we have to perform a numerical integration to compute
the low dimensional Hessian to be inverted in Newton's method. 
This can then be made with the same relative precision
and has therefore a computational cost in $O(\sqrt{sn})$,
provided that we can ensure fast convergence of Newton's method,
which requires to address two related issues. 

First, the constraint~\eqref{julie}
must be satisfiable by some measure
$\nu$ on $[a, b]$.
By using Lemma~\ref{frank} from Appendix~\ref{yasmina}
we can check recursively that 
$\hat m_1(q)$, \dots, $\hat m_k(q)$,
$k \leq l$, form an admissible moment sequence :
$\hat m_1(q)$ must lie in $[a, b]$
and $\hat m_k(q)$ has to be in an explicit interval $I_k$
computed from $\hat m_1(q)$,~\dots, $\hat m_{k - 1}(q)$. 
Since each estimate $\hat m_k(q)$
comes with its confidence interval,
our prediction will use only 
the longest sequence $\hat m_1(q)$,~\dots, $\hat m_k(q)$
for which $\hat m_1(q)$,~\dots, $\hat m_{k - 1}(q)$, $\tilde m_k$
form an admissible moment sequence for all $\tilde m_k$
in the Monte Carlo 95\% confidence interval associated with $\hat m_k(q)$.

Second, 
for such an admissible moment sequence $\hat m_1(q)$, \dots, $\hat m_k(q)$,
despite the fact that our maximal entropy estimator 
with $k$ non-trivial estimated moment is well-defined,
we should expect to run into numerical difficulties
in the case when our estimator would be associated 
with a large $\beta^*$, solution of
\begin{equation}\label{lea}
	\int_a^b y^j \nu_{\beta^*}(dy) =  \hat m_j(q),
	\qquad j \leq k.
\end{equation}
But in this case Markov's bounds on $F(q)$,
which are given in Appendix~\ref{louise}
and explicitly computed from $\hat m_1(q)$, \dots, $\hat m_k(q)$,
can already provide a good estimate of $F(q)$.
In practice, we will make use of Newton's method
only in the generic case
when Markov's bounds does not provide
an estimate on $F(q)$ with a relative error of 1\% or less.

In order to minimize the iteration number 
of each Newton's descent algorithm
to be run for a given value $q_i$
and with $k$ non-trivial moment estimated,
we will start the descent from the parameter $\beta^*$
computed at the previous value $q_{i + 1} > q_i$ with the same $k$,
if this parameter is available.
When this is not the case,
we use the parameter $\beta^*$ computed
with one moment less at the same $q_i$.
We will also fix a maximal iteration number equal to 50
to deal with the exceptional cases
when Markov's bounds are not sufficient for providing
a good estimate of $F(q)$ 
but Newton's algorithm does not reach 
an approximate solution of~\eqref{lea}
in a few iterations only.
This ensures that the extra computational cost
for giving Markov bounds
and computing the maximal entropy estimator of $F(q_i)$,
for all $i < 1 / \epsilon_0$, is in $O(\epsilon_0^{-1}\sqrt{sn})$
which is negligible with respect to the sampling cost
$O(\epsilon_0^{-1}sn\ln n)$.

We can finally precisely define the number of valid moment estimates
that appears for each value of $q$ through the different colors
in Figure~\ref{brune} and Figure~\ref{jean-francois}
as in the experimental results presented below.
For all $k \leq l$,
we say that $\hat m_1(q)$,~\dots, $\hat m_k(q)$
is a valid moment estimate 
\begin{itemize}
	\item if so is $\hat m_1(q)$,~\dots, $\hat m_{k - 1}(q)$; 
	\item if $\hat m_1(q)$,~\dots, $\hat m_{k - 1}(q)$, $\tilde m_k$
	form an admissible moment sequence for all $\tilde m_k$
	in the Monte Carlo 95\% confidence interval associated with $\hat m_k(q)$;
	\item and if $\hat m_1(q)$,~\dots, $\hat m_k(q)$ either provides 
	an estimate of $F(q)$ within a relative error of less than 1\%
	through Markov bounds or leads to the computation
	of its maximal entropy estimator in less than 50 iterations
	of Newton's algorithm.
\end{itemize}

\subsection{Experimental results}

In all pictures of this section,
we follow the same conventions as in Figure~\ref{jean-francois}
to show our experimental results.
Except when otherwise specified,
we considered unweighted graphs
with $n = 5000$ nodes, we took $\epsilon_0 = q_0 / \bar \lambda = .01$
and we sampled for each graph 400 coupled replicated forest trajectories
with $l = 4$ replicas, i.e., 1600 coupled forest trajectories
$(\Phi_q)_{2\alpha \geq q \geq q_0}$.

\subsubsection{What went right}

\begin{figure}[tbp]
	\caption{\footnotesize Two sparse Erdős-Rényi graphs.\label{laurent}}
	\sbox0{\includegraphics{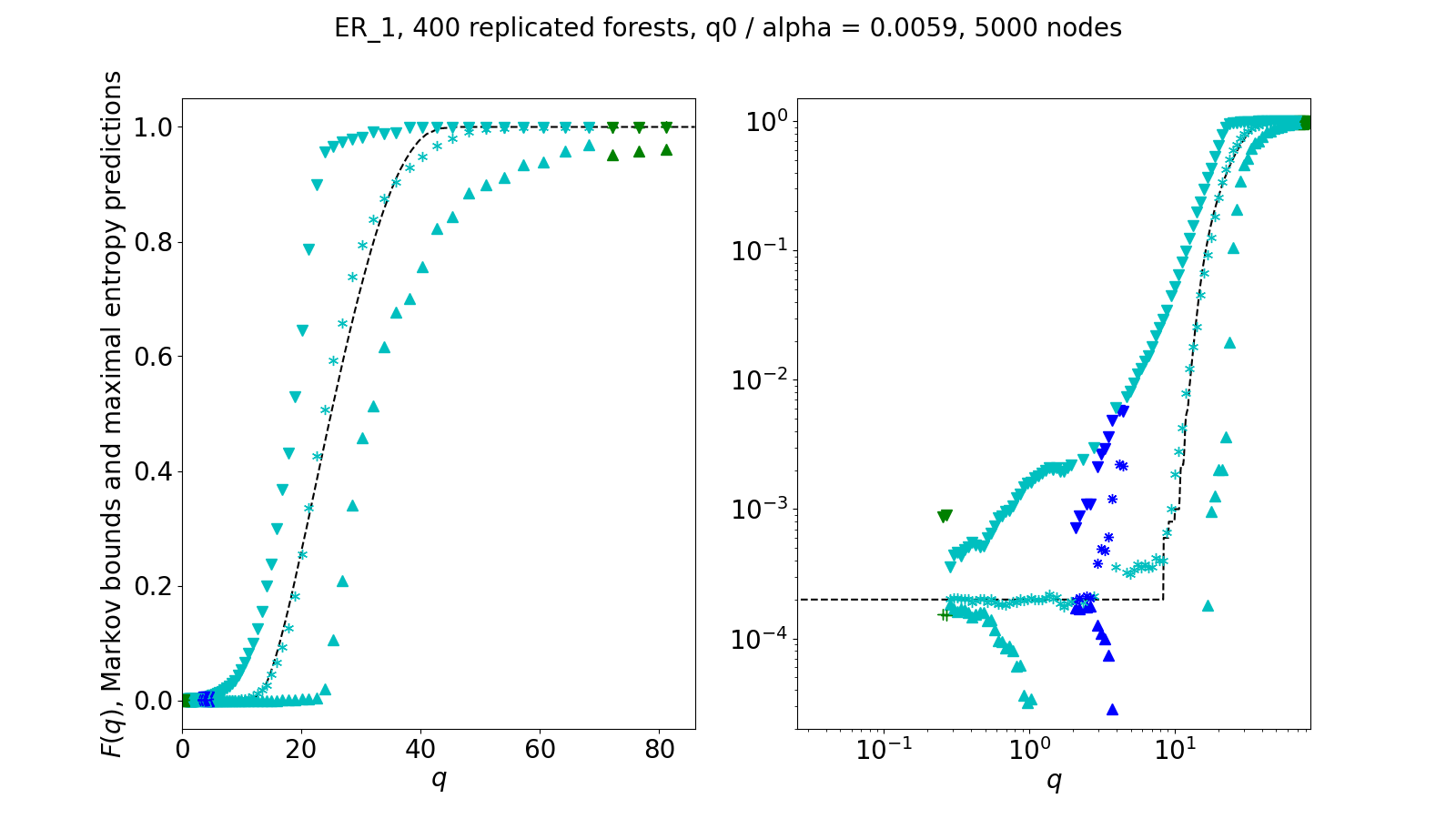}}
	\hbox to \hsize{%
		\includegraphics[clip, trim={.06\wd0} {.03\ht0} {.09\wd0} {.06\ht0}, width=2.8in]{./pictures/ER_1.png}%
		\hfill%
		\includegraphics[clip, trim={.06\wd0} {.03\ht0} {.09\wd0} {.06\ht0}, width=2.8in]{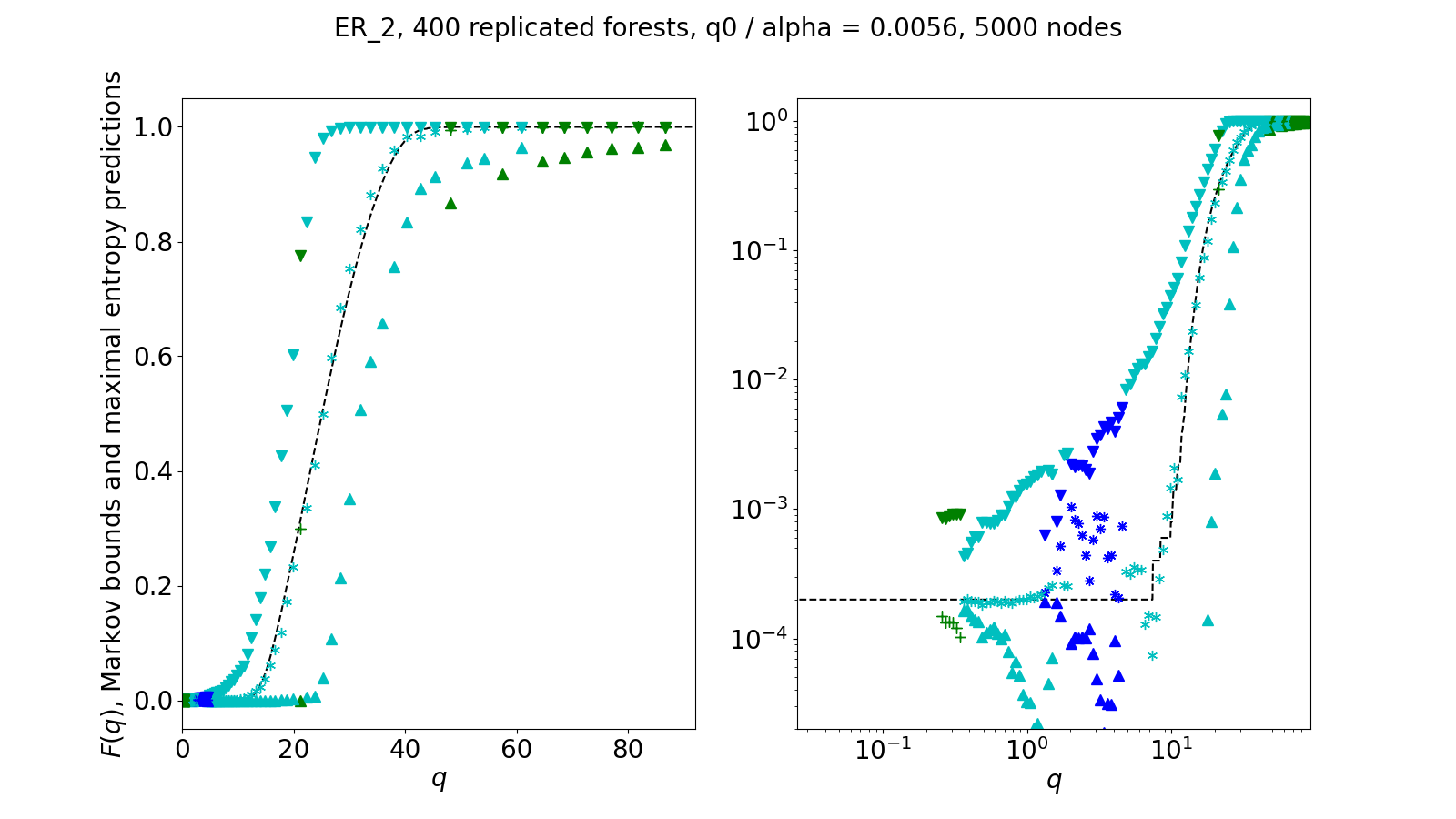}%
	}
\end{figure}

\begin{figure}[tbp]
	\caption{\footnotesize Two denser Erdős-Rényi graphs.\label{marie-laure}}
	\sbox0{\includegraphics{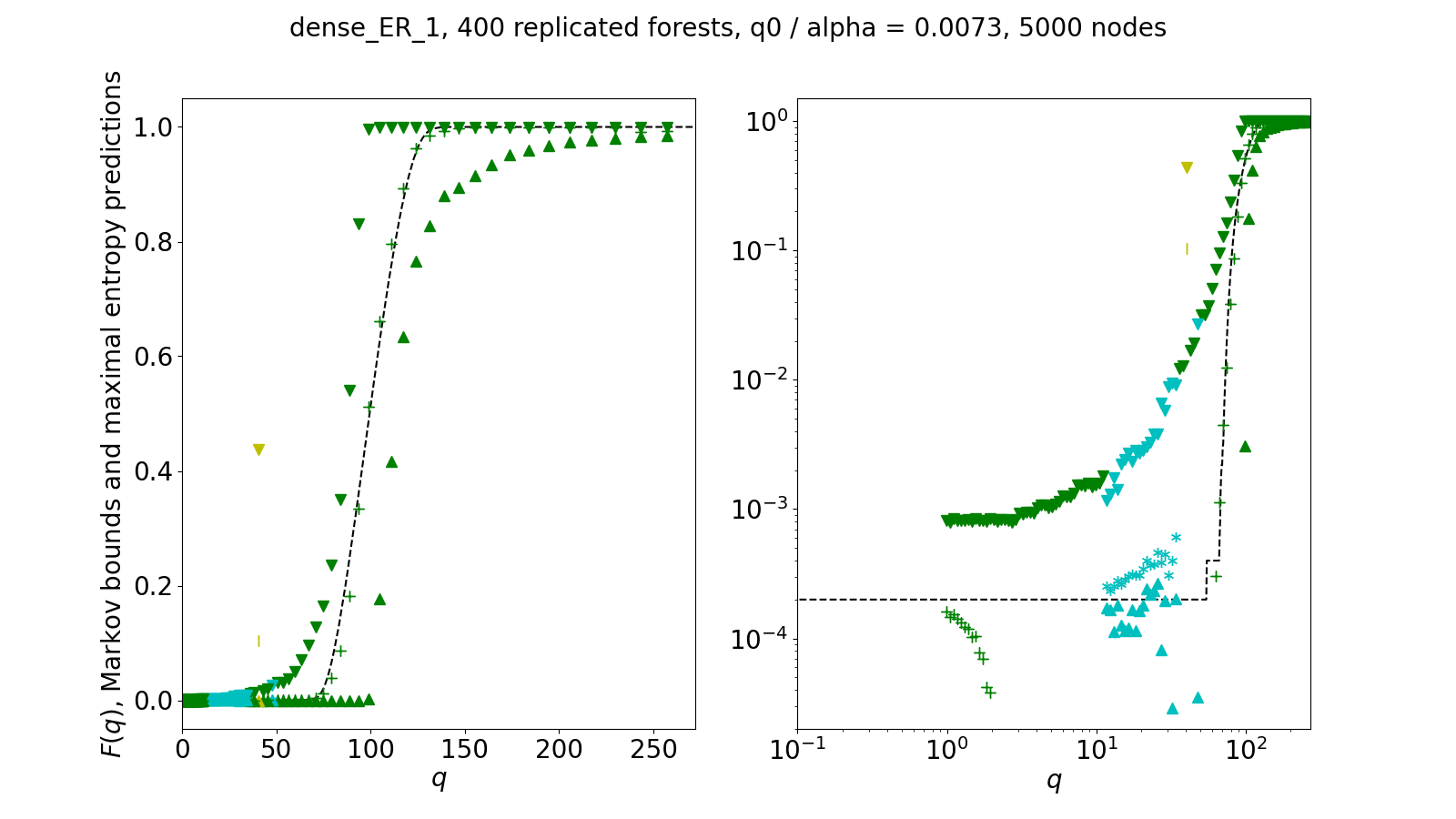}}
	\hbox to \hsize{%
		\includegraphics[clip, trim={.06\wd0} {.03\ht0} {.09\wd0} {.06\ht0}, width=2.8in]{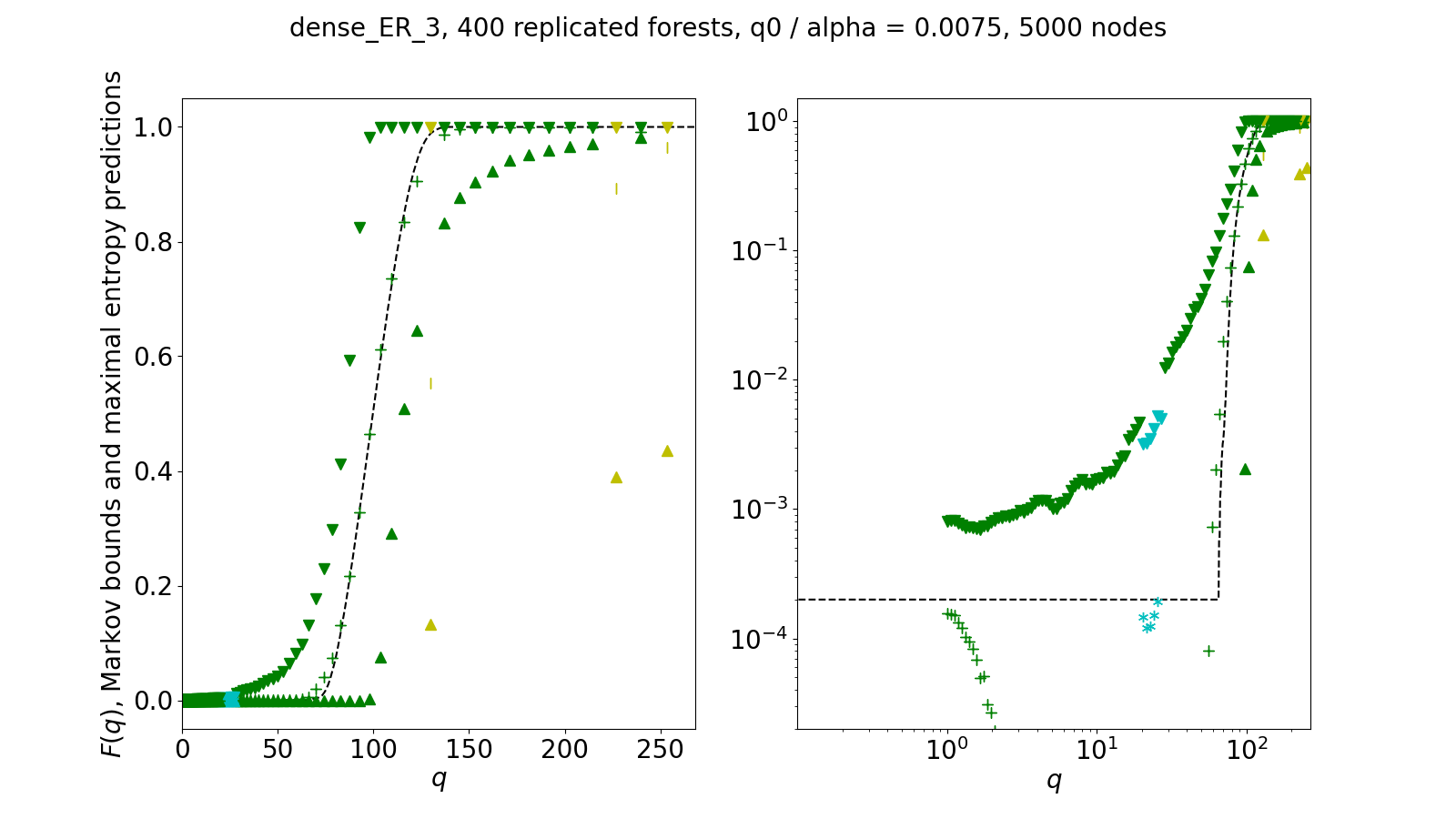}%
		\hfill%
		\includegraphics[clip, trim={.06\wd0} {.03\ht0} {.09\wd0} {.06\ht0}, width=2.8in]{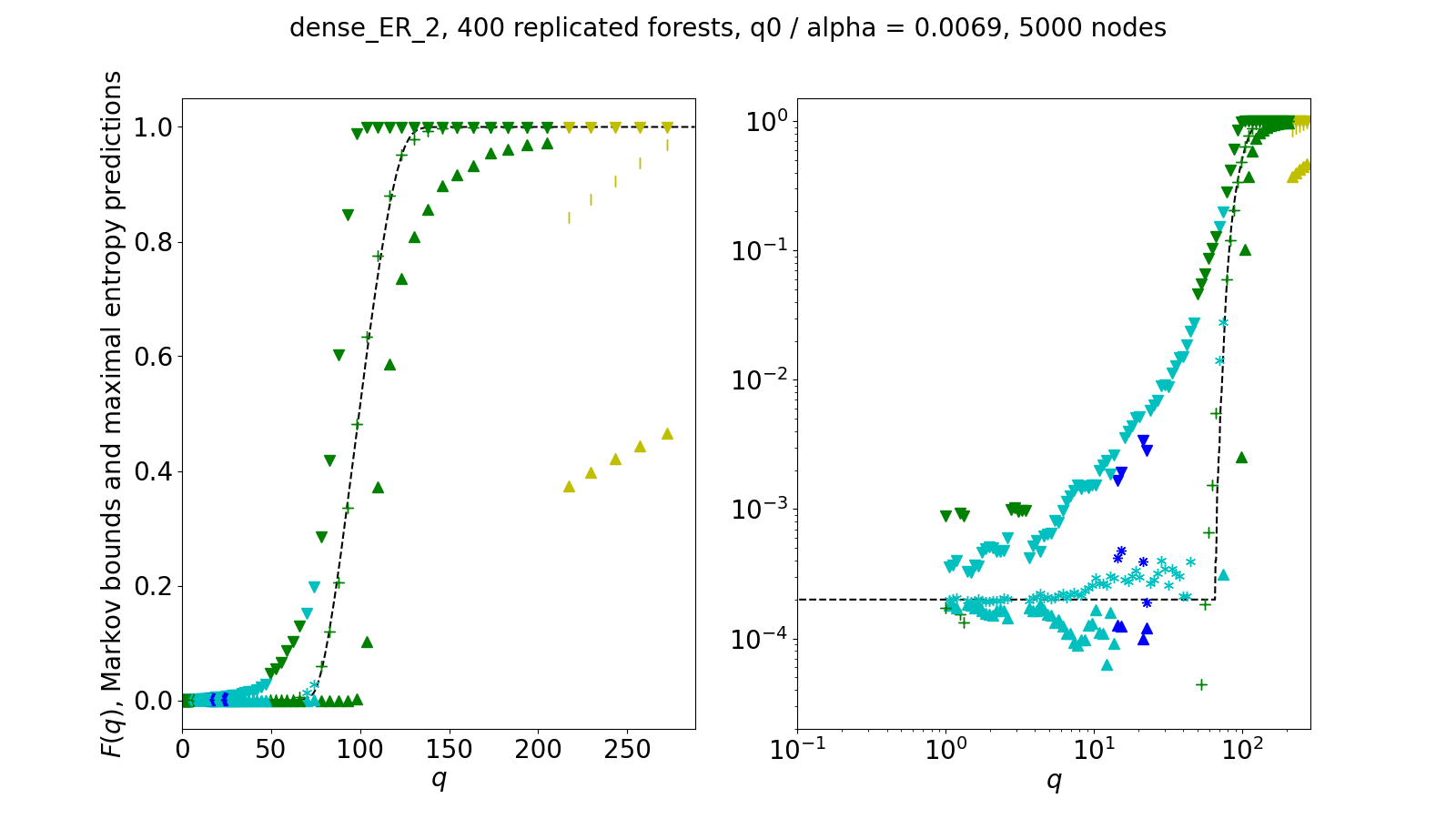}%
	}
\end{figure}

Figure~\ref{laurent} shows our results
for two different Erdős-Rényi graphs of size $n = 5000$
with small density $p = (3 \ln n) / n$;
Figure~\ref{marie-laure} shows our results
for two different Erdős-Rényi graphs with mean degree 100.

\begin{figure}[tbp]
	\caption{\footnotesize Two preferential attachment networks.\label{albane}}
	\sbox0{\includegraphics{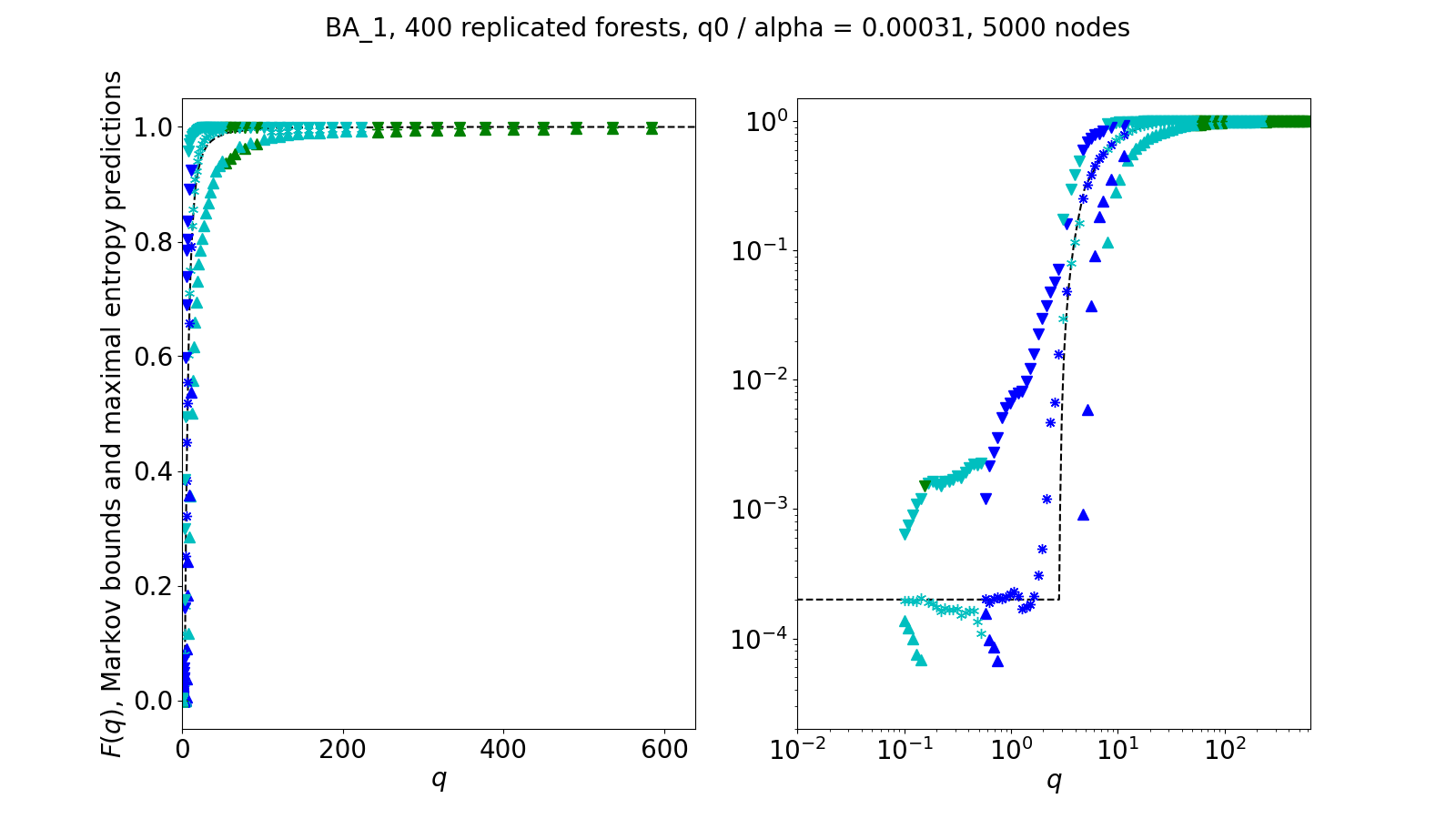}}
	\hbox to \hsize{%
		\includegraphics[clip, trim={.06\wd0} {.03\ht0} {.09\wd0} {.06\ht0}, width=2.8in]{./pictures/BA_1.png}%
		\hfill%
		\includegraphics[clip, trim={.06\wd0} {.03\ht0} {.09\wd0} {.06\ht0}, width=2.8in]{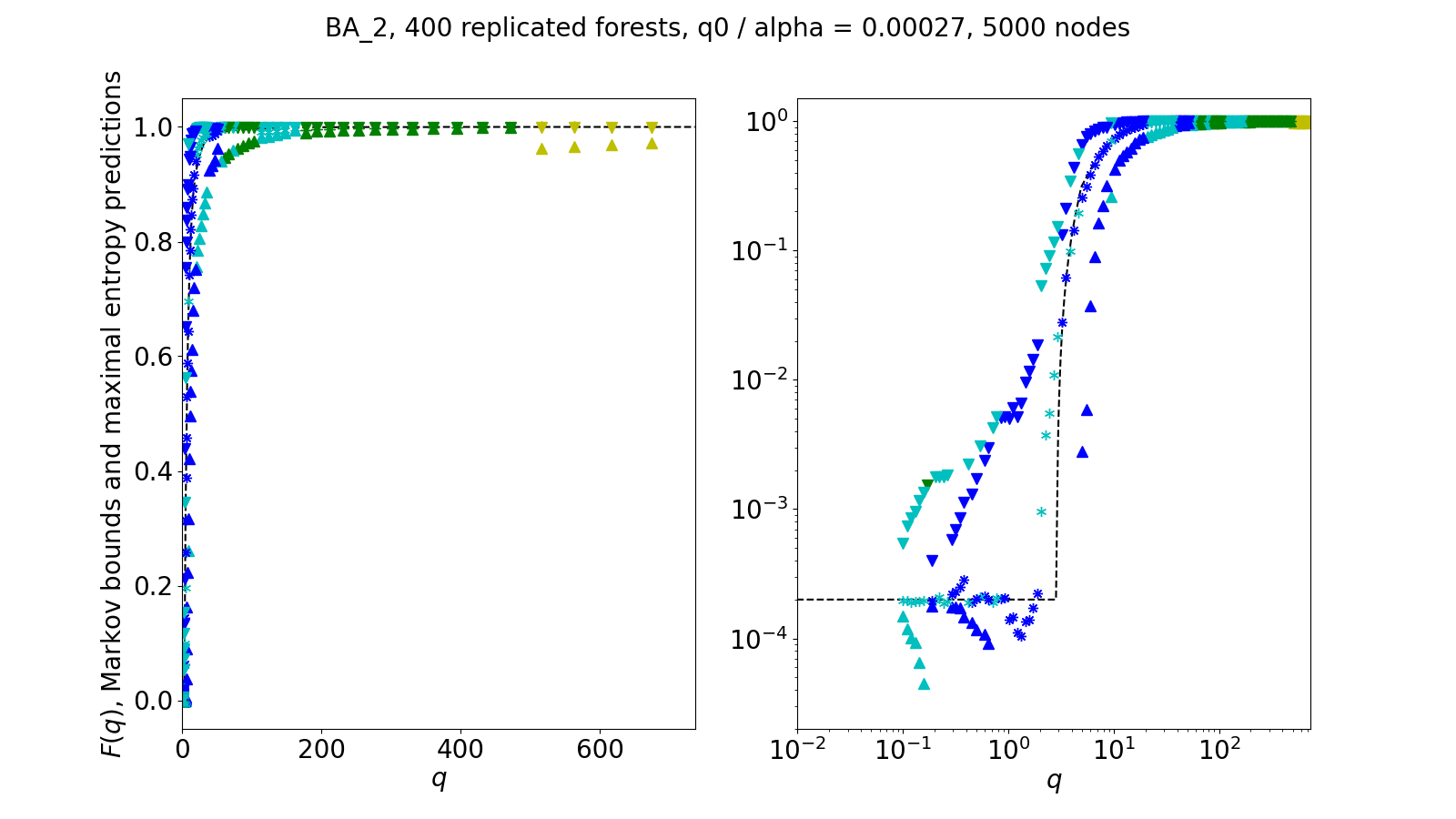}%
	}
\end{figure}

\begin{figure}[tbp]
	\caption{\footnotesize Two stochastic bloc models.\label{celine}}
	\sbox0{\includegraphics{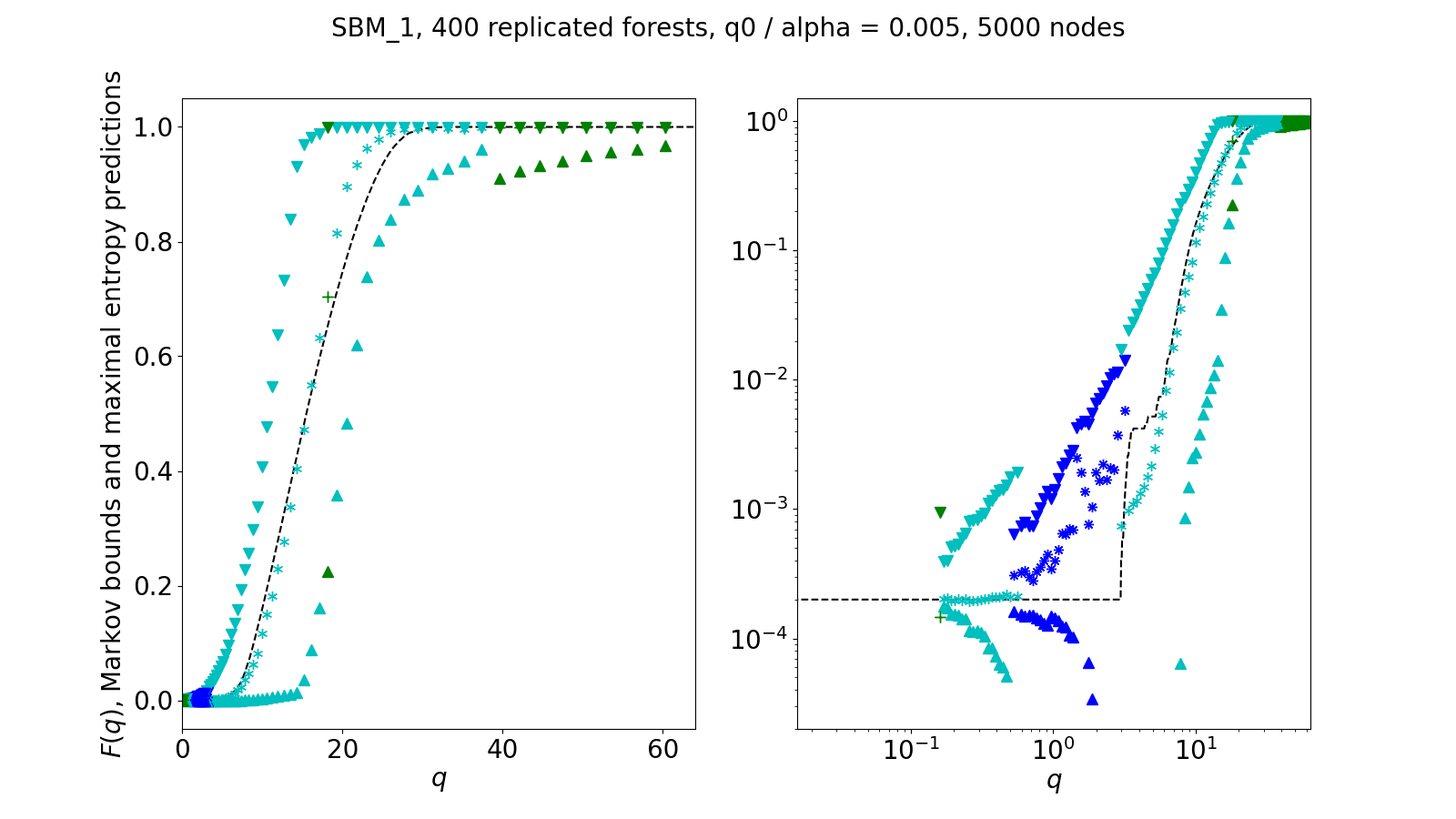}}
	\hbox to \hsize{%
		\includegraphics[clip, trim={.06\wd0} {.03\ht0} {.09\wd0} {.06\ht0}, width=2.8in]{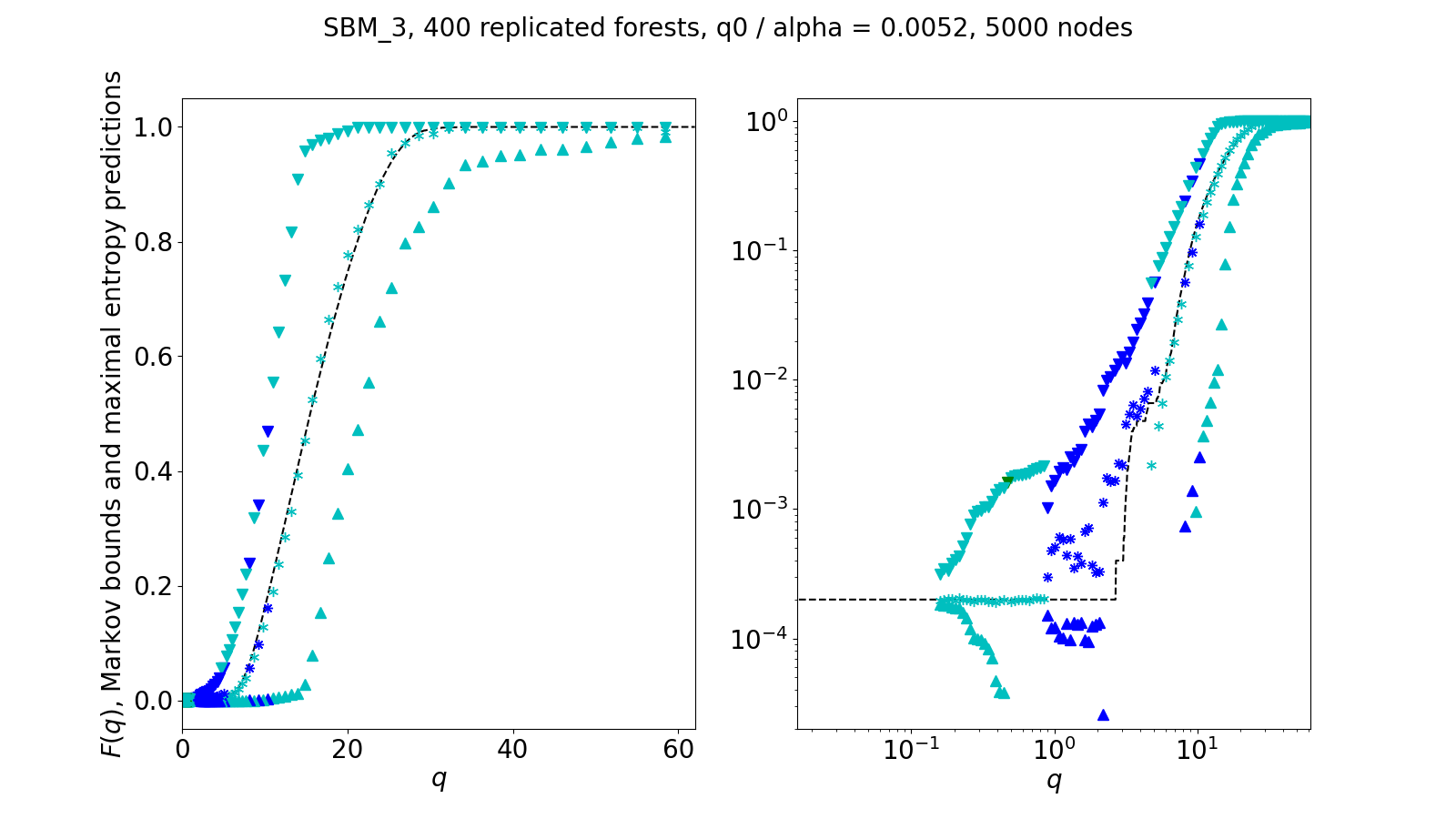}%
		\hfill%
		\includegraphics[clip, trim={.06\wd0} {.03\ht0} {.09\wd0} {.06\ht0}, width=2.8in]{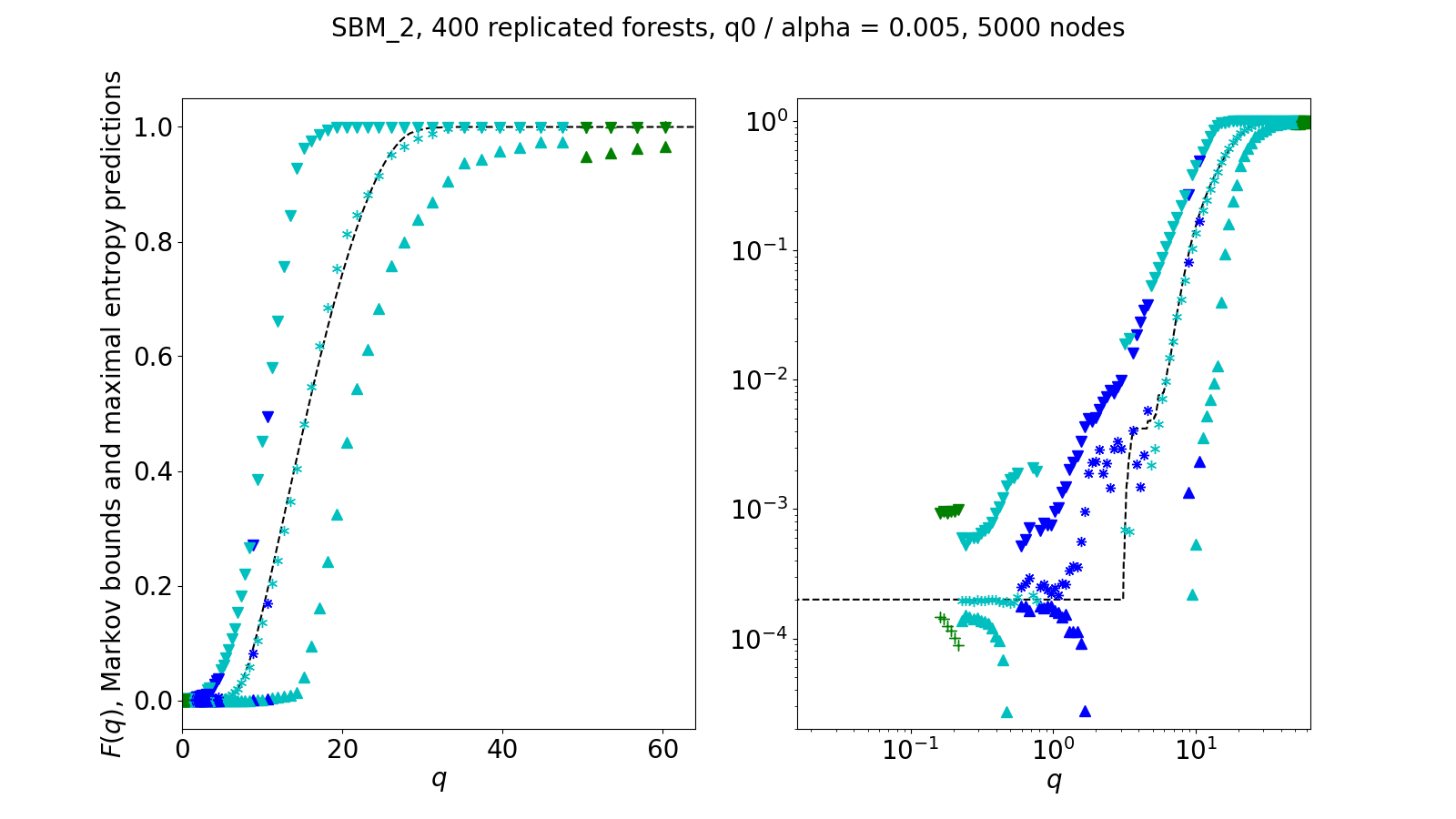}%
	}
\end{figure}

\begin{figure}[tbp]
	\caption{\footnotesize Two sensor graphs.\label{christiane}}
	\sbox0{\includegraphics{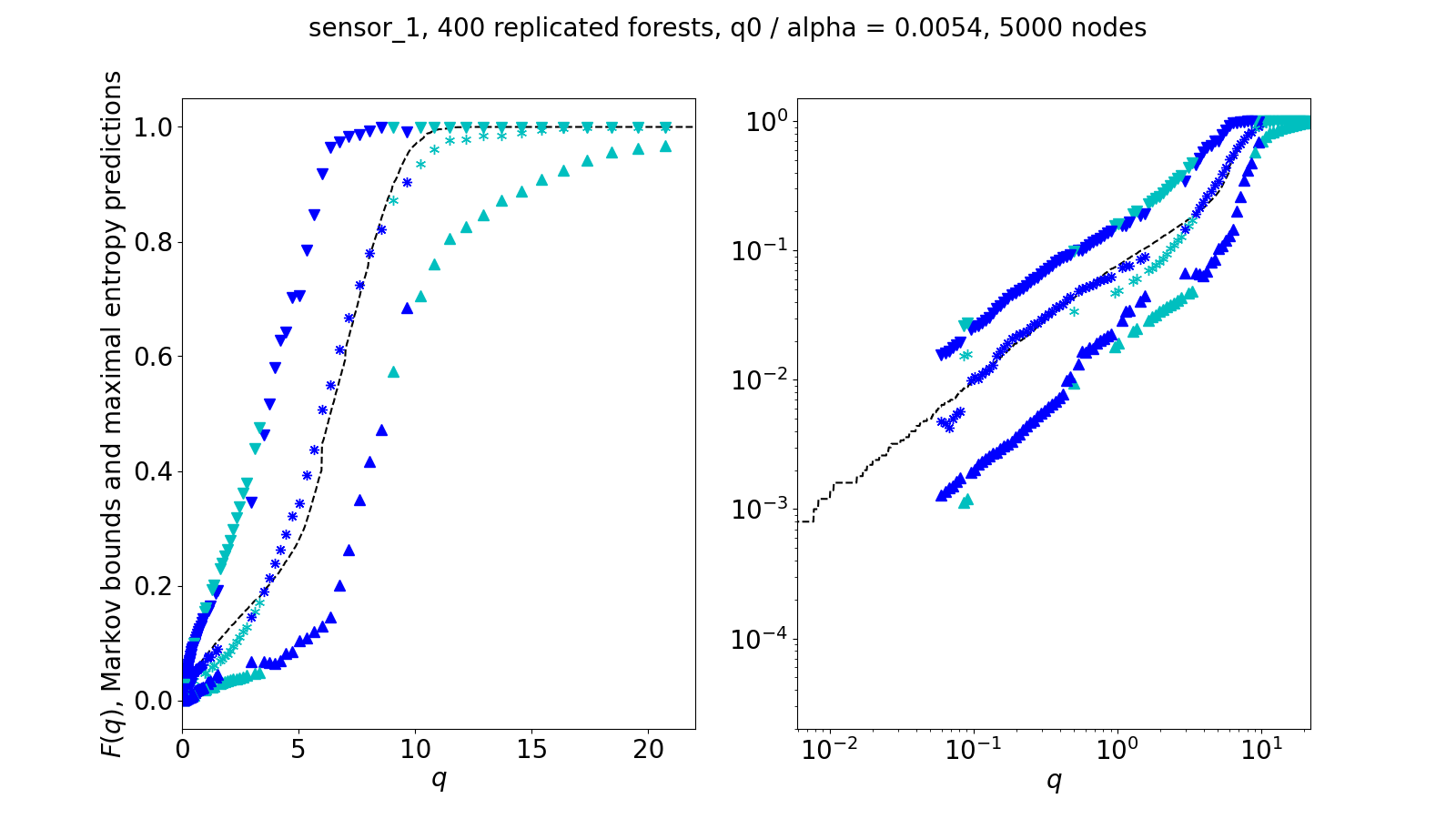}}
	\hbox to \hsize{%
		\includegraphics[clip, trim={.06\wd0} {.03\ht0} {.09\wd0} {.06\ht0}, width=2.8in]{./pictures/sensor_1.png}%
		\hfill%
		\includegraphics[clip, trim={.06\wd0} {.03\ht0} {.09\wd0} {.06\ht0}, width=2.8in]{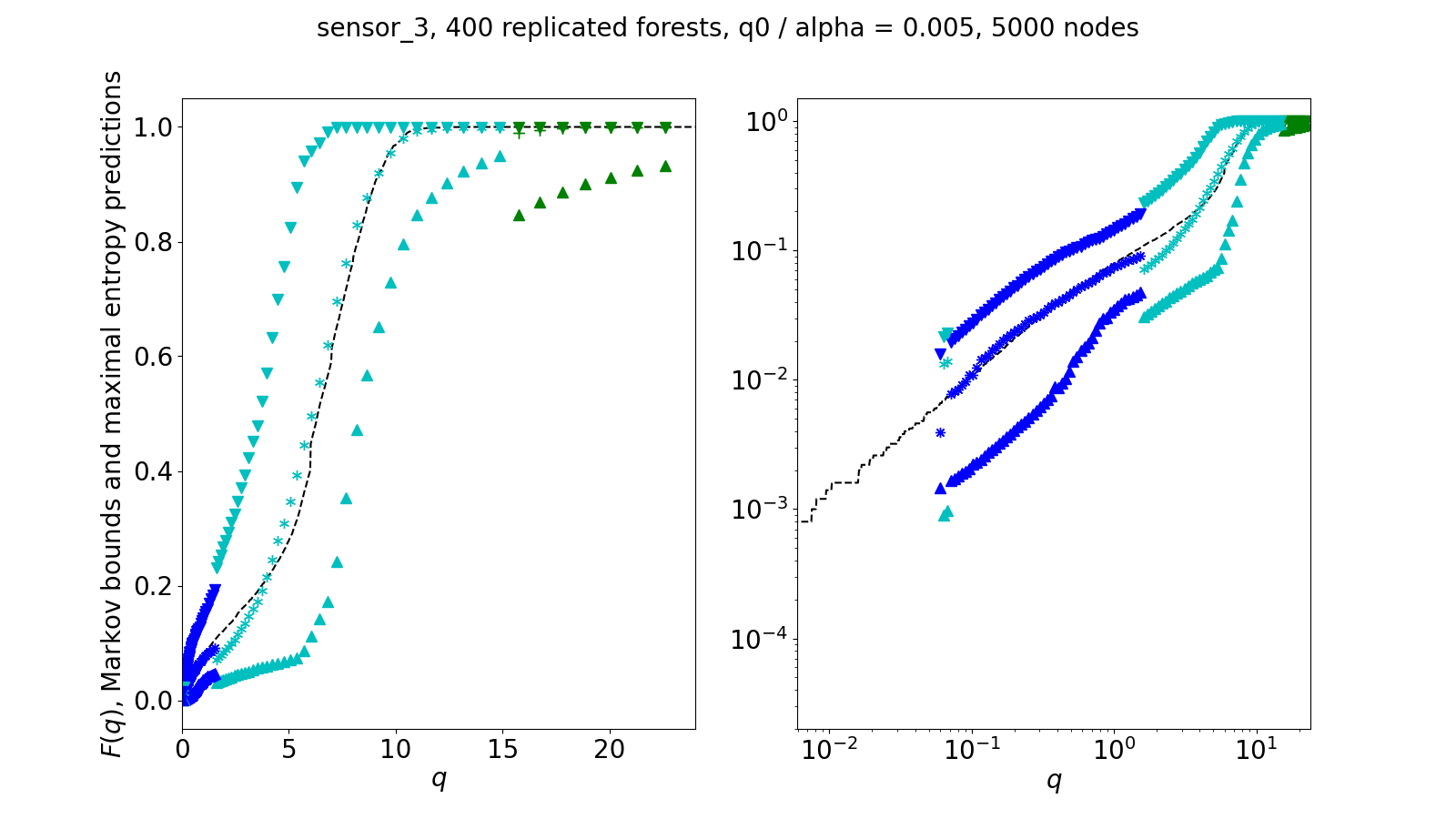}%
	}
\end{figure}

Figure~\ref{albane} shows our results
for two different Barabási-Albert preferential attachment random networks
with $m = 5$ new edges for each new node;
Figure~\ref{celine} shows our results for two different realisations
of a stochastic bloc model with 20 quite well formed communities
of the equal sizes, with intra and inter connection probabilities
equal to .049 and .0008, respectively,
resulting in an average degree equal to 16;
Figure~\ref{christiane} shows our results for two different 
sensor graphs, for which each node in a random set of points
is connected with its 5 nearest neighbours.

\begin{figure}[tbp]
	\caption{\footnotesize Comet and torus graphs.\label{annie}}
	\sbox0{\includegraphics{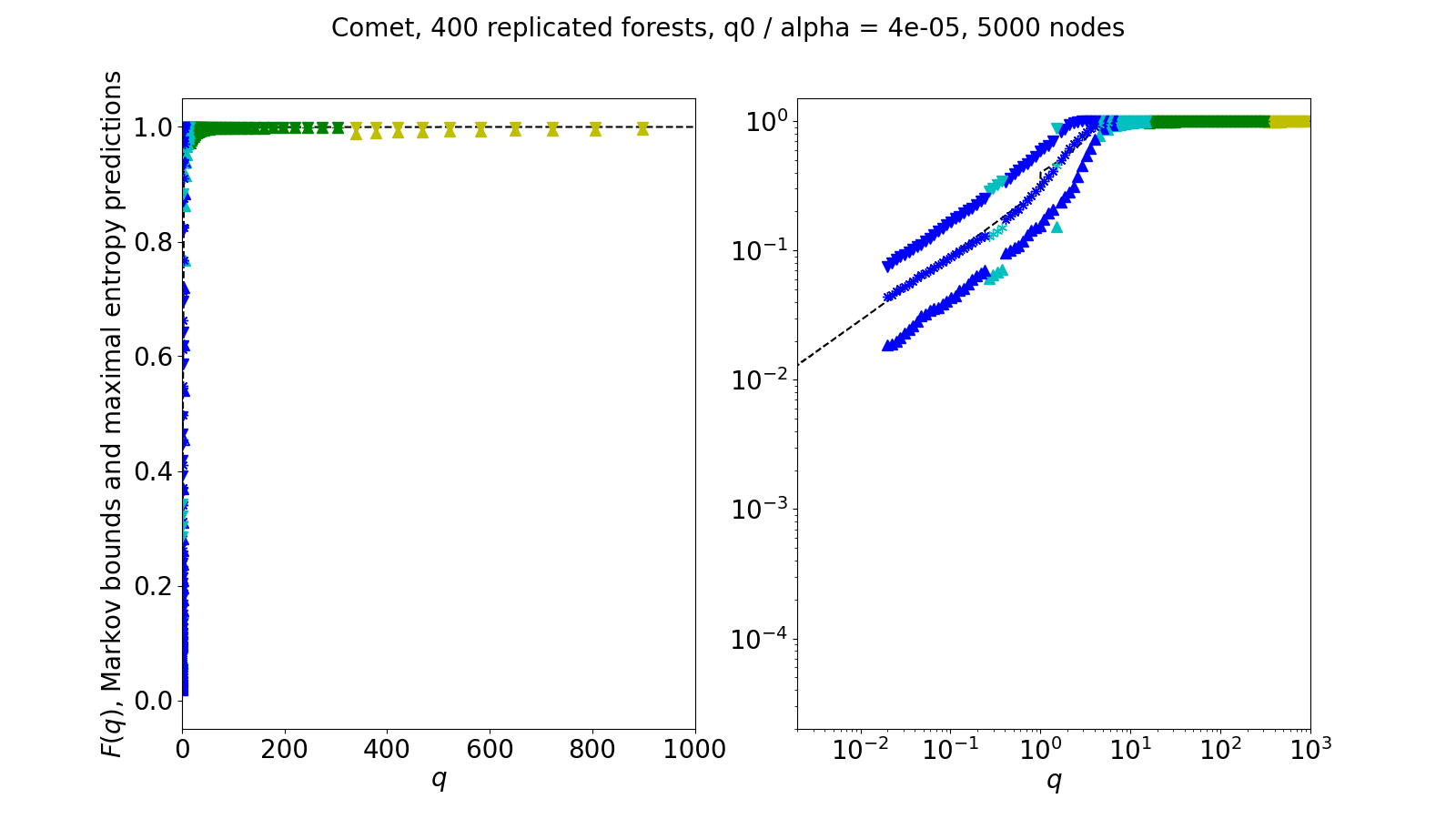}}
	\hbox to \hsize{%
		\includegraphics[clip, trim={.06\wd0} {.03\ht0} {.09\wd0} {.06\ht0}, width=2.8in]{./pictures/Comet.png}%
		\hfill%
		\includegraphics[clip, trim={.06\wd0} {.03\ht0} {.09\wd0} {.06\ht0}, width=2.8in]{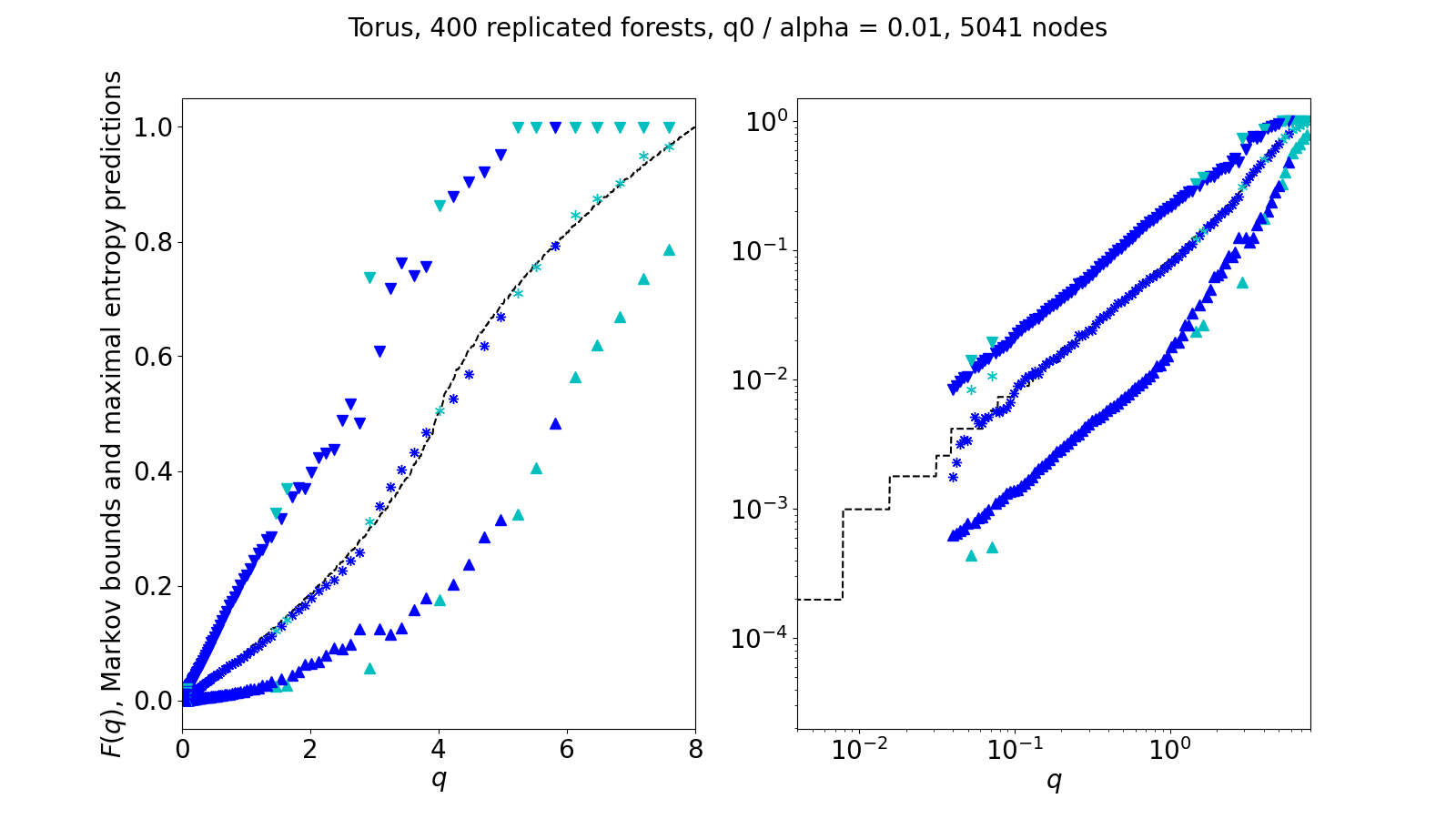}%
	}
\end{figure}

\begin{figure}[tbp]
	\caption{\footnotesize Karate club and Minnesota road networks.\label{antoine}}
	\sbox0{\includegraphics{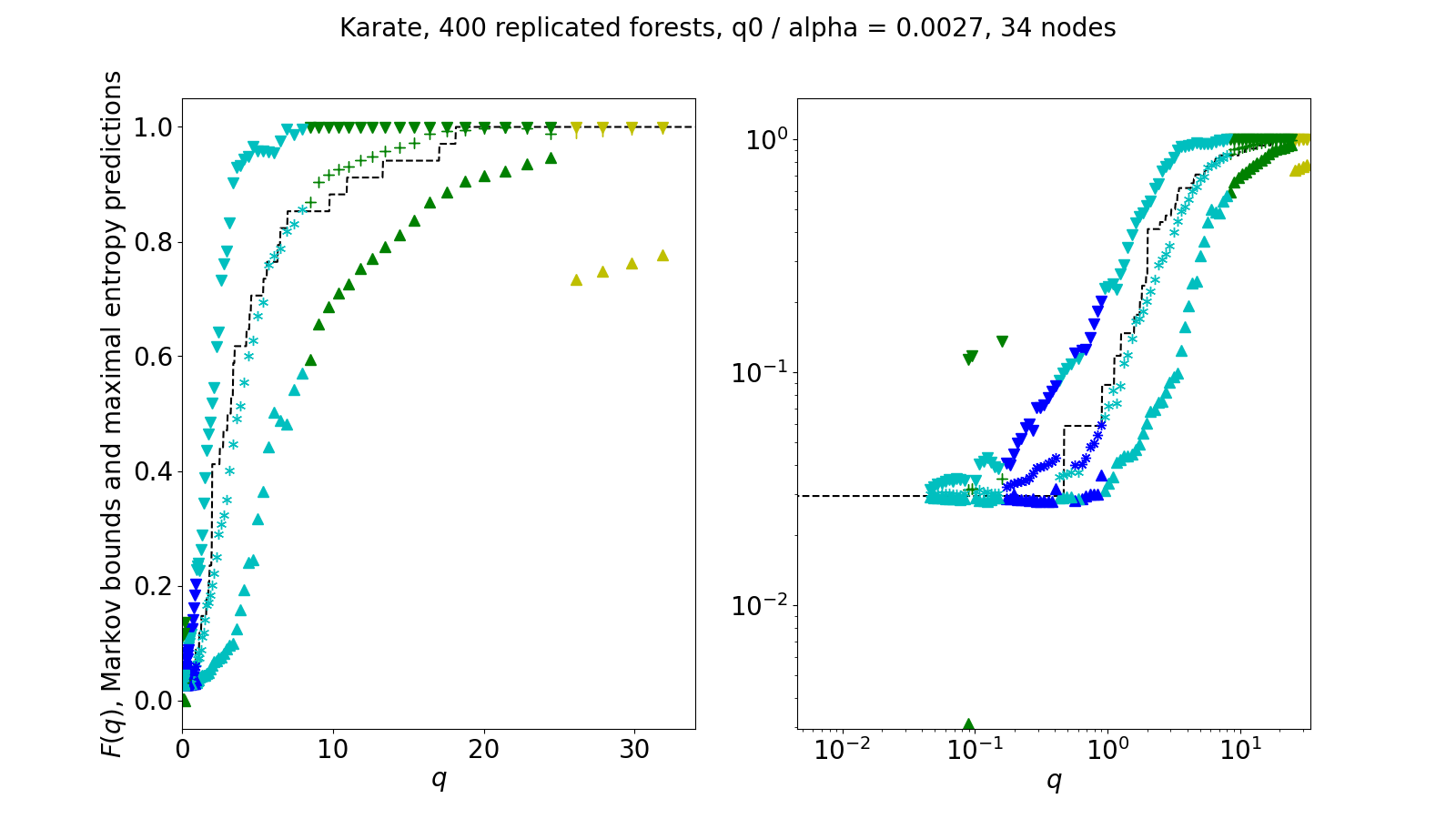}}
	\hbox to \hsize{%
		\includegraphics[clip, trim={.06\wd0} {.03\ht0} {.09\wd0} {.06\ht0}, width=2.8in]{./pictures/Karate.png}%
		\hfill%
		\includegraphics[clip, trim={.06\wd0} {.03\ht0} {.09\wd0} {.06\ht0}, width=2.8in]{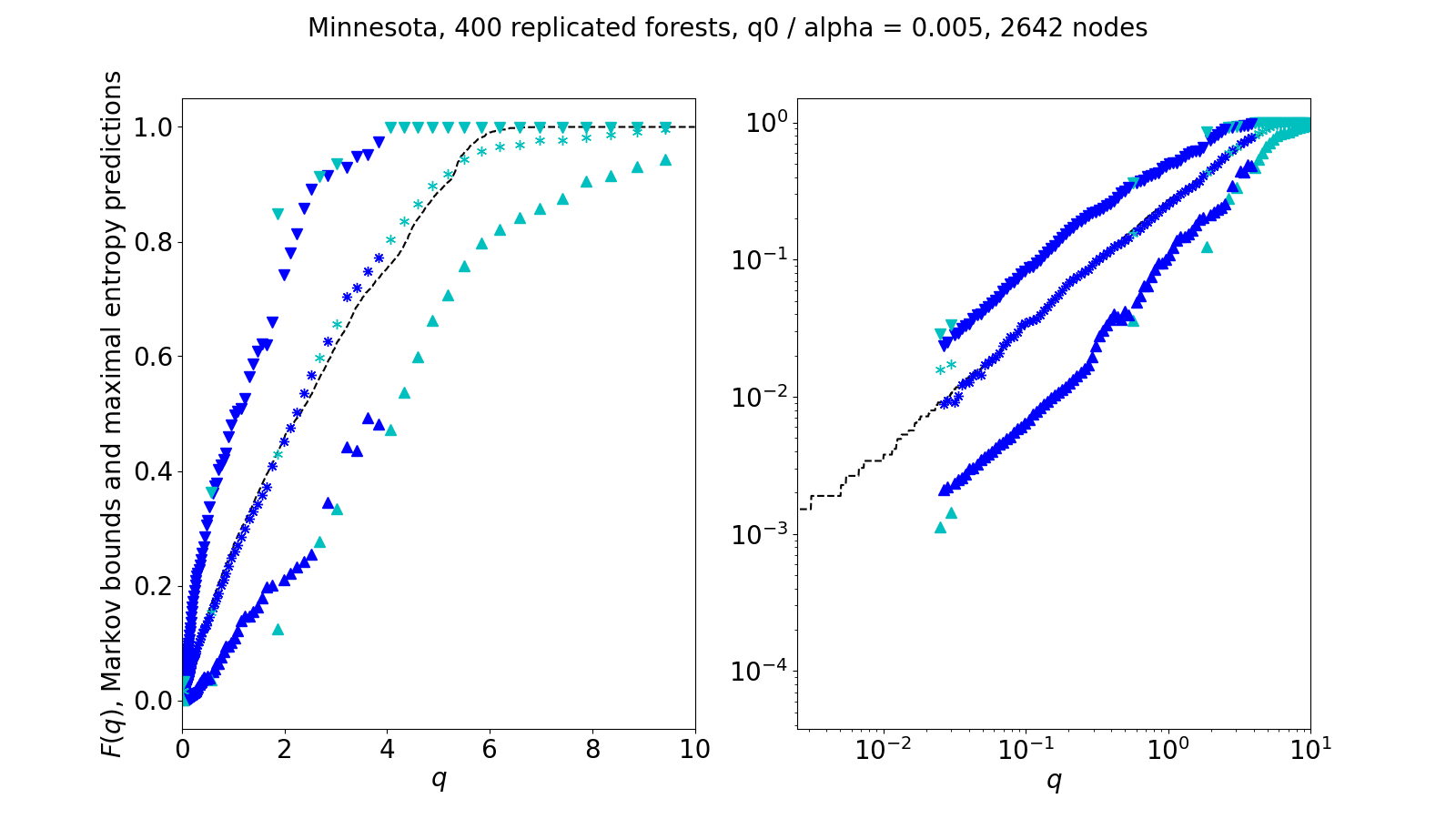}%
	}
\end{figure}

Figure~\ref{annie} shows our results for two deterministic graph.
On the left the comet graph, which is a rooted tree with 499 branches
of length 1 and 1 branch of length 4500.
On the right the $71\times 71$ torus.
Figure~\ref{antoine} shows our results
for the so-called karate club network with 34 vertices and 78 edges on the left,
and for the for the Minnesota road network with 2642 nodes and 3304 links on the right.

\subsubsection{What went wrong}
\label{alain}

\begin{figure}[tbp]
	\caption{\footnotesize 400 and 1600 replicated coupled forests
		for a swiss roll graph.\label{julien}}
	\sbox0{\includegraphics{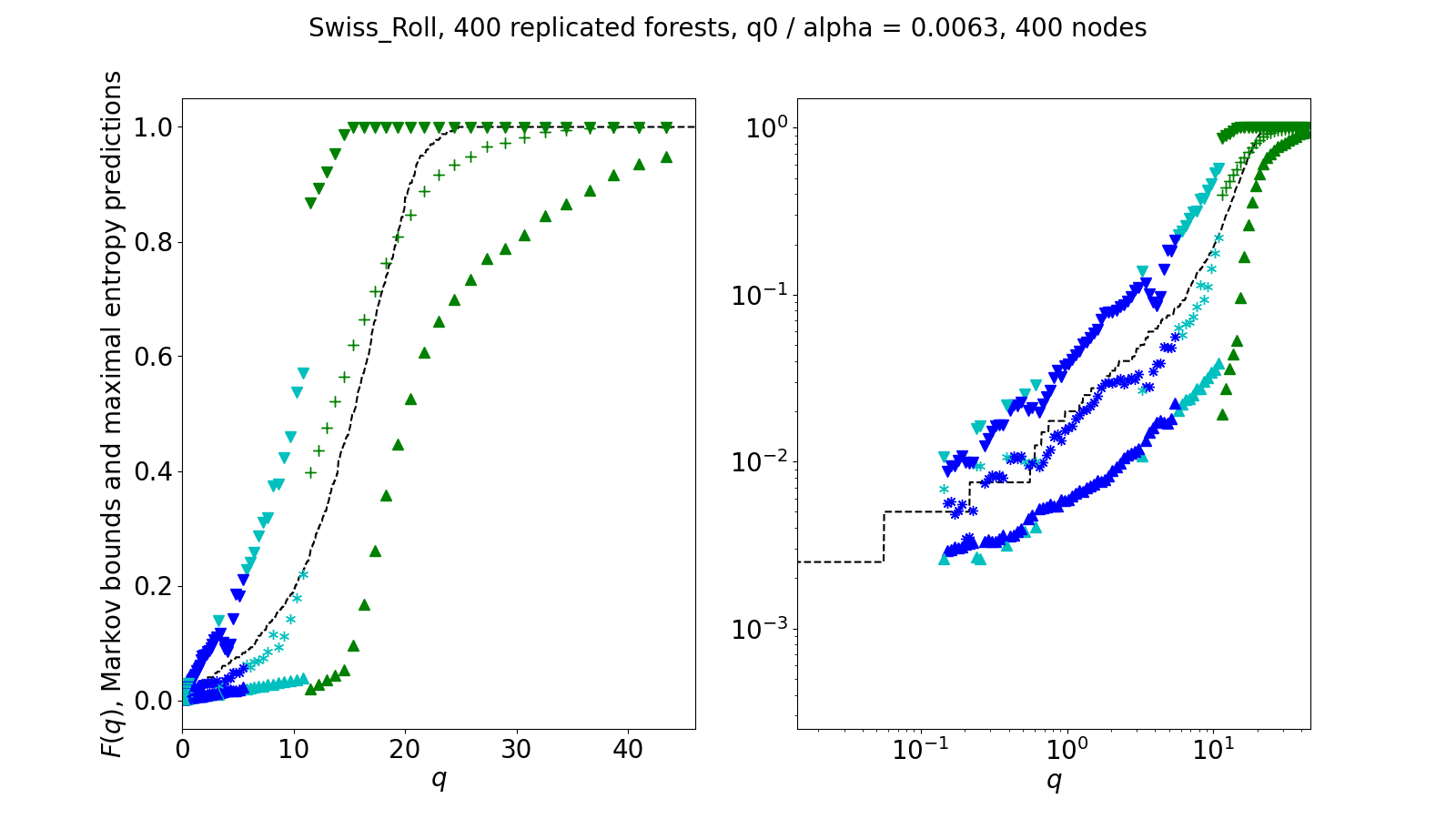}}
	\hbox to \hsize{%
		\includegraphics[clip, trim={.06\wd0} {.03\ht0} {.09\wd0} {.06\ht0}, width=2.8in]{./pictures/Swiss_roll.png}%
		\hfill%
		\includegraphics[clip, trim={.06\wd0} {.03\ht0} {.09\wd0} {.06\ht0}, width=2.8in]{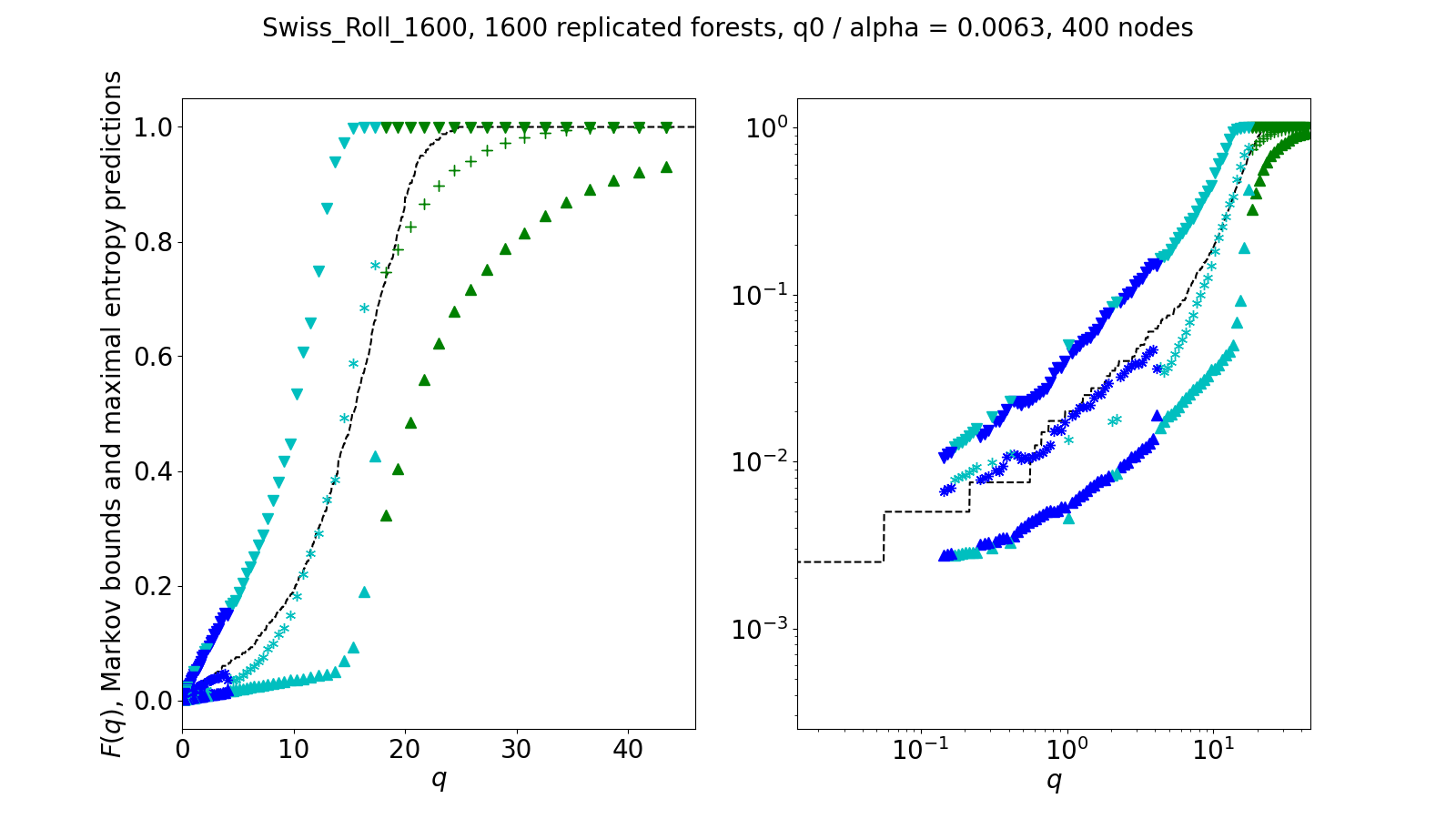}%
	}
\end{figure}

Figures \ref{jean-francois}--\ref{antoine}
already show that reaching a larger number of valid moment estimates
will not always improve the maximal entropy estimator.
It strengthens Markov's bound,
but these give a large window
that goes quite slowly to zero.
Getting more than a few valid moment estimates
is actually out of reach such a Monte Carlo method.
When $l$ increases,
the width $m_{l + 1}^+ - m_{l + 1}^-$ given by Lemma~\ref{frank}
for the interval where $m_{l + 1}$ must lie in goes dramatically to zero,
indeed.
Hence, increasing the forest number
when our microcanonical model gives a wrong prediction
with a few moments does not help.
Figure~\ref{julien} illustrates this
for the upper part of the spectrum of $-L$.
It shows our results with 400 replicated coupled forest on the left
and 1600 replicated coupled forest on the right for a ``swiss roll'' graph
with 400 vertices and 2883 edges.

\begin{figure}[tbp]
	\caption{\footnotesize 400 and 1600 replicated coupled forests for the star graph.\label{rachel}}
	\sbox0{\includegraphics{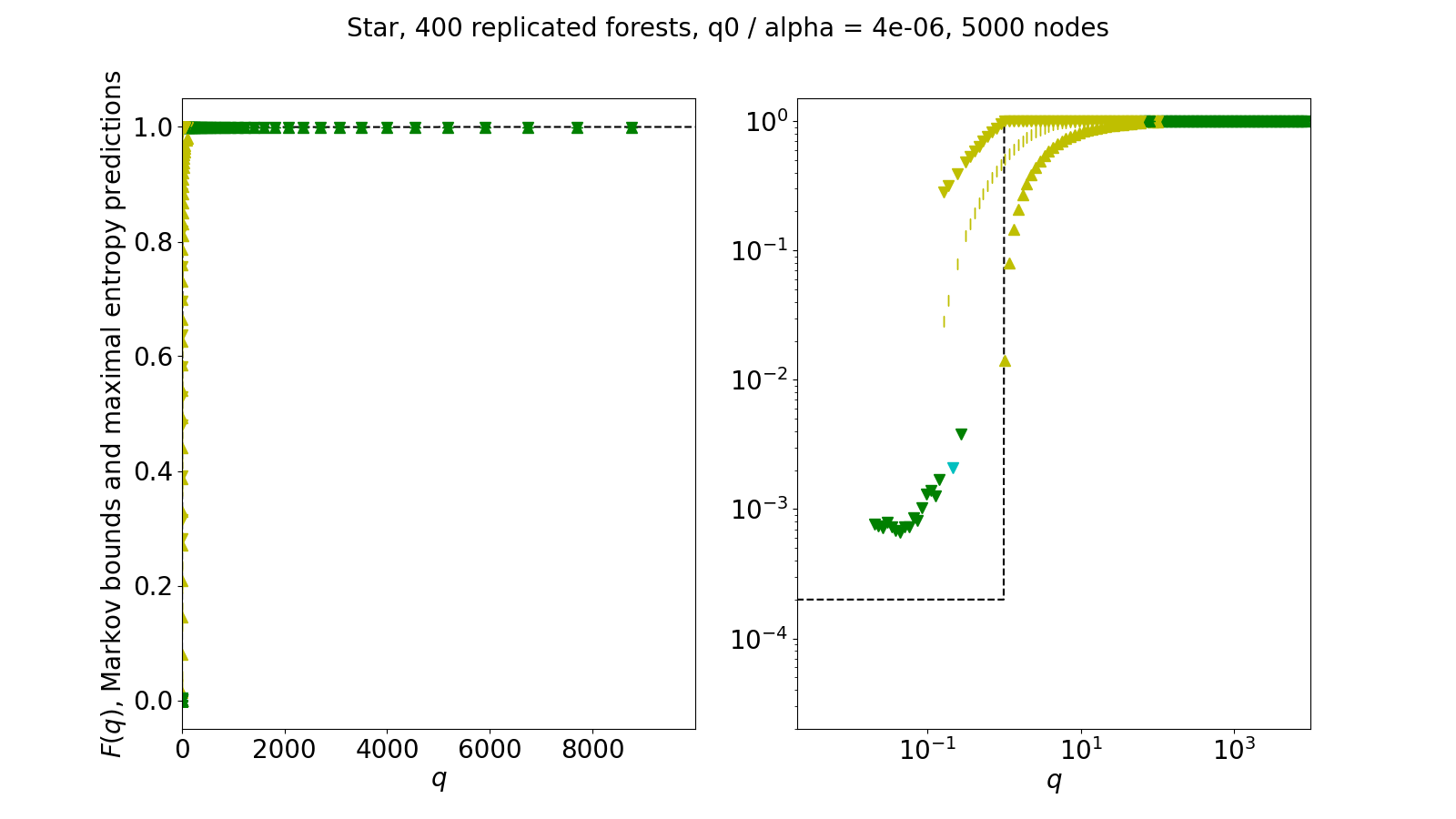}}
	\hbox to \hsize{%
		\includegraphics[clip, trim={.06\wd0} {.03\ht0} {.09\wd0} {.06\ht0}, width=2.8in]{./pictures/Star.png}%
		\hfill%
		\includegraphics[clip, trim={.06\wd0} {.03\ht0} {.09\wd0} {.06\ht0}, width=2.8in]{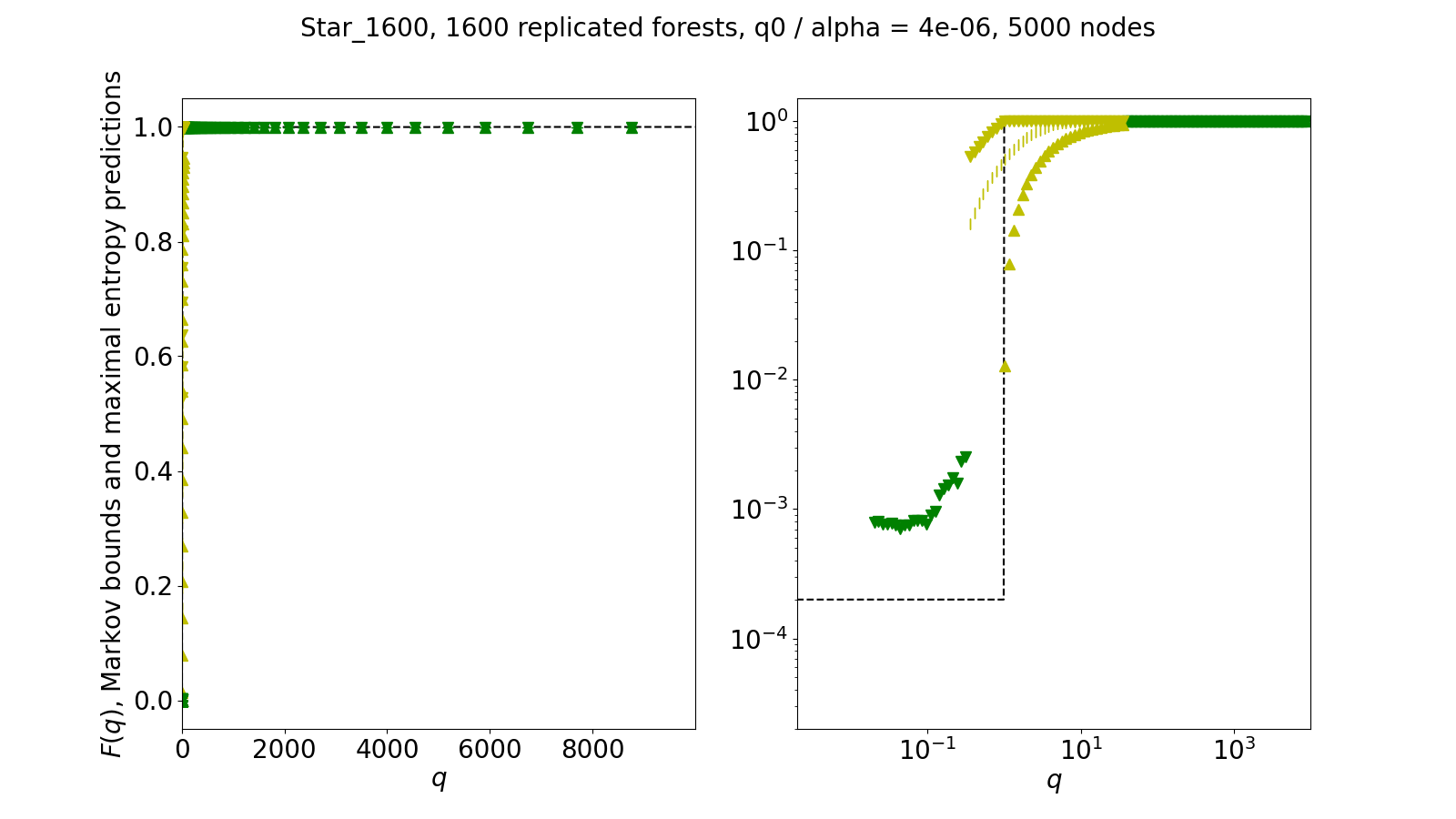}%
	}
\end{figure}

This problem can also occur when only one moment is available,
which is the case for the star graph, rooted tree with $n - 1 = 4999$ branches
of length 1.
As shown by Figure~\ref{rachel},
we cannot rely on the microcanonical ensemble with only one moment
and simply ignore higher (non-admissible) moment estimates.

\begin{figure}[tbp]
	\caption{\footnotesize 100 replicated coupled forests
		for our two denser Erdős-Rényi graphs.\label{raphael}}
	\sbox0{\includegraphics{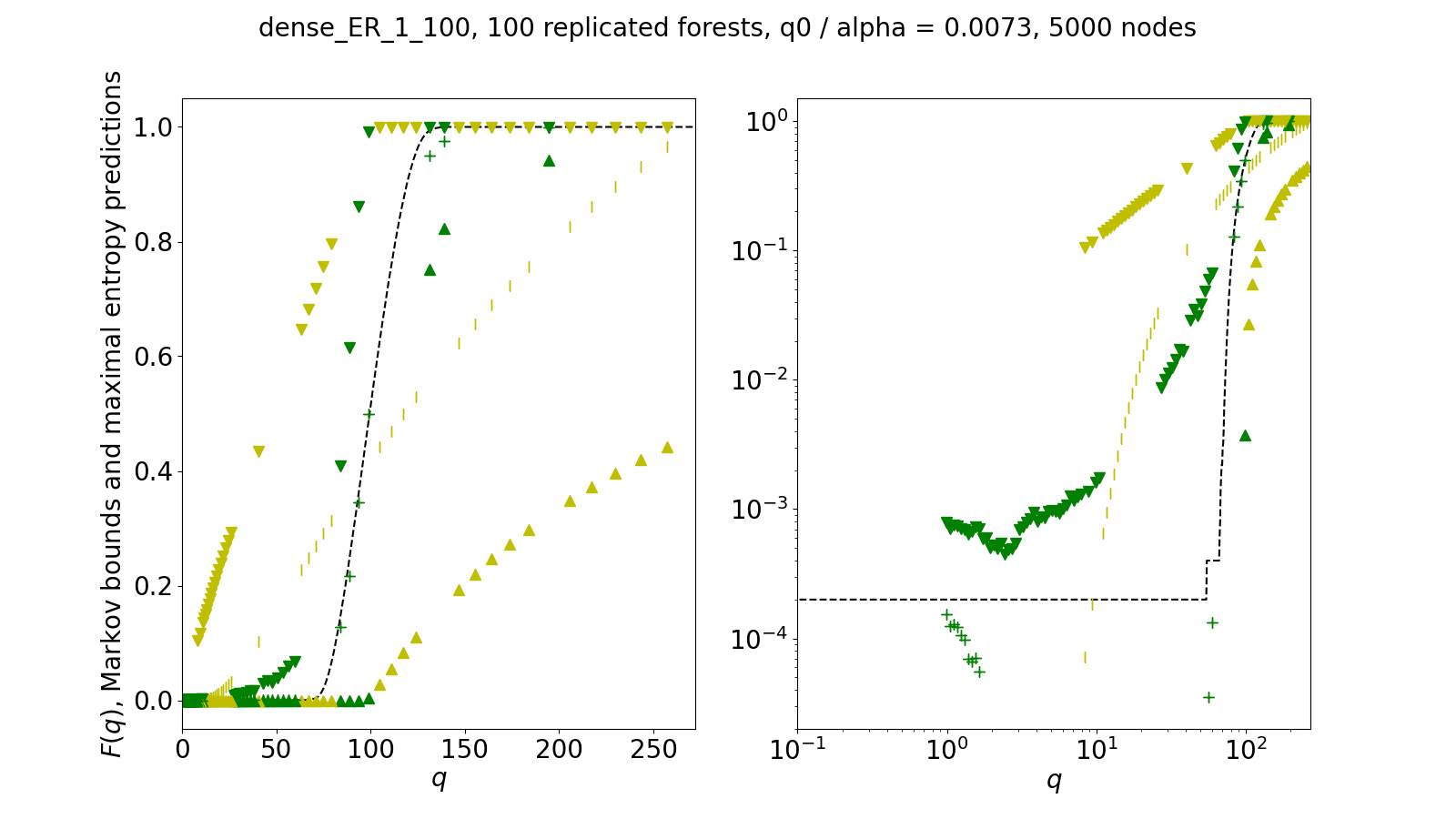}}
	\hbox to \hsize{%
		\includegraphics[clip, trim={.06\wd0} {.03\ht0} {.09\wd0} {.06\ht0}, width=2.8in]{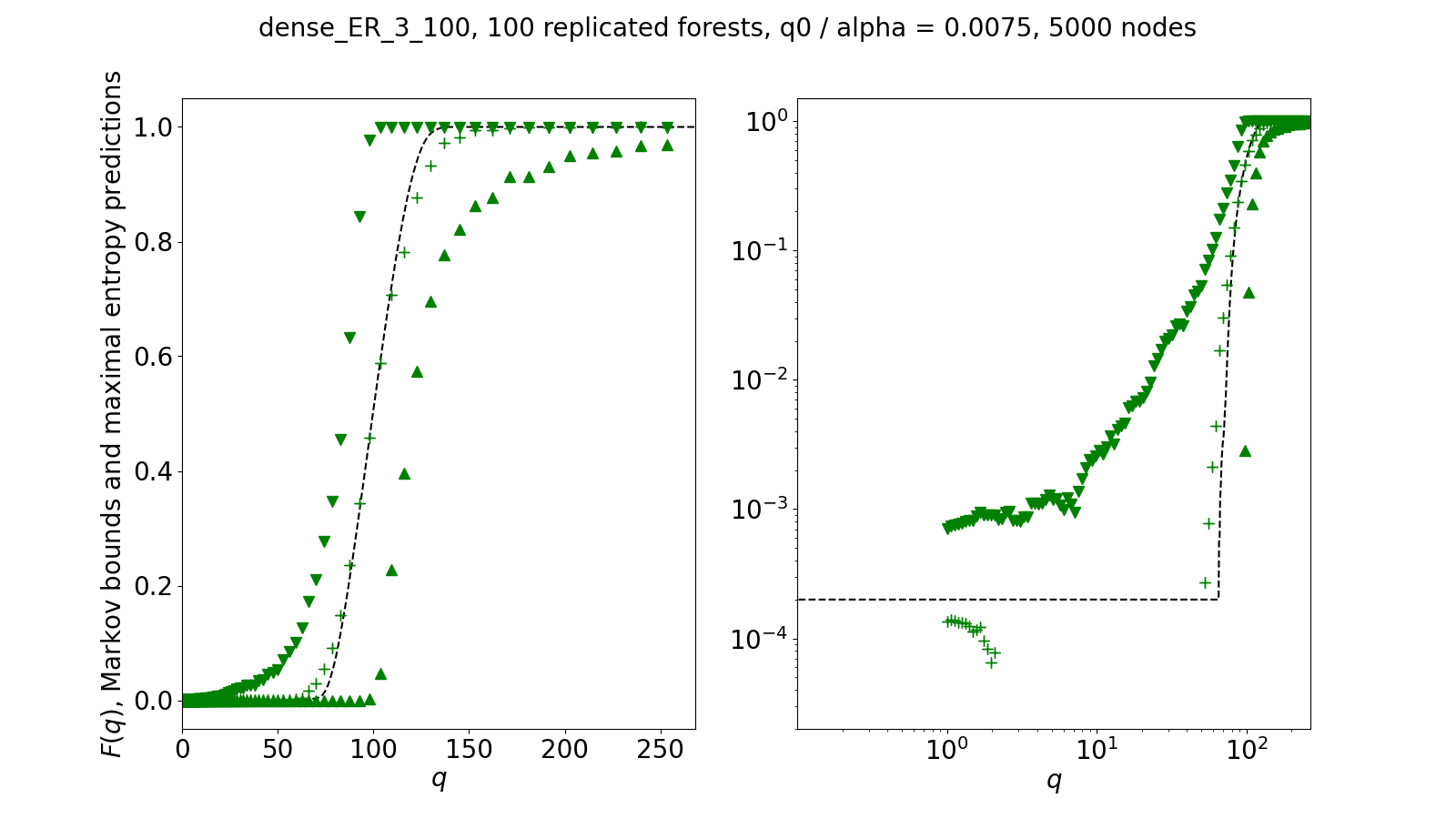}%
		\hfill%
		\includegraphics[clip, trim={.06\wd0} {.03\ht0} {.09\wd0} {.06\ht0}, width=2.8in]{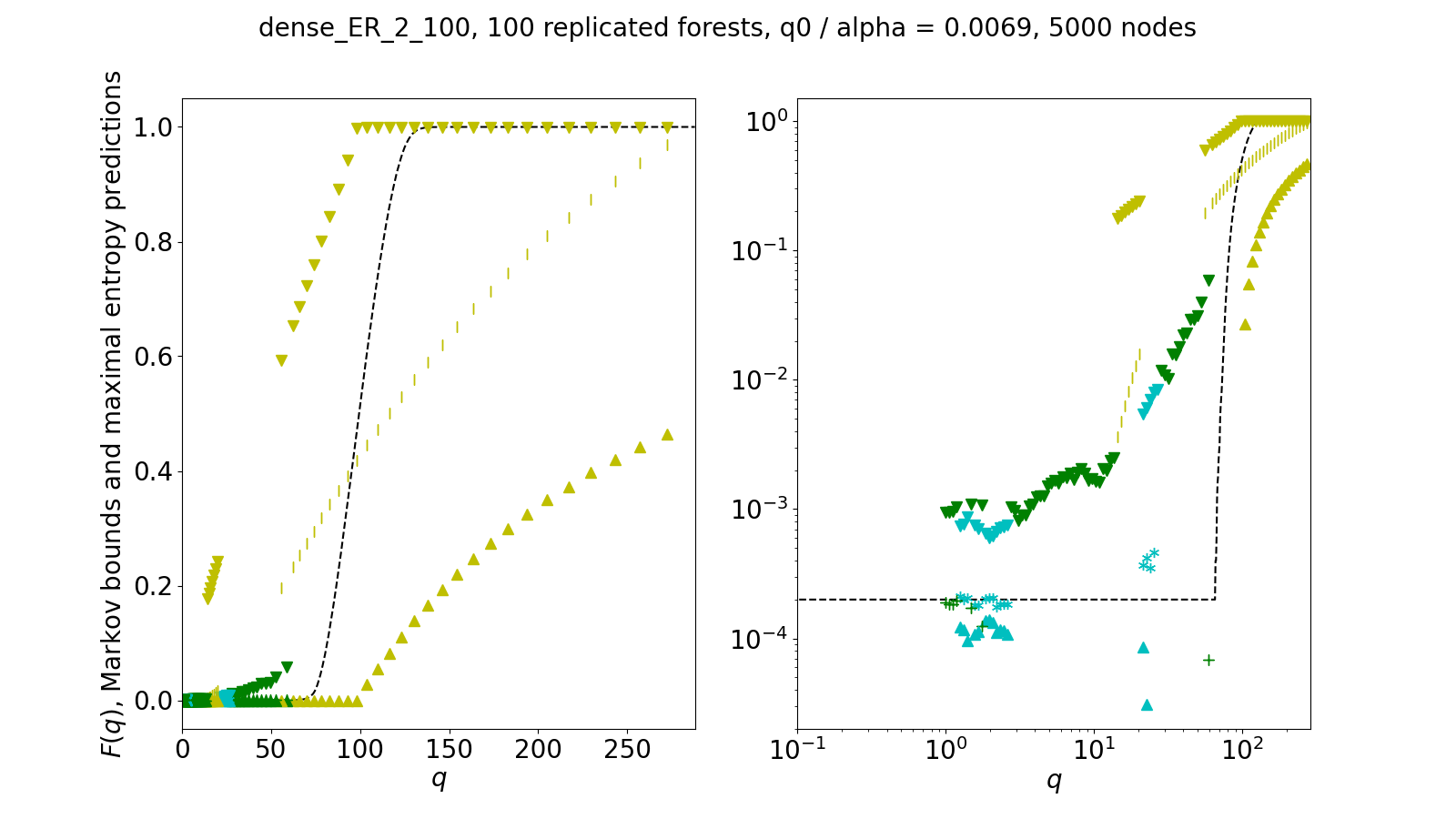}%
	}
\end{figure}

This is an extreme case and we can usually go beyond one valid moment estimate by sampling more forests.
However, when reaching one more valid moment estimate, our maximum entropy estimator changes drastically
in this case. This is a source of instability as illustrated by Figure~\ref{raphael},
which shows our results with only 100 replicated forests for the previous two denser Erdős-Rényi graphs:
in one case we got a second valid moment estimate for the relevant values of $q$,
in one case we did not.

\subsection{Room for improvement}

\subsubsection{Monte Carlo measures}

We decided to make moment estimates and compute Markov bounds 
as well as maximal entropy estimators in only $1 / \epsilon_0$ values 
of $q$  because each measurement of $|\xi^k_q|$, $k \leq l$,
for a given value of $q$ has a computational cost in $O(n)$.
But we could instead measure the difference
$|\xi^k_{q_{i + 1}}| - |\xi^k_{q_i}|$
between two successive values of $q$ along the coupled forest algorithm.
Since each forest update is only local, we might be able to compute
this difference at a lower numerical cost.
We could then perform spectral estimation on a finer scale. 

We could also rely on other observables along coupled forest trajectories,
still related with the Laplacian spectrum.
For example, by adding self loops into the arrow stacks
to get space-homogeneous unfreezing rates $1 / (1 + \alpha t)$, $t \geq 0$,
instead of $1 / (1 + w(x) t)$, $x \in \cl X$,
the total number of arrows one has to read to update the current forest
at an unfreezing time $t$ is a mixture of geometric random variables
with $n$ mixture weights and success probabilities,
each of them depending on a single eigenvalue $\lambda_j$ of $-L$.
Working with a mixture rather than a sum of random variables
---like $|\rho(\Phi_q)|$,
which is a sum of $0-1$ Bernoulli random variables
with mean $q / (q + \lambda_j)$, $j < n$---
can give better insights into the spectrum.
Those geometric laws could also be related
with better scale separations than these Bernoulli laws.

We finally note that in our experiments
we did not take advantage of the fact
that each permutation of the $l$ forests 
needed to build $\xi^k_q$, $k \leq l$,
for a given value of $q$ 
lead to a different set
with the same mean size.
This could be used to implement
some variance reduction
in our Monte Carlo measurements.

\subsubsection{Using more information for Stieltjes transform inversion}

Our last observation in Section~\ref{alain}
strongly suggests that we should not simply 
get rid of higher moment estimates
$\hat m_{k + 1}(q)$, $\hat m_{k + 2}(q)$,~\dots\
when we only have $k$ valid moment estimates
at a given $q$ in $[q_0, 2\alpha]$.
In the comparison with state of the art spectrum estimation methods
that we presented in~\cite{BCGQT},
we already included some projection 
on the set of admissible moment sequences.
It is however clear
that in order to estimate $F(q)$
there is still much more information to be used:
that contained in $\hat m_k(q')$, $k \leq l$, for $q' \neq q$.
How to use it is outside the scope of this paper.
Here we simply wanted to point out
that our low-cost pointwise inversion,
which relates the law of $\Phi_q$ for a given $q$
with the number of eigenvalues below $q$,
already provides a sound estimation of this number.
Together with the coupled forest process
that covers all values of $q \geq \epsilon_0\bar\lambda$
at a sampling cost in $O(\epsilon_0^{-1}n\ln n)$,
this seems to us quite promising
for spectrum estimation in large networks.
By using the double cover trick
that we recall in Appendix~\ref{olga},
we can extend this perspective
to spectrum estimation of all real symmetric matrices.\\

\noindent\textbf{Acknowledgments}. This work was partially funded by the Région Auvergne-Rhône-Alpes project
TrustGNN, the ANR project GRANOLA (ANR-21-CE48-0009) and the European Union (Marie Curie
project Rand4TrustPool, 101148828)

\bibliographystyle{plainnat} 

\bibliography{biblio.bib}   

\newpage

\begin{appendix}
	\section{Working with any symmetrical real matrix}
	\label{olga}
	
	A symmetric real matrix $M = (M(x, y))_{x, y \in \cl X}$
	has to satisfy two constraints for being interpreted
	as the Laplacian $L$ of a graph with vertex set $\cl X$.
	First, its off-diagonal entries must be non-negative,
	\begin{equation}\label{eric}
		M(x, y) \geq 0, 
		\qquad x \neq y,
	\end{equation}
	second, each diagonal entry on a line $x$ in $\cl X$
	must be the opposite of the sum all off-diagonal entries
	on this line,
	$$
	M(x, x) = -\sum_{y \neq x} M(x, y),
	\qquad x \in \cl X.
	$$
	As described in Section~\ref{anne},
	our spectrum estimation method immediately extend
	to sub-Laplacians $L = M$, which satisfy~\eqref{eric}
	but only 
	$$
	M(x, x) \leq - \sum_{y \neq x} M(x, y),
	\qquad x \in \cl X.
	$$
	In Section~\ref{georges} we recall the double cover trick,
	which appears in \cite{Gr}.
	It reduces to the sub-Laplacian case
	that of any diagonally dominant symmetric matrix,
	which satisfies
	\begin{equation}\label{sebastien}
		M(x, x) \leq - \sum_{y \neq x} |M(x, y)|,
		\qquad x \in \cl X.
	\end{equation}
	Since by subtracting a large enough multiple of the identity
	one can always reduce the case of any real symmetric matrix
	to this diagonally dominant case with a shifted spectrum,
	our spectrum estimation method extends to all real symmetric matrices.
	
	\subsection{Sub-Markovian generators}
	\label{anne}
	
	In the case of a sub-Laplacian $L = M$ we only have to consider
	loop-erased random walks killed in each $x$ in $\cl X$
	at rate $q + \delta(x)$ instead of $q$, with
	$$
	\delta(x) = w(x) - \sum_{y \neq x} w(x, y) \geq 0,
	\qquad x \in \cl X,
	$$
	where $w(x) = -M(x, x)$ and $w(x, y) = M(x, y)$ for all $x \neq y$.
	As far as stacks of marks and arrows are concerned,
	this amounts to sample at each level associated with a vertex $x$
	a uniform mark as previously and an arrow $(x, y)$, with $x \neq y$,
	with the same probability $w(x, y) / w(x)$
	or a ``killing arrow'' with the remaining probability $\delta(x) / w(x)$.
	To build the coupled forest process
	$(\Phi_q)_{q > 0}$, we just have to replace each arrow by a stop
	whenever the associated mark is below $q / (q + w(x))$.
	We then have two kinds of roots in each $\Phi_q$,
	those associated with such a low mark,
	which can be unfrozen,
	and those associated with a killing arrow,
	which cannot.
	For replicated forests $\Phi_{k, q}$, $k < l$,
	we denote by $\tilde\rho(\Phi_{k, q})$ 
	the set of roots of the first kind in $\Phi_{k, q}$,
	we define (recall~\eqref{felix})
	$$
	\tilde\xi_q^k = \bigl\{x \in \cl X :
	R^1(x) \in \tilde\rho(\Phi_{0, q}),\, \dots,\,
	R^k(x) \in \tilde\rho(\Phi_{k - 1, q}),\,
	R^k(x) = x\bigr\}
	$$
	we get
	$$
	\ds E\bigm[|\tilde\xi_q^k|\bigr] = \sum_{j < n} \left(q \over q + \lambda_j\right)^k
	$$
	where the $\lambda_j$, $j < n$, are the eigenvalues of $-L$,
	and the whole analysis can be repeated.
	
	\subsection{The double cover trick}
	\label{georges}
	
	For any diagonally dominant real symmetric matrix
	$M = (M(x,y))_{x, y \in \cl X}$ satisfying~\eqref{sebastien},
	we call $\cl X_1 = \cl X$ a vertex set, and
	we associate with each vertex $x$ in $\cl X$ another vertex $-x$
	in a one-to-one correspondence with another and disjoint vertex set $-\cl X$.
	We define $\cl X_2 = \cl X \cup -\cl X$ 
	in order to build two sub-Laplacians
	$L_1 = (L_1(x, y))_{x, y \in \cl X_1}$
	and $L_2 = (L_2(x, y))_{x, y \in \cl X_2}$
	such that
	\begin{equation}\label{adore}
		2\sigma_2 = \sigma_1 + \sigma
	\end{equation}
	with $\sigma_2$, $\sigma_1$ and $\sigma$
	the spectral measures\footnote{The factor 2 accounts
		for the normalisations of spectral measures.}
	associated with $L_2$, $L_1$ and $M$ respectively.
	With $F_2$, $F_1$ and $F$ the associated cumulative distribution functions,
	it will then hold
	$$
	F = 2 F_2 - F_1,
	$$
	so that an estimates of $F_1$ and $F_2$ provides an estimate on $F$.
	
	To this end we set for all {\it distinct} $x$ and $y$ in $\cl X$
	$$
	L_2(x, y) = L_2(-x, -y) = [M(x, y)]_+ \geq 0,
	$$
	$$
	L_2(x, -y) = L_2(-x, y) = [M(x, y)]_- \geq 0,
	\qquad L_2(x, -x) = L_2(-x, x) = 0
	$$
	and 
	$$
	L_2(x, x) = L_2(-x, -x) = M(x, x)
	\leq -\sum_{z \neq x} [M(x, z)]_+ + [M(x, z)]_-
	$$
	so that $L_2$ is a sub-Laplacian, indeed.
	We then defined $L_1$ as the sub-Markovian generator
	that describes the projection on $\cl X = \cl X_1$
	of the killed random walk described by $L_2$.
	For distinct $x$ and $y$ in $\cl X$, 
	it holds
	$$
	L_1(x, y) = L_2(x, y) + L_2(x, -y) = L_2(-x, y) + L_2(-x, -y) = |M(x, y)|
	$$
	and 
	$$
	-L_1(x, x) = -L_2(x, x) = -L_2(-x, -x) = -M(x, x) \geq \sum_{z \neq x} L_1(x, z).
	$$
	
	To prove~\eqref{adore} 
	it suffices then to observe
	that the orthogonal vector spaces of even and odd functions $f$ in $\ds R^{\cl X_2}$
	---such that, for all $x \in \cl X$, $f(-x)$ is equal to $f(x)$ and $-f(x)$, respectively---
	are stable by $L_2$, which acts on them as $L_1$ and $M$, respectively.
	$L_2$ has indeed a block matrix representation
	$$
	L_2
	= \begin{pmatrix}
		M_+ &  M_- \\
		M_- &  M_+ 
	\end{pmatrix}
	$$
	such that for all row vector $U$ it holds
	$$
	L_2 \begin{pmatrix}
		U \\
		U
	\end{pmatrix}
	= \begin{pmatrix}
		M_+U + M_-U \\
		M_-U + M_+U
	\end{pmatrix}
	= \begin{pmatrix}
		L_1 U\\
		L_1 U 
	\end{pmatrix},
	\quad
	L_2 \begin{pmatrix}
		U \\
		-U
	\end{pmatrix}
	= \begin{pmatrix}
		M_+U - M_-U \\
		M_-U - M_+U
	\end{pmatrix}
	= \begin{pmatrix}
		MU\\
		-MU 
	\end{pmatrix}.
	$$
	
	\section{Pseudocodes}
	\label{sacha}
	
	\nt{Algorithm~\ref{alg:RSF} is simply Wilson's algorithm \cite{Wi} adapted to the forest case, and with an added priority queue to store the values of $q$ at which nodes need to be reactivated when going down the coupled forest trajectory. 
	Algorithm~\ref{alg:coupled_RSF} is the coupled forest algorithm, that outputs one vector $\texttt{Next}$ (storing the state of the forest) at each value of $q\in[q_{\text{min}},q_{\text{max}}]$ for which it changes. In practice, we never run this algorithm as is as this list of forests can be unnecessary large. We run $r$ (the number of moment estimates we need) coupled forests in parallel, stall the trajectories at the values of $q$ in the grid we chose, compute the number of roots of roots to obtain the moment estimates we need, and reactivate all forests until they stop at the next value of $q$ in the grid, \textit{etc}.}
		
	\nt{Also, we give here an implementation with a priority queue. In the main text, we discussed the computation time for an alternative implementation that consists in adding self-loops until all nodes have the same degree, and running Wilson on this regularized graph. On the positive side, one does not need to store a priority queue as, at the reactivation step, the value of $q$ and its associated node $i$ can be randomly drawn from known distributions (a Beta distribution for $q$ and a uniform distribution for $i$). On the downside, the added self-loops create possibly many extra steps in the random walk. As a rule of thumb, in cases where the maximum degree $\alpha$ is much larger than the average degree $\bar\lambda$ (cases where many self-loops need to be added to regularize the graph), then the priority queue version should be preferred, as its computation time will be even better than the one discussed in Theorem \ref{clo}. We do not delve more into details here: precise implementations strategies and actual time comparisons will be the object of a later work.}
	
	\begin{algorithm}
		\caption{$\left(\text{Next}, \mathcal{P}\right)$ = \texttt{RandomForest}($\mathcal{G},q$)}
		\label{alg:RSF}
		\begin{algorithmic}[1]
			\Inputs{$\mathcal{G}=(\cl X,w)$\\
				$q>0$}
			\Initialize{
				\textit{\# Initially, the forest is empty}\\	
				$\forall i\in \cl X,\quad\text{InForest}[i] \leftarrow \texttt{false} $ 
				\\ $\forall i\in\cl X, \quad\text{Next[i]} \leftarrow \bm{\text{nil}} $ 
				\\ $\forall i\in\cl X, \quad w[i] \leftarrow \sum_{j} w(i,j)$\textit{ \# Weighted degrees} \\
				$\mathcal{P}\leftarrow\emptyset$ \textit{\# Priority queue}
			} 
			\For{$i \in\cl X$}            
			\State{$u \leftarrow i$}
			\State{\textit{\# Start a random walk to create a forest branch}}
			\While{\textbf{not} $\text{InForest}[u]$} \textit{ \# Stop if $u$ is in the forest}
			\State{$\mathcal{U} \leftarrow Unif([0,1])$}
			\If{$\mathcal{U} \leq \frac{q}{q+ w[u] }   $} \textit{\#If true, $u$ becomes a root}
			\State{$\text{InForest}[u] \leftarrow \texttt{true}$ \textit{ \# Add $u$ to the forest}} 
			\State{$\text{Next}[u] \leftarrow 0$ \textit{ \# Set next of $u$ to null}}
			\State{Add $\left(q'=\frac{\mathcal{U} w[u]}{1-\mathcal{U}},u\right)$ to the priority queue $\mathcal{P}$}
			\Else \textit{ \# If false, continue the random walk}
			
			\State{$\text{Next[u]} \leftarrow \texttt{RandomSuccessor}(u,\mathcal{G})$}
			\State{$u\leftarrow\text{Next}[u]$}
			\EndIf
			\EndWhile 
			
			\State{$u\leftarrow i$ \textit{ \# Go back to the initial node}}
			\State{\textit{\# Add the newly created branch to the forest}}
			\While{\textbf{not} $\text{InForest}[u]$}
			\State{$\text{InForest}[u] \leftarrow \texttt{true}$}
			\State{$u\leftarrow\text{Next}[u]$}
			\EndWhile
			\EndFor
			\State{\texttt{return} Next, $\mathcal{P}$}
		\end{algorithmic}
	\end{algorithm}

	\begin{algorithm}
		\caption{$\left(Q, \{\text{Next}_q\}_{q\in Q}\right)$ = \texttt{CoupledRandomForest}($\mathcal{G},q_{\text{min}},q_{\text{max}}$)}
		\label{alg:coupled_RSF}
		\begin{algorithmic}[1]
			\Inputs{$\mathcal{G}=(\cl X,w)$\\
				$q_{\text{max}}\geq q_{\text{min}}>0$}
			\nt{\Initialize{
				$Q\leftarrow\emptyset$ \textit{\# Initialize an empty list}\\
				$\left(\text{Next}, \mathcal{P}\right)$ = \texttt{RandomForest}($\mathcal{G},q_{\text{max}}$) \textit{\# sample a random forest at $q_{\text{max}}$}
				\\
				Save: add $q_{\text{max}}$ to $Q$, $\text{Next}_{q_{\text{max}}}\leftarrow\text{Next}$
			} 
			\While{$q>q_{\text{min}}$}
				\State{$(q, i) \leftarrow$ \texttt{PopLargest}$(\mathcal{P})$} \textit{\# find the couple $(q,i)$  in $\set{P}$ with maximal $q$ and remove it from $\mathcal{P}$}
				\State{$\text{Next[i]} \leftarrow \texttt{RandomSuccessor}(i,\mathcal{G})$} \textit{\# Reactivate: draw next of $i$, $i$ is no longer a root}
				\State{$r \leftarrow \texttt{LookAhead}(i,\text{Next})$} \textit{\# $r$ is the first node down the Next vector that is a root or equals $i$}
				\If{$r \neq i$} 
					\State{\textbf{pass}} \textit{~~~\# there is nothing to do ! Go directly to line~\ref*{line:end} }
				\Else \textit{~~~\# if $r==i$: a cycle is created}
					\State{$\mathcal{A} \leftarrow \texttt{GetCycle}(i,\text{Next})$} \textit{\# $\mathcal{A}$ is the set of nodes in the cycle: they all need to be reactivated}
					\While{$\mathcal{A}$ is not empty}
						\State{$u \leftarrow \texttt{Pop}(\mathcal{A})$} \textit{\# pop an arbitrary node from $\mathcal{A}$ ($u$ is removed from $\mathcal{A}$)}
						\State{$\mathcal{U} \leftarrow Unif([0,1])$}
						\If{$\mathcal{U} \leq \frac{q}{q+ w[u]}$} \textit{\# If true, $u$ becomes a root}
							\State{$\text{Next}[u] \leftarrow 0$ \textit{ \# Set next of $u$ to null}}
							\State{Add $\left(q'=\frac{\mathcal{U} w[u]}{1-\mathcal{U}},u\right)$ to the priority queue $\mathcal{P}$}
						\Else \textit{ \# If false, continue the random walk}
							\State{$\text{Next[u]} \leftarrow \texttt{RandomSuccessor}(u,\mathcal{G})$}
							\State{$r' \leftarrow \texttt{LookAhead}(u, \text{Next}, \mathcal{A})$} \textit{\# $r'$ is the  first node down the Next vector that is a root, is in $\mathcal{A}$, or equals $u$}
							\If{$r' == u$} \textit{~~~\# a new cycle is created}
								\State{$\mathcal{B} \leftarrow \texttt{GetCycle}(u,\text{Next})$} \textit{\# $\mathcal{B}$ is the set of all nodes in this new cycle}
								\State{$\mathcal{A} \leftarrow \mathcal{A}\cup\mathcal{B}$}
							\EndIf
						\EndIf
					\EndWhile
				\EndIf
				\State{Save: Add $q$ to $Q$, $\text{Next}_q \leftarrow \text{Next}$\label{line:end}} \textit{    \# Save the state of the forest at $q$}
			\EndWhile
			\State{\texttt{return} $Q, \{\text{Next}_q\}_{q\in Q}$}}
		\end{algorithmic}
	\end{algorithm}

	\clearpage
	
	\section{Canonical representations and Markov's bounds
		for the truncated moment problem on a bounded interval}
	\label{daniele}
	
	We collect here some basics of truncated moments problems,
	which provide our needed Markov's bounds
	and that we essentially learned from \cite{KN}.
	Although they rely on well-know quadrature formulas,
	one important difficulty we faced in making this work
	was to recover these bounds,
	that we felt to have been 
	at least partially forgotten
	in our community.

	\subsection{Admissible moment sequences and quadrature formulas}
	
	Let $[a, b] \subset \ds R$ and $m_0$, $m_1$, \dots, $m_l$ be $n = l + 1$ real numbers.
	We say that $m_0$, \dots, $m_l$ is an {\bf admissible} moment sequence
	if there exists a non-negative measure $\mu$ on $[a, b]$ such that
	$$
	m_k = \int_a^b x^k \mu(dx), 
	\qquad 0 \leq k \leq l.
	$$
	In this case we say that $\mu$ is a {\it representation\/}
	of $m_0$, \dots, $m_l$.
	We denote by $\cl D^l_{a, b} \subset \ds R^n$
	the set of all such admissible moment sequences
	and we set $\cl D_{a, b} = \cup_{l \geq 0} \cl D^l_{a, b}$
	
	For $m \in \cl D^l_{a, b}$ with representation $\mu$,
	the linear form on the $n$-dimensional space $\ds R_l[X]$
	of the polynomials with real coefficients
	and degree less than of equal to $l$
	$$
	L_m: P = \sum_{k=0}^l p_k X^k  \in \ds R_l[X]
	\mapsto \int_a^b P(x) \mu(dx) = \sum_{k=0}^l p_k m_k
	$$
	is completely determined by $m$ only.
	In particular,
	if an atomic measure $\nu = \sum_{i \in I} w_i \delta_{x_i}$
	with atoms $x_i$, $i \in I$, in $[a, b]$
	provides another representation of $m$,
	i.e,
	\begin{equation}\label{adel}
		\sum_{i \in I} w_i x_i^k = m_k,
		\qquad 0 \leq k \leq l,
	\end{equation}
	then we have the quadrature formula
	\begin{equation}\label{guillaume}
		L_m(P) =  \int_a^b P(x)\mu(dx) = \sum_{i \in I} w_i P(x_i),
		\qquad P \in \ds R_l[X].
	\end{equation}
	Equations~\eqref{adel} and~\eqref{guillaume} 
	are actually equivalent.

	\subsection{The first orthogonal polynomials}
	
	Let $\mu$ be a measure on $[a, b]$ with $n = l + 1$
	first moments $m_0$, $m_1$, \dots, $m_l$.
	This moment sequence $m$ does in general {\it not\/}
	determine completely the bilinear form on $\ds R_l[X]$
	$$
	\langle P, Q\rangle_\mu
	= \int_a^b P(x)Q(x)\mu(dx),
	\qquad P, Q \in \ds R_l[X],
	$$
	since such quantities usually depend on moments of larger order
	$l + 1$, \dots, $2l$.
	As a consequence it does not determine the unitary and normalized orthogonal polynomials
	$U_k$ and $V_k = U_k / \langle U_k, U_k\rangle_\mu^{1 / 2}$, $k \geq 0$,
	associated with a regular enough measure $\mu$.
	The moments $m_0$, \dots, $m_l$ do however determine the first of them,
	which can be computed inductively with the formulas
	$$
	U_0 = 1,
	\qquad V_0 = {1 \over \langle U_0, U_0\rangle^{1 / 2}_\mu} U_0;
	$$
	$$
	\beta_0 = {\langle XU_0, U_0\rangle_\mu \over \langle U_0, U_0\rangle_\mu},
	\qquad U_1 = (X - \beta_0)U_0,
	\qquad V_1 = {1 \over \langle U_1, U_1\rangle^{1 / 2}_\mu} U_1;
	$$
	and, for $k \geq 1$,
	$$
	\beta_k = {\langle XU_k, U_k\rangle_\mu \over \langle U_k, U_k\rangle_\mu},
	\qquad \gamma_k = {\langle XU_k, U_{k - 1}\rangle_\mu \over \langle U_{k - 1}, U_{k -1}\rangle_\mu},
	$$
	$$
	U_{k + 1} = (X - \beta_k)U_k - \gamma_k U_{k - 1},
	\qquad V_{k + 1} = {1 \over \langle U_{k + 1}, U_{k + 1}\rangle^{1 / 2}_\mu} U_{k + 1};
	$$
	provided that none of these denominators vanishes.
	If some of them vanishes
	we say that $m$ is a {\bf singular\/} moment sequence,
	otherwise we say that $m$ is {\bf regular}
	---so that for any admissible moment sequence $m_0$,~\dots, $m_l$
	there is $0 \leq k \leq l$ such that $m_0$, \dots, $m_j$
	is regular for all $j < k$ and singular for all $j \geq k$.
	
	These formulas show that
	all moments up to order $2k + 1 = 2(k + 1) - 1$
	have to be known to compute $U_{k + 1}$,
	and all moments up to order $2k + 2 = 2(k + 1)$
	have to be known to compute $V_{k + 1}$.
	As a consequence :
	\begin{itemize}
		\item if $l = 2r$ is even,
		then $m$ determines the first $r + 1$ unitary and normalized orthogonal polynomials
		$U_0$, \dots, $U_r$ and $V_0$, \dots, $V_r$ ;
		\item if $l = 2r + 1$ is odd,
		then $m$ determines the first $r + 2$ unitary orthogonal polynomials
		$U_0$, \dots, $U_{r + 1}$,
		but only the first $r + 1$ normalized orthogonal polynomials
		$V_0$, \dots, $V_r$.
	\end{itemize}
	We recall that orthogonal polynomials of degree $k \geq 0$
	have $k$ distinct roots in $]a, b[$.

	\subsection{Principal representations and quadratures}
	
	When looking for a quadrature formula on $[a, b]$
	or an atomic representation $\nu$ of $n = l + 1$ moments
	of a measure $\mu$,
	we will distinguish four different cases depending whether 
	neither $a$ nor $b$ are atoms of $\nu$,
	both $a$ and $b$ are atoms of $\nu$,
	only $a$ is an atom of $\nu$
	or only $b$ is an atom of $\nu$.
	Solving Equation~\eqref{adel} with $k$ atoms outside $\{a, b\}$
	amounts then to solve a problem with $n$ equations
	and $2k$, $2k + 2$, $2k + 1$ or $2k + 1$ degrees of freedom,
	respectively.
	For an atomic distribution $\nu$ on $[a, b]$
	with $j$ atoms in $\{a, b\}$
	and $k$ atoms outside $\{a, b\}$ 
	we then define the {\bf index\/} of $\nu$ as
	$$
	{\rm ind}(\nu) = j + 2k.
	$$
	
	We will treat separately the ``odd case'' $l = 2r + 1$
	and the ``even case'' $l = 2r$.
	For a regular moment sequence,
	in both cases we will get two uniquely defined atomic representations
	of index $n$ and known as the {\bf upper\/}
	and {\bf lower principal representations}
	depending whether $b$ is an atom or not.
	In all cases this is nothing more than rebuilding
	some well-known quadrature formulas.

	\subsubsection {Gauss and Lobatto quadratures for the odd case}
	
	\noindent\textbf{Lower principal representation.}
	For $l = 2r + 1$
	and a regular sequence $m$
	of $n = l + 1 = 2(r + 1)$ moments $m_0$, \dots, $m_l$,
	let us first look for
	an atomic representation of $m$
	with index $n$ and without mass in $\{a, b\}$,
	i.e., for $x_0$, \dots, $x_r$ in $]a, b[$
	and $w_0$, \dots, $w_r > 0$,
	such that 
	\begin{equation}\label{kai}
		\sum_{i = 0}^{r} w_i x_i^k = m_k,
		\qquad 0 \leq k \leq l.
	\end{equation}
	Let $Q = \prod_{i = 0}^r (X - x_i)$. 
	Since $Q$ is of degree $r + 1$,
	Equation~\eqref{guillaume} gives then
	$$
	\langle Q, X^j\rangle_\mu = 0,
	\qquad 0 \leq j \leq l - (r + 1) = r,
	$$
	which identifies $Q$ as the unitary orthogonal polynomial $U_{r + 1}$,
	and the $x_i$, $0 \leq i \leq r$, as its roots.
	We can then solve the Vandermonde system
	given by the first $r + 1$ equations in~\eqref{kai}
	to identify the $w_i$, $0 \leq i \leq r$.
	
	To check that these $x_i$ and $w_i$, $0 \leq i \leq r$,
	form a solution indeed,
	we first note that the quadrature formula~\eqref{guillaume}
	holds by construction for 
	$$
	P = 1,\, X,\, \dots,\, X^r,\, Q,\, XQ,\, \dots,\, X^r Q,
	$$
	that is on a whole basis of $\ds R_l[X]$.
	It only remains to check
	that our weights $w_i$, $0 \leq i \leq r$,
	are positive.
	To this end we consider for any $0 \leq j \leq r$,
	the Legendre polynomial $L_j$ such that $L_j(x_j) = 1$,
	and $L_j(x_i) = L'_j(x_i) = 0$ for $i \neq j$.
	This is a positive polynomial of degree $2r = l - 1$
	and our quadrature formula gives
	$$
	w_j = \int_a^b L_j(x) \mu(dx) \geq 0
	$$
	for any representation $\mu$ of $m$.
	And had some $w_j$ been equal to zero,
	we would have get an atomic representation of $m$
	with index less than $n$,
	so that $m$ would have been singular.
	
	A regular moment sequence of $n = 2(r + 1)$ moments
	admits then a unique atomic representation
	with $r + 1$ atoms in $]a, b[$, which gives rise
	to the standard Gaussian quadrature.
	We call it its {\bf lower principal representation}.
	
	\noindent\textbf{Upper principal representation.}
	For $l = 2r + 1$
	and a regular sequence $m$
	of $n = l + 1 = 2r + 2$ moments $m_0$, \dots, $m_l$,
	we now look for
	an atomic representation of $m$
	with index $n$ and mass both in $a$ and $b$,
	i.e., for $r$ atoms $x_i$, $i < r$, in $]a, b[$
	as well as $w_a$, $w_b$ and $r$ weights $w_i$, $i < r$,
	in $\ds R_+ \setminus \{0\}$
	such that 
	\begin{equation}\label{glenn}
		w_a a^k + w_b b^k + \sum_{i < r} w_i x_i^k = m_k,
		\qquad 0 \leq k \leq l.
	\end{equation}
	For $\mu$ a representation of $m$,
	let $\mu^{a, b}$ be the measure on $[a, b]$
	defined by 
	\begin{equation}\label{domi}
		\mu^{a, b}(dx) = (b - x)(x - a) \mu(dx),
		\qquad x \in [a, b], 
	\end{equation}
	and let 
	$m_0^{a, b}$, \dots, $m^{a, b}_{l^{a, b}}$,
	with $l^{a,b} = l - 2 = 2r - 1 = 2(r - 1) + 1$,
	be its first $1 + l^{a, b}$ moments,
	which depend on $m$ only.
	The polynomial $Q = \prod_{i < r} (X - x_i)$ has degree $r$
	and Equation~\eqref{guillaume} gives
	$$
	\langle Q, X^j\rangle_{\mu^{a, b}}
	= \langle (b - X)(X - a)Q, X^j\rangle_\mu
	= 0,
	\qquad 0 \leq j \leq l - (2 + r) = r - 1,
	$$
	which identifies $Q$ with the unitary polynomial $U_r^{a, b}$
	associated with $m_0^{a, b}$, \dots, $m_{l^{a,b}}^{a, b}$.
	As previously this identifies the $x_i$, $i < r$,
	as the roots of $U^{a,b}_r$
	and the first $r + 2$ equations in~\eqref{glenn}
	identify $w_a$, $w_b$ and $w_i$, $i < r$.
	The quadrature formula~\eqref{guillaume}
	holds for 
	$$
	P = 1,\, X,\, \dots,\, X^{r + 1},\, (b - X)(X -a)Q,\, X(b -X)(X - a)Q,\, \dots,\, X^{r - 1}(b - X)(X - a)Q,
	$$
	again a whole basis of $\ds R_l[X]$.
	And we can conclude in the same way by considering the Legendre polynomials $L_j$, $j < r$,
	such that $L_j(x_j) = 1$ and $L_j(a) = L_j(b) = L_j(x_i) = L_j'(x_i) = 0$ for $i \neq j$
	as well as $L_a$ and $L_b$ such that $L_a(a) = L_b(b) = 1$
	and $L_a(b) = L_b(a) = L_a(x_i) = L_b(x_i) = L_a'(x_i) = L_b'(x_i) = 0$, $i < r$,
	which are non-negative on $[a, b]$ and have degree $2(r - 1) + 2 = 2r = l - 1$, 
	$2r + 1 = l$  and $2r + 1 = l$, respectively.
	
	A regular moment sequence of $n = 2(r + 1)$ moments
	admits then a unique atomic representation
	with $r$ atoms in $]a, b[$, one in $a$ and one in $b$,
	which gives rise to the Gauss-Lobatto quadrature.
	We call it its {\bf upper principal representation}.

	\subsubsection {Gauss-Radau quadratures for the even case}
	
	For $l = 2r$
	and a regular sequence $m$
	of $n = l + 1 = 2r + 1$ moments $m_0$, \dots, $m_l$,
	we can look for an atomic representation of $m$
	with index $n$ and some mass either in $a$ or $b$.
	As previously, this leads to the Gauss-Radau quadratures
	associated with the {\bf lower} or {\bf upper} {\bf principal representations},
	respectively.
	
	In the first case we look for
	$x_i$, $i < r$, in $]a, b[$
	as well as $w_a$ and  $w_i$, $i < r$, in $\ds R_+ \setminus \{0\}$
	such that 
	\begin{equation}\label{nicolas}
		w_a a^k + \sum_{i < r} w_i x_i^k = m_k,
		\qquad 0 \leq k \leq l.
	\end{equation}
	The $x_i$, $i < r$, are the roots
	of the unitary orthogonal polynomial $U^a_r$ with degree $r$
	and associated with $\mu^a$ defined by 
	\begin{equation}\label{pierre}
		\mu^a(dx) = (x - a)\mu(dx),
		\qquad x \in [a, b]. 
	\end{equation}
	The last known moment or $\mu^a$
	is of order $l^a = l - 1 = 2r - 1 =2(r - 1) + 1$. Then, 
	$w_a$ as well as $w_i$, $i < r$,
	are identified by inverting the Vandermonde system
	given by the first $r + 1$ equations in~\eqref{nicolas}.
	
	In the second case we solve in the same way
	$$
	w_b b^k + \sum_{i < r} w_i x_i^k = m_k,
	\qquad 0 \leq k \leq l.
	$$
	The $x_i$, $i < r$, are the roots
	of the unitary orthogonal polynomial $U_r^b$ with degree $r$
	and associated with $\mu^b$ defined by
	\begin{equation}\label{bruno}
		\mu^b(dx) = (b - x)\mu(dx),
		\qquad x \in [a, b].
	\end{equation}

	\subsection{The truncated moment problem}
	\label{yasmina}
	
	The truncated moment problem
	is that of identifying $\cl D_{a, b}$.
	After observing that $\cl D^0_{a, b} = \ds R_+$
	and that $m_0 = 0$ constitutes
	the only singular sequence in $\cl D^0_{a, b}$,
	the problem can be solved with the following
	\begin{lmm}\label{frank}
		If $m = (m_0, \dots, m_l)$ is a regular moment sequence
		with lower and upper principal representations $\nu^-$ and $\nu^+$,
		then their respective next moments $m^-_{l + 1}$ and $m^+_{l + 1}$
		satisfy
		$$
		m^-_{l + 1} = \int_a^b x^{l + 1} \nu^-(dx) 
		< m^+_{l + 1} = \int_a^b x^{l + 1} \nu^+(dx)
		$$
		and it holds
		$$
		m^-_{l + 1}
		\leq \int_a^b x^{l + 1} \mu(dx) 
		\leq m^+_{l + 1}
		$$
		for any representation $\mu$ of $m$.
		
		In addition, if $\mu$ has a moment $m_{l + 1}$ of order $l + 1$
		that coincides with $m^-_{l + 1}$ or $m^+_{l + 1}$,
		then $\mu$ coincides with $\nu^-$ or $\nu^+$,
		respectively.
		In this case $m_0$, \dots, $m_{l + 1}$
		is a singular moment sequence.
		Conversely, any sequence $m_0$, \dots, $m_{l + 1}$
		is a regular moment sequence in $\cl D^{l + 1}_{a, b}$
		provided that $m^-_{l + 1} < m_{l + 1} < m^+_{l + 1}$.
	\end{lmm}
	
	\noindent{\bf Proof:}
	Consider in the odd case $l = 2r + 1$,
	with $x_0$, \dots, $x_r$ the atoms of $\nu^-$,
	the positive polynomial 
	$$
	P = \prod_{i = 0}^r (X - x_i)^2
	= X^{2r + 2} + Q.
	$$
	Writing $m_{l + 1}$ for the moment of order $l + 1$ 
	of a representation $\mu$ of $m$,
	it holds 
	$$
	0 \leq \int_a^b P(x)\mu(dx) 
	= m_{l + 1} + \int_a^b Q(x)\mu(dx)
	= m_{l + 1} + \int_a^b Q(x) \nu^-(dx),
	$$
	since $Q \in \ds R_l[X]$.
	The last term in the previous equation is $- m^-_{l + 1}$,
	since it also holds
	$$
	0 = \int_a^b P(x)\nu^-(dx) 
	= m^-_{l + 1} + \int_a^b Q(x) \nu^-(dx).
	$$
	This proves $m^-_{l + 1} \leq m_{l + 1}$.
	And in case of equality,
	it holds
	$$
	\int_a^b P(x)\mu(dx) = 0,
	$$
	which implies that the support of $\mu$
	is contained in that of $\nu^-$.
	Writing $\mu = \sum_{i = 0}^r w_i \delta_{x_i}$
	the weights $w_i$ are then identified
	by the Vandermonde equations
	$$
	\sum_{i = 0}^r w_i x_i^k = m_k,
	\qquad 0 \leq k \leq r,
	$$
	which implies that $\mu$ and $\nu^-$ coincide, indeed.
	
	The upper bound $m_{l + 1} \leq m^+_{l + 1}$
	is proved similarly by considering 
	the positive polynomial on $[a, b]$
	$$
	P = (b - X)(X - a)\prod_{i < r}(X - x_i)^2
	= - X^{2r + 2} + Q,
	$$
	with $a$, $b$ and $x_i$, $i < r$, the atoms of $\nu^+$.
	And we can deal with the case of equality in the same way.
	
	Since lower and upper principal representations
	are distinct representations of a regular moment sequence,
	this also proves that $m^-_{l + 1} < m^+_{l + 1}$.
	Finally, if $m_{l + 1} = \alpha m_{l + 1}^- + (1 - \alpha)m^+_{l + 1}$
	for some $\alpha \in\; ]0, 1[$, then $(m_0, \dots, m_{l + 1})$ 
	admits $\alpha \nu^- + (1 - \alpha)\nu^+$ as a representation:
	it forms a regular moment sequence.
	
	We can deal with the even case $l = 2r$ in the same way
	by considering the polynomial
	$$
	P = (X - a)\prod_{i < r} (X - x_i)^2
	$$
	and
	$$
	P = (b - X)\prod_{i < r}(X - x_i)^2,
	$$
	with $x_i$, $i < r$, the atoms
	of the lower and upper principal representation, respectively.
	\qed
	
	We note, as the consequence of this lemma,
	that a sequence of $n = l + 1$ moments 
	$m_0$, \dots, $m_l$ forming a singular sequence
	has a unique representation $\nu$,
	which is atomic and satisfies ${\rm ind}(\nu) < n$.
	
	\subsection{Canonical representations}
	
	Consider $n = l + 1$ moments $m_0$, \dots, $m_l$
	forming a regular sequence $m$.
	For each $\xi \in \;]a, b[$
	there is a unique atomic measure $\nu^\xi$
	with an atom at $\xi$
	and of index strictly less than $n + 2$
	---actually $n$ or $n + 1$---
	known as the {\bf canonical representation\/} of $m$
	associated with $\xi$.
	
	To see that such a representation exists,
	we can consider
	$$
	t = \sup\bigl\{s \geq 0 : m - (s, s\xi, s\xi^2, \dots, s\xi^l) \in \cl D_{a,b}^l\bigr\} \leq m_0.
	$$
	The previous lemma implies
	that $t > 0$ and $m - (t, \dots, t\xi^l)$ is singular,
	i.e., has an atomic representation $\nu$ with index strictly less that $n$.
	So that $\nu + t\delta_\xi$ is an atomic representation of $m$ 
	with index $n + 1$ at most and with an atom in $\xi$.
	We show how to build it in this section.

	\subsubsection{The odd case}
	
	In the odd case $l = 2r + 1$, i.e.,
	$$
	n + 1 = l + 2 = 2r + 1 + 2,
	$$
	we are looking either for $x_i \in [a, b] \setminus \{\xi\}$
	and $w_i \geq 0$, $i < r$,
	as well as $w_a \geq 0$ and $t > 0$ such that
	\begin{equation}\label{veronique}
		t \xi^k + w_a a^k + \sum_{i < r} w_i x_i^k = m_k,
		\qquad k < 2r + 2,
	\end{equation}
	or for $x_i \in [a, b] \setminus \{\xi\}$
	and $w_i \geq 0$, $i < r$,
	as well as $w_b \geq 0$ and $t > 0$ such that
	\begin{equation}\label{blandine}
		t \xi^k + w_b b^k + \sum_{i < r} w_i x_i^k = m_k,
		\qquad k < 2r + 2.
	\end{equation}
	In the first case,
	consider some representation $\mu$ of $m$
	and the measure $\mu^a$ defined by~\eqref{pierre}.
	The first $l + 1$ moments $m_0$, \dots, $m_l$ of $\mu$
	define the first moments $m^a_0$, \dots, $m^a_{l^a}$ of $\mu^a$,
	with $l^a = l - 1 = 2r$,
	hence the first $r + 1$ associated normalized orthogonal polynomials
	$V^a_0$, \dots, $V^a_r$.
	With $P = \prod_{i < r}(X - x_i)$, Equation~\eqref{guillaume} gives then
	$$
	\langle V_k^a, P\rangle_{\mu^a}
	= \langle (X - a)V_k^a, P\rangle_\mu
	= t (\xi - a)P(\xi)V_k^a(\xi),
	\qquad 0 \leq k \leq r.
	$$
	Since $P \in \ds R_r[X]$,
	this shows that $P$ and 
	$$
	Q^a = \sum_{i = 0}^r V_i^a(\xi)V_i^a
	$$
	are proportional,
	namely
	$$
	P = t(\xi - a)P(\xi)Q^a,
	$$
	and gives
	$$
	t = {1 \over (\xi - a)Q^a(\xi)},
	$$
	since $P(\xi)$ and $Q^a(\xi)$ are not zero.
	Provided that~\eqref{veronique} does have a solution
	in the prescribed domain,
	this solution is unique,
	since the singular moment sequence $m - (t, \dots, t\xi^l)$
	admits one representation only.
	Our $x_i$, $i < r$, are in practice the roots of $Q^a$;
	$w_a$ and $w_i$, $i < r$, can then be found
	by solving the Vandermonde system
	given by the first equations in~\eqref{veronique}. 
	
	In the second case
	we get in the same way
	$$
	t = {1 \over (b - \xi)Q^b(\xi)},
	$$
	with
	$$
	Q^b = \sum_{i = 0}^r V_k^b(\xi)V_k^b,
	$$
	where $V_0^b$, \dots, $V_r^b$ are the first $r + 1$
	normalized orthogonal polynomials associated
	with the first $l^b + 1  = 2r + 1$ moments
	$m^b_0$, \dots, $m^b_{l^b}$
	of $\mu^b$ defined by~\eqref{bruno}.
	Again,
	our $x_i$, $i < r$, are the roots of $Q^b$,
	and the weights
	$w_i$, $i < r$, as well as $w_b$
	can be found by  solving the Vandermonde system
	given by the first equations in~\eqref{blandine}.
	
	This proves that $\nu^\xi$ must be either
	the unique solution of~\eqref{veronique},
	if it exists,
	or that of~\eqref{blandine},
	if it exists
	---we already know that at least one of them exists.
	This does not completely settle the uniqueness issue,
	except when $(\xi -a)Q^a(\xi) = (b - \xi)Q^b(\xi)$,
	i.e., when~\eqref{veronique} and~\eqref{blandine}
	implies a same value for $t$,
	so that $m - (t, \dots, t\xi^l)$
	is a singular moment sequence
	with only one representation.
	We will see later 
	that when $(\xi -a)Q^a(\xi) \neq (b - \xi)Q^b(\xi)$,
	only one of the two equations~\eqref{veronique} and~\eqref{blandine}
	can lead to a solution inside the desired domain.

	\subsubsection{The even case}
	
	For $l = 2r$ we have
	$$
	n + 1 = l + 2 = 2r + 2
	$$
	and we are looking either for 
	$x_i \in [a, b] \setminus \{\xi\}$
	and $w_i \geq 0$, $i < r$, as well as $t > 0$
	such that 
	\begin{equation}\label{lina}
		t \xi^k + \sum_{i < r} w_i x_i^k = m_k,
		\qquad k < 2r + 1,
	\end{equation}
	or for
	$x_i \in [a, b] \setminus \{\xi\}$
	and $w_i \geq 0$, $i < r - 1$, as well as
	$w_a \geq 0$, $w_b \geq 0$ and $t > 0$,
	such that 
	\begin{equation}\label{marc-alain}
		t \xi^k + w_a a^k + w_b b^k + \sum_{i < r - 1} w_i x_i^k = m_k,
		\qquad k < 2r + 1 = 2(r - 1) + 2 + 1,
	\end{equation}
	Similarly to the odd case,
	given a representation $\mu$ of $m$,
	its first $l + 1$ moments $m_0$, \dots, $m_l$ define 
	the first $r + 1$ normalized orthogonal polynomials
	$V_0$, \dots, $V_r$ associated with $\mu$,
	as well as the first $l^{a, b} + 1 = l - 2 + 1 = 2(r - 1) + 1$ moments
	$m^{a, b}_0$, \dots, $m^{a, b}_{l^{a, b}}$ 
	and the first $r$ normalized orthogonal polynomials
	$V^{a, b}_0$,~\dots, $V^{a, b}_{r - 1}$ associated with $\mu^{a, b}$
	defined by~\eqref{domi}.
	We can then define
	$$
	Q = \sum_{k = 0}^r V_k(\xi)V_k
	$$
	and 
	$$
	Q^{a, b} = \sum_{k < r} V^{a, b}_k(\xi)V^{a, b}_k
	$$
	to get that~\eqref{lina} leads to 
	$$
	t = {1 \over Q(\xi)}
	$$
	and identifies the $x_i$, $i < r$,
	as the roots of $Q$,
	then the weights $w_i$, $i < r$, as the solution of the Vandermonde system
	given by the first equations in~\eqref{lina},
	while~\eqref{marc-alain} leads to 
	$$
	t = {1 \over (b - \xi)(\xi - a)Q^{a, b}(\xi)}
	$$
	and identifies $x_i$, $i < r - 1$,
	as the roots of $Q^{a, b}$,
	then the weights $w_i$, $i < r$, as well as $w_a$ and $w_b$
	as the solution of the Vandermonde system
	given by the first equations in~\eqref{marc-alain}.
	
	Again, we know that one at least of~\eqref{lina}
	and~\eqref{marc-alain} gives a solution in the prescribed domain,
	which shows that the canonical representation $\nu^\xi$
	is uniquely defined in the special case where
	$$
	Q(\xi)= (b - \xi)(\xi - a)Q^{a, b}(\xi).
	$$
	We will see in the next section that
	equations~\eqref{lina} and~\eqref{marc-alain}
	cannot both lead to a solution in the desired domain 
	when
	$$
	Q(\xi) \neq (b - \xi)(\xi - a)Q^{a, b}(\xi),
	$$
	so that the canonical representation $\nu^\xi$
	is always uniquely defined, indeed.

	\subsection{Markov's bounds}\label{louise}
	Markov's bounds are tight lower an upper bounds
	on $\mu(]\xi, b])$ and $\mu([\xi, b])$
	for a given $\xi \in \;]a, b[$ and any measure $\mu$
	the first $n = l + 1$ moments of which form a given regular moment sequence.
	These bounds are derived from an associated canonical representation
	$$
	\nu^\xi = t \delta_\xi + \nu'
	$$
	with $t > 0$ and ${\rm ind}(\nu') \leq n + 1 - 2 = l$.
	They read
	\begin{equation}\label{sophie}
		\nu^\xi(]\xi, b]) 
		\leq \mu(]\xi, b])
		\leq \mu([\xi, b])
		\leq \nu^\xi([\xi, b])
	\end{equation}
	and 
	\begin{equation}\label{etienne}
		\mu(\{\xi\}) = \mu([\xi, b]) - \mu(]\xi, b]) \leq t.
	\end{equation}
	
	When ${\rm ind}(\nu') = 2k + 1$  is odd
	and $c$ is the only one atom of $\nu'$ in $\{a, b\}$,
	call $x_i$, $i < k$, the interior atoms of $\nu'$
	and consider the Legendre polynomials $L_-$ and $L_+$
	defined by 
	$$
	L_-(\xi) = 0,
	\qquad L_+(\xi) = 1,
	$$
	$$
	L_-(c) = L_+(c) = \ds 1_{\{c > \xi\}},
	\qquad L_-(x_i) = L_+(x_i) = \ds 1_{\{x_i > \xi\}},
	\qquad i < k,
	$$
	and
	$$
	L'_-(x_i) = L'_+(x_i) = 0,
	\qquad i < k.
	$$
	$L_-$ and $L_+$ satisfy
	$$
	L_- \leq \ds 1_{]\xi, b]} \leq \ds 1_{[\xi, b]} \leq L_+
	$$
	and have degree $2k + 1 \leq l$.
	It follows that
	$$
	\begin{aligned}
		\nu'(]\xi,b]) = \nu^\xi(]\xi, b]) = \int_a^b L_-(x)\nu^\xi(dx) &= \int_a^b L_-(x)\mu(dx)\\
		&\leq \int_a^b \ds 1_{]\xi, b]}(x) \mu(dx) = \mu(]\xi, b])
	\end{aligned}
	$$
	and
	$$
	\begin{aligned}
		\mu([\xi, b]) &= \int_a^b \ds 1_{[\xi, b]}(x) \mu(dx) \\
		&\leq \int_a^b L_+(x) \mu(dx) = \int_a^b L_+(x) \nu^\xi(dx) = \nu^\xi([\xi, b]) = t + \nu'(]\xi, b]),
	\end{aligned}
	$$
	which proves~\eqref{sophie} and~\eqref{etienne}.
	
	When ${\rm ind}(\nu') = 2k$ and $x_i \in \;]a, b[$, $i < k$,
	are all the atoms of $\nu'$,
	consider the Legendre polynomials $L_-$ and $L_+$ defined by
	$$
	L_-(\xi) = 0,
	\qquad L_+(\xi) = 1,
	$$
	$$
	L_-(x_i) = L_+(x_i) = \ds 1_{\{x_i > \xi\}},
	\qquad i < k,
	$$
	and
	$$
	L'_-(x_i) = L'_+(x_i) = 0,
	\qquad i < k.
	$$
	Since these have degree $2k \leq l$,~\eqref{sophie} and~\eqref{etienne}
	follow in the same way.
	
	Finally, when ${\rm ind}(\nu') = 2k$ and $x_i \in \;]a, b[$, $i < k - 1$,
	as well as $a$ and $b$ are the atoms of $\nu'$,
	the same result is obtained by considering the Legendre polynomials
	of degree $2k \leq l$ defined by
	$$
	L_-(\xi) = 0,
	\quad L_+(\xi) = 1,
	$$
	$$
	L_-(a) = L_+(a) = 0,
	\quad L_-(b) = L_+(b) = 1,
	\quad L_-(x_i) = L_+(x_i) = \ds 1_{\{x_i > \xi\}},
	\qquad i < k - 1,
	$$
	and
	$$
	L'_-(x_i) = L'_+(x_i) = 0,
	\qquad i < k - 1.
	$$
	
	Now, Equation~\eqref{etienne} shows that there cannot be
	two canonical representations associated with a same $\xi \in \;]a, b[$ 
	and two different weights in $\xi$: 
	we would get a contradiction by taking that with the larger weight
	in place of $\mu$ and that with the smaller weight $t$ in place of $\nu^\xi$.
	This concludes the proof of uniqueness for canonical representations.
	
	Krein and Nudel'man even show in \cite{KN} 
	that the weight in any $\xi \in \;]a, b[$ of the canonical representation $\nu^\xi$
	is always, with the notation of the previous section,
	$$
	t = {1 \over \max\bigl((\xi -a)Q^a(\xi),\, (b - \xi)Q^b(\xi)\bigr)}
	$$
	in the odd case or
	$$
	t = {1 \over \max\bigl(Q(\xi),\, (b - \xi)(\xi - a)Q^{a, b}(\xi)\bigr)}
	$$
	in the even case.
	Here, we simply conclude:
	\begin{thm*}{\bf (Markov bounds)}
		For any $\xi \in ]a, b[$, each regular moment sequence 
		$m = (m_0, \dots, m_l)$ admits a unique canonical representation $\nu^\xi$
		with an atom in $\xi$. 
		Then, for any representation $\mu$ of $m$, it holds
		$$
		\nu^\xi(]\xi, b]) 
		\leq \mu(]\xi, b])
		\leq \mu([\xi, b])
		\leq \nu^\xi([\xi, b]).
		$$
	\end{thm*}
	
\end{appendix}

\end{document}